\documentclass[journal,draftcls,onecolumn,letterpaper,twoside,11pt]{IEEEtran}

\usepackage[margin=1in]{geometry}
\usepackage{amssymb, mathrsfs,verbatim} 
\usepackage[cmex10]{amsmath}
\usepackage{cite}
\RequirePackage[colorlinks,linkcolor=black,citecolor=black,urlcolor=blue]{hyperref}
\usepackage{enumerate}
\usepackage[dvips]{graphicx}
\usepackage{verbatim}
\usepackage{bbm}
\usepackage{subfig}
\usepackage{epsfig}
\usepackage{epstopdf}
\usepackage{graphicx} 

\usepackage{float,latexsym,amsfonts,amstext,times}

\usepackage{chngcntr}
\usepackage{apptools}

\usepackage[lined,linesnumbered,algoruled,noend]{algorithm2e}
\usepackage{cases}

\usepackage{color}

\usepackage{multicol}

\newtheorem{theorem}{Theorem}[section]
\newtheorem{definition}{Definition}[section]
\newtheorem{proposition}{Proposition}[section]
\newtheorem{corollary}{Corollary}[section]
\newtheorem{lemma}{Lemma}[section]
\newtheorem{remark}{Remark}[section]

\newcommand{\cS}{\mathcal{S}}
\newcommand{\cJ}{\mathcal{J}}
\newcommand{\cD}{\mathcal{D}}
\newcommand{\cC}{\mathcal{C}}
\newcommand{\cF}{\mathcal{F}}

\newcommand{\cL}{\mathcal{L}}

\newcommand{\cV}{\mathcal{V}}

\newcommand{\Pro}{\mathsf{P}}
\newcommand{\Exp}{\mathsf{E}}
\newcommand{\bN}{\mathbb{N}}

\newcommand{\bS}{\mathbb{S}}

\begin{document}

\title{Sequential anomaly identification with observation control under generalized error metrics
}

\author{Aristomenis Tsopelakos, Georgios Fellouris
\thanks{This work was presented in part in ISIT '20, \cite{af_20}. Aristomenis Tsopelakos, Georgios Fellouris are with the Coordinated Science Laboratory, University of Illinois at Urbana-Champaign.}
}


\maketitle

\begin{abstract}
The problem of sequential anomaly detection and  identification is considered, where multiple data sources are simultaneously monitored and the goal is to  identify in real time  those, if any, that exhibit  ``anomalous'' statistical behavior. An upper bound is postulated on the number of  data sources  that can be sampled at each sampling instant, but the decision maker selects which ones to sample based on the already collected data. Thus, in this context, a policy consists  not only of a stopping rule  and a decision rule that  determine when sampling should be terminated and which  sources to identify as anomalous upon stopping, but also of a sampling rule that determines which sources to sample at each time instant subject to the sampling constraint. Two distinct formulations   are considered, which require control of  different, ``generalized''  error metrics. The first one  tolerates  a certain user-specified number of errors, of any kind,  whereas  the second tolerates   distinct, user-specified numbers of false positives and  false negatives. For each of them, a universal asymptotic lower bound on  the expected time for stopping is established as the error probabilities go to 0, and it shown to be  attained by a policy that combines the stopping and decision rules proposed in the full-sampling case with a  probabilistic sampling rule that  achieves a specific long-run sampling frequency for each source. Moreover, the optimal to a first-order asymptotic approximation expected time for stopping is compared  in simulation studies with the corresponding factor in a finite regime, and the impact of the sampling constraint and tolerance to errors is  assessed.
\end{abstract}


\section{Introduction}\label{sec:Intro}
In many scientific and engineering applications, measurements from various data sources are collected sequentially, aiming to identify  in real-time the sources, if any,  with  ``anomalous" statistical behavior. In internet security systems \cite{bol_han02}, for example, the data streams may refer to the transition rate of a link, where a meager rate warns of a possible intrusion. In finance, the data streams may refer to prices in the stock market \cite{dimson1988stock}, where we need to detect an unusual rate of return of a stock price.

Such applications, among many others, motivate the formulation of  sequential multiple testing problems  where multiple data sources are monitored sequentially, a binary hypothesis testing problem is formulated for each of them, and the goal is to simultaneously solve these testing problems as quickly as possible, while  controlling certain error metrics.  In some  works, e.g.,\cite{bar12,brt14,song2017asymptotically}, it is assumed that all sources can be monitored at each sampling instant upon stopping. In others, e.g.,  \cite{cohen2015active,huang2017active,a_prob,gur19,lambez2021anomaly,hemo_g,deshmukh,prab22,ts2023}, a sampling constraint is imposed, according to which it is  possible  to observe only a subset of sources at each time instant, but the decision maker selects which ones to sample based on the already collected data. In the latter case,  in addition to a stopping rule  and a decision rule that determine  when to stop sampling and  which  sources to identify as anomalous upon stopping,  it is also required to specify a sampling rule, which  determines the sources to be sampled at each time instant until stopping. This sampling constraint leads to a sequential multiple testing problem with adaptive sampling design, which lies in the field of ``sequential testing with controlled sensing", or ``sequential design of experiments" \cite{chernoff1959sequential}, \cite{nitinawarat2013controlled}, \cite{bessler60}, \cite{veer15}. Hence, methods and results from these works are applicable as well.

In the above works, the problem formulation requires, implicitly or explicitly,  control of the ``classical'' misclassification error rate,  i.e., the probability of at least \textit{one} error, of any kind,  or of the two  ``classical'' familywise error rates, i.e., the  probabilities of at least one false alarm and at least one missed detection. However, such error metrics can be impractical, especially  when there is a common stopping time at which the decision is made for all data sources, as in the above papers. Indeed, even with a  relatively small number of  data sources, a single ``difficult" hypothesis may determine,  even inflate, the overall time required for the  decisions to be made.  

This  inflexibility has motivated the adoption  of more lenient error metrics, which are prevalent in the fixed-sample-size  literature of multiple testing   \cite{Bartof18,Bartof20,Bartof21}. As it was shown in \cite{Bartof21}, the methodology and asymptotic optimality theory  in \cite{song2017asymptotically}, \cite{a_prob} for the sequential identification problem under classical familywise error rates remain valid with other error metrics as long as 
 \textit{the latter are  bounded above and below up to a multiplicative constant by the corresponding classical familywise error rates}. 
\footnote{The authors in \cite{Bartof21} focus in the full sampling case, but exactly the same arguments apply in  the case of sampling constraints.} This is indeed the case for many error metrics, such as  the  false discovery rate (FDR) and the false non-discovery rate (FNR) \cite{ben95}. However, it is not the case for  others, such as  the 
\textit{generalized misclassification error rate}, which
 requires control of the probability of at least $k$ errors of any kind,  and the \textit{generalized familywise  type-I and type-II errors} \cite{leh05,rom06,wol07}, which  require control of the probabilities of at least $k_1$ false alarms and at least $k_2$ missed detections.  As it was shown in  \cite{song_fel_gener}  \textit{in  the case of full sampling}, that is, when all data sources are continuously monitored and there are no sampling constraints,  the sequential identification problem with these generalized error metrics when   $k>1$ or $k_1+k_2>2$ requires a distinct methodology and pose additional mathematical challenges compared to their corresponding classical versions, where $k=1$ or $k_1=k_2=1$.

In the present work, we address the same sequential anomaly identification problem as in \cite{song_fel_gener}, but \textit{in the presence of sampling constraints}. That is, unlike \cite{song_fel_gener}, we assume that it is not possible to observe all data sources at all times. Instead, we impose an upper bound on the average number of samples collected from all sources up to stopping. For each of the two resulting problem formulations, i.e., for each of the two generalized error metrics under consideration, (i) we establish a universal lower bound on the optimal expected time for stopping, to a first-order asymptotic approximation as the corresponding error probabilities go to zero, (ii) we show that this lower bound is attained by a policy that  utilizes the stopping and decision rules in \cite{song_fel_gener}, as long as the 
long-run sampling frequency of each source is greater than or equal to a certain value, which is computed explicitly 
and depends both on the source and on the true subset of anomalous sources, (iii)  we show that the latter can be achieved by a probabilistic sampling rule. 

The methodology and the theory that we develop in the
present work is based on an interplay of  techniques,  ideas, and methods from the sequential multiple testing problem with full sampling and generalized error control in \cite{song_fel_gener}, and the sequential multiple testing problem with sampling constraints and ``classical" familywise error control in \cite{a_prob}. However, there are various interesting special features that arise in the present work, both from a technical but also from a methodological point of view.  First, there are certain technical challenges which force us to  strengthen some of our  assumptions compared to \cite{a_prob}.  Specifically, 
we  require  a stronger moment assumption on the  log-likelihood ratio statistics than the finiteness of the Kullback-Leibler information numbers, 
and  we have to  postulate a harder sampling constraint  (which is still weaker than the usual  constraint in the literature that  fixes  the  number of sources that are sampled at each time instant).    Second, while in  \cite{a_prob} it is shown that  there is no need for ``forced exploration'',  as in the  general sequential testing problem with controlled sensing \cite{chernoff1959sequential,nitinawarat2013controlled}, this turns out to be needed in this work.


The remainder of the paper is organized as follows. In Section \ref{sec:formul}, we formulate  the two problems we consider in this work. In Section \ref{prelim}, we state and  solve two auxiliary  max-min optimization problems, which play a key role in the formulation of our main results. In Section \ref{sec:lb_mis}, we state the universal asymptotic lower bound for each of the two problems. In Section \ref{sec:famil}, for each of the two problems, 
we introduce  a family of policies that satisfies the error  constraints, and we state a  criterion on the sampling rule that guarantees the asymptotic optimality of such a policy. In Section \ref{sec:ao}, we 
present a class of probabilistic sampling rules for which the aforementioned criterion is satisfied. In Section \ref{sec:sim}, we present a simulation study that illustrates our theoretical results, and in Section \ref{sec:con}, we state our conclusions and discuss future research directions.

We end this section with some notations we use throughout the paper. We use $:=$ to indicate the definition of a new quantity and $\equiv$ to indicate the equivalence of two notions. We set  $\bN := \{1, 2 \ldots, \}$, $\bN_0 := \{0\} \cup \bN$, and  $[n] := \{1, \ldots, n\}$ for  $n \in \bN$. We denote by $|A|$ the size and by $2^A$ the powerset of a set $A$, and by $A \triangle B$ the symmetric difference of two sets $A,B$. The $\lfloor  a  \rfloor$, $\lceil a  \rceil$ stand for the floor and the ceiling of a positive number $a$. The indicator function is denoted by $\mathbf{1}$. In a summation, if the lower limit is larger than the upper limit, then the summation is assumed to be equal to $0$. For positive sequences $(x_n)$ and $(y_n)$, we write $x_n\sim y_n$ when  $\lim_n (x_n/y_n)=1$, $x_n \gtrsim y_n$ when  $\liminf_n (x_n/y_n) \geq 1$, and $x_n \lesssim y_n$ when  $\limsup_n (x_n/y_n) \leq 1$. Finally, \textit{iid} stands for independent and identically distributed.  

\section{Problem formulation}\label{sec:formul}
Let $(\bS,  \cS)$ be a   measurable space and  let   $(\Omega, \cF, \Pro)$  be  a probability space that hosts $M$ independent  sequences of \textit{iid}, $\bS$-valued random elements,  $\{X_i(n),   n \in \bN\}$,   $ i \in [M]$, which are generated by $M$ distinct data sources, and $M$ independent sequences of random variables,  $\{Z_{i}(n) : n \in \bN\}$, $i \in [M]$, uniformly distributed in $(0,1)$,  to be used for randomization purposes. For each $i \in [M]$, we assume that each $X_i(n)$ has a density with respect to some $\sigma$-finite measure $\nu_i$  that is equal to either $f_{1i}$ or $f_{0i}$, and we say that  source $i$ is  ``anomalous'' if  its density is  $f_{1i}$.  We  denote by $\Pro_A$  the underlying probability measure, and by $\Exp_A$ the corresponding expectation, when the subset of anomalous sources is  $A \subseteq [M]$. We simply write $\Pro$ and $\Exp$ whenever the identity of the subset of anomalous sources is not relevant. 

The  problem we consider is the identification of the anomalous sources, if any, on the basis of sequentially acquired observations from all sources \textit{when it is not possible to observe all of them at every sampling instant}. Specifically,  we have to specify (i) the random  time,  $T$, at which sampling is terminated, (ii) the subset of sources, $\Delta$, that are declared as anomalous upon stopping, and, for each  $n \leq T$,  (iii) the subset, $R(n)$,  of sources that are sampled  at time $n$.  At any time instant, the decision whether to stop or not,  as well as the subsets of sources that are identified as anomalous in the former case, or sampled next in the latter, must be determined based on the already collected data. 

Therefore, we say that $R:=\{R(n): n \in \bN\}$  is a \textit{sampling rule}  if  $R(n)$ is $\cF^{R}_{n-1}$-measurable  for every $n \in \bN$,  where
\begin{align} \label{filtration}
	\cF^R_n := \begin{cases}
		\sigma\left( \cF^R_{n-1},\, Z(n),\, \{X_i(n) \,:\,  i \in R(n)\} \right),\quad  &n \in \bN, \\
		\sigma(Z(0)), \quad  &n=0,
	\end{cases}
\end{align}  
and $Z(n):=(Z_1(n),\ldots, Z_{M}(n))$. 
Moreover,  we say that the triplet $(R, T, \Delta)$ is a \textit{policy} if 
\begin{itemize}
\item[(i)] $R$ is a sampling rule, 
\item[(ii)]  $T$ is a stopping time with respect to filtration $\{\cF^R_n\,:\, n \in \mathbb{N} \}$,  
\item[(iii)]  $\Delta$ is $\cF^R_T$-measurable,  i.e.,
\begin{equation*}
	\{T=n, \Delta =D\}\in \cF^R_n, \quad \forall \; n \in \bN \quad \text{and} \quad  D \subseteq [M],
\end{equation*}
\end{itemize}
in which case we refer to  $T$ as a \textit{stopping rule}  and to $\Delta$ as  a \textit{decision rule}. 

We denote by $\cC$ the family of all policies, and we focus on  policies that satisfy a sampling constraint and control the probabilities of certain types of error.  To define these, we need to introduce some additional notation.  Thus, for any sampling rule $R$ and time instant $n \in \bN$, we denote by $R_i(n)$ the indicator of whether source $i$ is sampled at time $n \in \bN$, i.e., $$R_i(n) := \mathbf{1}\{ i \in R(n) \},\quad i \in [M],$$ 
by  $\pi_{i}^R(n)$ the proportion of times source $i$ has been sampled in the first $n$ time instants, i.e.,
\begin{align}\label{nri}
	\pi_i^R(n) :=  \frac{1}{n}\sum_{m=1}^n  R_i(m),
\end{align}
and we note that the average number of observations from all sources in the first $n$ time instants is 
\begin{equation}\label{nri2}
	\frac{1}{n}  \sum_{m=1}^n  |R(m)|=  \frac{1}{n} \sum_{i=1}^{M} \sum_{m=1}^n  R_i(m)= 
	\sum_{i=1}^{M} \pi_i^R(n).
\end{equation}

\subsection{The sampling constraint}

For any real number $K$ in $(0,M]$, we say that a policy $(R,T,\Delta)$ belongs to $\cC(K)$ if the average number of observations from all sources until stopping is less than or equal to $K$, i.e.,
\begin{equation}\label{samp_cons}
	\frac{1}{T}  \sum_{m=1}^T  |R(m)|= 
	\sum_{i=1}^{M} \pi_i^R(T) \leq  K.
\end{equation}

This sampling constraint is clearly satisfied when at most $\lfloor K \rfloor$ sources are sampled at each time instant up to stopping, i.e., when 
\begin{equation}\label{flk}
	 |R(m)| \leq \lfloor K \rfloor, \quad \forall \, m \leq T,
\end{equation}
whereas it implies 
the sampling constraint considered in  \cite{a_prob}, 
\begin{equation}\label{samp_cons_looser}
	 \Exp \left[  \sum_{m=1}^T  |R(m)| \right] \leq K \,  \Exp[T].
\end{equation}
In other words, the sampling constraint  \eqref{samp_cons} is looser than \eqref{flk}, but stricter than \eqref{samp_cons_looser}.

\subsection{The error constraints}
We consider two  types of error control,  which lead to two  distinct problem formulations. We characterize them both as  ``generalized'',  as they generalize their corresponding ``classical'' versions of misiclassification error rate and familywise error rates.

\subsubsection{Control of  generalized misclassification error rate}
For any $K \in (0,M]$, $k \in [M]$, and  $\alpha \in (0,1)$, we say that a policy $(R,T, \Delta)$ in $\cC(K)$ belongs to $\mathcal{C}(\alpha; k,K)$  if the probability of at least $k \in [M]$ errors of any kind is at most $\alpha$, i.e.,
\begin{equation}\label{error_const}
	\Pro_{A}(|A \triangle \Delta | \geq  k ) \leq \alpha, \quad \forall \, A \subseteq [M],
\end{equation}
and we denote by $\mathcal{J}_{A}(\alpha; k,K)$ the smallest possible expected time for stopping in
$\mathcal{C}(\alpha; k, K)$, when the subset of anomalous sources is $A \subseteq [M]$, i.e., 
\begin{equation} \label{optimal_unknown}
	\mathcal{J}_{A}(\alpha; k,K):=\inf \limits_{(R,T,\Delta) \in \cC(\alpha; k,K)}\Exp_{A}[T].
\end{equation}
The first problem we consider in this paper is  to evaluate 
$\mathcal{J}_{A}(\alpha; k,K)$  to a first-order asymptotic approximation as  $\alpha \to 0$ for any $A \subseteq[M]$, and to find a policy that achieves $\mathcal{J}_{A}(\alpha; k,K)$   in this  asymptotic sense, simultaneously for every $A \subseteq [M]$.\\

\subsubsection{Control of  generalized familywise error rates}
For any $K \in (0, M]$,  $k_{1},\, k_{2} \in [M]$,  and $\alpha, \, \beta \in (0,1)$,  such that  $\alpha + \beta < 1$ and $k_{1} + k_{2} \leq M$, we say that  a policy $(R,T, \Delta)$ in $\cC(K)$ belongs to   $\mathcal{C}(\alpha,\beta; k_{1},k_{2}, K)$ if  the probability of at least $k_{1}$ false positives does not exceed  $\alpha$ and the probability of at least $k_{2}$  false negatives does not exceed $\beta$, i.e., 
\begin{equation}\label{b_error_const}
	\Pro_{A}(|\Delta \setminus A| \geq k_{1}) \leq \alpha, \qquad \text{and}  \qquad   \Pro_{A}(|A \setminus \Delta| \geq k_{2}) \leq \beta, \qquad \forall \, A \subseteq [M],
\end{equation}
and we denote by $\mathcal{J}_{A}(\alpha,\beta; k_1,k_2,K)$ the smallest expected time for stopping in $\mathcal{C}(\alpha,\beta; k_{1},k_{2}, K)$ when the subset of anomalous sources is $A\subseteq [M]$, i.e.,
\begin{equation}\label{optimal_both_kin}
	\mathcal{J}_{A}(\alpha,\beta;k_1,k_2,K):=\inf\limits_{(R,T,\Delta) \in \mathcal{C}(\alpha,\beta;k_{1},k_{2},K)}\Exp_{A}[T].
\end{equation}
The second problem we consider in this paper is to  evaluate   $\mathcal{J}_{A}(\alpha,\beta;k_1,k_2,K)$ to a first-order asymptotic approximation as  $\alpha, \beta  \to 0$, for any $A \subseteq [M]$,  and to find a family of  policies that achieves   $\mathcal{J}_{A}(\alpha,\beta;k_1,k_2,K)$ in this asymptotic sense, simultaneously for every $A \subseteq [M]$.\\

\begin{remark}
	As we mentioned in the Introduction, both problems have been solved for the full sampling case, i.e., \textit{when all sources are observed at each instant}, in \cite{song_fel_gener}. On the other hand, in the presence of sampling constraints, neither of them  has  been considered beyond the special case of ``classical" misclassification error rate ($k=1$) and ``classical" familywise error rates ($k_1=k_2=1$) in \cite{a_prob}.
\end{remark}

\subsection{Distributional Assumptions}
For each $i \in [M]$, the Kullback-Leibler (KL) divergences of $f_{1i}$ and $f_{0i}$ are assumed to be positive and finite, i.e., for each $i \in [M]$ we have:
\begin{equation}\label{IJ}
\begin{aligned}
I_{i}& := \int_{\bS}  \log \left( f_{1i} / f_{0i}  \right) \, f_{1i} \, d \nu_i \, \in (0, \infty), \\ 
J_{i} &:= \int_{\bS}  \log \left( f_{0i} / f_{1i}  \right) \, f_{0i} \,  d \nu_i \, \in (0, \infty).
\end{aligned}	
\end{equation}

To establish asymptotic lower bounds for \eqref{optimal_unknown} and \eqref{optimal_both_kin}, 
we assume that
\begin{align}\label{momr}
\begin{split}
\sum_{i=1}^M \int_{\bS} |g_{i}|\log^{+} (|g_{i}|)\, f_{1i} \,  d \nu_i  &< \infty, \qquad \\
\sum_{i=1}^M \int_{\bS} |g_{i}|\log^{+}(|g_{i}|)\, f_{0i} \,  d \nu_i  &< \infty,
\end{split}
\end{align}
where $g_i :=  \log \left( f_{1i} / f_{0i}  \right)$.  This assumption is not needed neither in the full sampling case considered in \cite{song_fel_gener}, nor in the case of classical error control ($k=1$ or $k_1=k_2=1$) under sampling constraints in \cite{a_prob}.  Nevertheless,  it is weaker than the typical assumption in the sequential controlled sensing literature (e.g.,\cite{chernoff1959sequential}, \cite{nitinawarat2013controlled}), according to which it is required that 
\begin{align} \label{momr0}
\begin{split}
\sum_{i=1}^M \int_{\bS}   |g_i|^{\mathfrak{p}} \, f_{1i} 
\; d\nu_i  &< \infty, \\
\sum_{i=1}^M \int_{\bS}   |g_i|^{\mathfrak{p}}  \, f_{0i} \; d\nu_i & < \infty,
\end{split}
\end{align}
holds for  $\mathfrak{p}=2$. However, in order to show that the asymptotic lower bound in \eqref{optimal_unknown} and  \eqref{optimal_both_kin}  can be achieved, in certain cases
we will   need to require that \eqref{momr0} holds for some  $\mathfrak{p}>4$.

\section{Two  max-min optimization problems}\label{prelim}
In this section, we formulate and solve two auxiliary max-min optimization problems. These will be  used in Section \ref{sec:lb_mis} to express the asymptotic lower bounds for \eqref{optimal_unknown} and \eqref{optimal_both_kin}, and in Section \ref{sec:famil} to design procedures that achieve these lower bounds. 

To be specific,  we denote by  
$\boldsymbol{L}:=\{ L_{i} \,:\, i \in [|\boldsymbol{L}|] \}$ an ordered set  of positive numbers, i.e.,
	      \begin{equation} \label{L1}
	      L_{0}:=0 < L_{1} \leq \ldots \leq L_{|\boldsymbol{L}|},
	      \end{equation}
 for  each $i \in [|\boldsymbol{L}|]$, we denote by $\widetilde{L}_{i}$ the harmonic mean of the $|\boldsymbol{L}|-i+1$ largest elements in $\boldsymbol{L}$, i.e., 
	      \begin{equation} \label{L2}
	      \widetilde{L}_{i}:= \frac{|\boldsymbol{L}|-i+1}{\sum_{u=i}^{|\boldsymbol{L}|} (1/L_{u})}, \quad i \in [|\boldsymbol{L}|],
	      \end{equation}
and we also set $\widetilde{L}_{|\boldsymbol{L}|+1}:=\infty$. Then, assuming that  $|\boldsymbol{L}| \leq M$, we introduce the following function,
\begin{equation}\label{mn_m}
\cV(\boldsymbol{c};\kappa,\boldsymbol{L}) := \min_{U \subseteq [|\boldsymbol{L}|]: \, |U|= \kappa} \, \sum_{i \in U} c_i \, L_{i},
\end{equation}
where
\begin{itemize}	
   \item $\kappa$ is a positive integer in $[|\boldsymbol{L}|]$,
	
   \item $\boldsymbol{c}:=(c_1,\ldots,c_{|\boldsymbol{L}|},0,\ldots,0)$ is a vector in 
	\begin{equation}\label{domain}
	\cD(K) := \left\{(c_1, \ldots, c_M) \in [0,1]^M: \sum_{i=1}^{M}c_{i} \leq K \right\}.
	\end{equation}
\end{itemize}

The constants $M$ and $K$ are  defined  as in the previous section, i.e., $M$ is a  positive integer, and $K$ a real number in $(0,M]$. 

\subsection{Optimization Problem I}
Let $\boldsymbol{L}$ be an ordered set of size $|\boldsymbol{L}| \leq M$, and $\kappa \in [|\boldsymbol{L}|]$. The first max-min optimization problem we consider is
\begin{equation}\label{mm_st}
V(\kappa,K,\boldsymbol{L}):= \max_{\boldsymbol{c} \in \cD(K)} \cV(\boldsymbol{c};\kappa,\boldsymbol{L}).
\end{equation}
In the following lemma, we provide an expression for the value $V(\kappa,K,\boldsymbol{L})$ of the max-min optimization problem \eqref{mm_st}, as well as for the maximizer of \eqref{mm_st} with the minimum $\cL^{1}$ norm.\\

\begin{lemma}\label{mm_g}
	The value of the max-min optimization problem \eqref{mm_st} is equal to the expression
	\begin{equation}\label{st0mb}
	V(\kappa,K,\boldsymbol{L})=
	\begin{cases}	
	(\kappa-u)\,\frac{K}{|\boldsymbol{L}|-u}\,\widetilde{L}_{u+1}, \quad &\mbox{if}\quad v=0,\\
	x\, L_{v-1} + \sum_{i=v}^{u} L_{i} + (\kappa-u)\, y\, \widetilde{L}_{u+1}, \quad &\mbox{if}\quad v \geq 1,
    \end{cases}
    \end{equation}
	where $x, y$ are real numbers in $[0,1)$, and $u, v$ are integer numbers in $[0,\kappa]$ such that $v \leq u$. The values of $x$, $y$, $v$, $u$ are determined by Algorithm \ref{alg} presented in Appendix \ref{miscl}. The maximizer of \eqref{mm_st} with the minimum $\cL^{1}$ norm is given by
	\begin{equation}\label{c_pr}
	\boldsymbol{c}'(\kappa,K,\boldsymbol{L}):=(c'_1,\ldots,c'_{|\boldsymbol{L}|},0,\ldots,0),
	\end{equation}
	where 
	\begin{itemize}
		\item for all $i \in \{u+1,\ldots,|\boldsymbol{L}|\}$ we have
		\begin{equation}\label{c1} 
		c'_{i} \, L_i =
		\begin{cases}
		y \, \widetilde{L}_{u+1},\;\, \qquad \text{if} \quad u < \kappa, \\
		L_{\kappa}, \qquad  \qquad \, \text{if}  \quad  u=\kappa,
		\end{cases}
		\end{equation}
		
		\item if $v=0$ then
		\begin{equation}
		c'_{i} = 0, \qquad \mbox{for all} \quad i \, \in \, \{1,\ldots,u\},
		\end{equation}
		
		\item if $v \geq 1$ then 
		\begin{equation}\label{c2}
		c'_{i} = 1, \qquad \mbox{for all} \quad i \, \in \, \{v,\ldots,u\},
		\end{equation}
		
		\item if $v \geq 2$ then
		\begin{equation}
		c'_{v-1} =
		\begin{cases}
		x, \qquad \text{if} \quad x >0,\\
		0, \qquad \text{if} \quad x =0,
		\end{cases}
		\end{equation}

        \item if $v \geq 3$ then
		\begin{equation}\label{c4}
		c'_{i} = 0,  \qquad \mbox{for all} \quad i \leq v-2.
		\end{equation}
	\end{itemize}
\end{lemma}

\begin{IEEEproof}
Appendix \ref{miscl}.	
\end{IEEEproof}

\begin{remark}\label{rem}
	In the symmetric case where $L_{i}=L$ for all $i \in [|\boldsymbol{L}|]$, then 
	\begin{equation*}
	V(\kappa ,K,\boldsymbol{L})=\kappa \, (K/|\boldsymbol{L}|) \, L,
	\end{equation*}
	and $c'_{i} = (K/|\boldsymbol{L}|) \wedge 1$ for all $i \in [|\boldsymbol{L}|]$.
\end{remark}

\begin{remark}
Based on the size of $K$ we distinguish the following cases on the form of $V(\kappa,K,\boldsymbol{L})$.
\begin{itemize}
	\item If $K$ is relatively large, i.e.,
	      \begin{equation}\label{x1}
	      K \geq \kappa + L_{\kappa}\sum_{i=\kappa +1}^{|\boldsymbol{L}|} 1/L_i,
	      \end{equation}
	      then
	     \begin{equation}\label{vf}
	     V(\kappa,K,\boldsymbol{L})= \sum_{i=1}^{\kappa} L_i,
	     \end{equation}
	     which is the largest possible value that $V(\kappa,K,\boldsymbol{L})$ can take over all possible values of $K \leq M$. In this case $x=y=0$, $v=1$, $u=\kappa$.
	 
	\item If $K$ is relatively small, i.e.,
	      \begin{equation}\label{x2}
	      K < L_{u^* +1}\sum_{i=u^* +1}^{|\boldsymbol{L}|} 1/L_i,
	      \end{equation}
	      then
	      \begin{equation*}
	      V(\kappa,K,\boldsymbol{L})=(\kappa -u^*)\; \frac{K}{|\boldsymbol{L}|-u^*} \;\widetilde{L}_{u^* +1},
	      \end{equation*}    
	      where $u^*$ is a quantity defined as
	      \begin{equation*}
	      u^* := \max\left\{u \in \{0,\ldots,\kappa-1\} \,:\, \frac{\kappa - u}{|\boldsymbol{L}|-u}\;\widetilde{L}_{u+1} \geq L_{u}\right\}.
	      \end{equation*}
	      In this case $x=0$, $u=u^*$, $v=0$, and $y=K/(|\boldsymbol{L}|-u^*)$.
	      
	\item If $K$ is between the values given in \eqref{x1}, \eqref{x2} the $V(\kappa,K,\boldsymbol{L})$ has the general form described in Lemma \ref{mm_g}.
\end{itemize}
\end{remark}

\begin{remark}
The only case we can have a maximizer $(\tilde{c}_1,\ldots,\tilde{c}_{|\boldsymbol{L}|})$ which is not equal to the maximizer with the minimum $\cL^1$ norm, is when \eqref{x1} holds by strict inequality. The only difference between the two maximizers is that for $(\tilde{c}_1,\ldots,\tilde{c}_{|\boldsymbol{L}|})$ there is $i \in \{ \kappa +1,\ldots,|\boldsymbol{L}|\}$ such that
\begin{equation}\label{cli}
\tilde{c}_i L_i > L_{\kappa}.
\end{equation}
Although \eqref{cli} can hold for some $i \in \{ \kappa +1,\ldots,|\boldsymbol{L}|\}$, this does not change the value of $V(\kappa,K,\boldsymbol{L})$ in \eqref{vf}.
\end{remark}

\subsection{Optimization Problem II}
Let $\boldsymbol{L}_{1}:= \{ L_{1,i}: i \in [|\boldsymbol{L}_{1}|]\}$ and $\boldsymbol{L}_{2}:=\{ L_{2,i}: i \in [|\boldsymbol{L}_{2}|] \}$ be two ordered sets such that $|\boldsymbol{L}_{1}|+|\boldsymbol{L}_{2}| \leq M$,  let $\kappa_{1}$, $\kappa_{2}$ be two positive integers such that $\kappa_{1} \in [|\boldsymbol{L}_{1}|]$, $\kappa_{2} \in [|\boldsymbol{L}_{2}|]$, and let $r$ be an arbitrary positive number. The second max-min optimization problem we consider in this section is more complex than the first, and it has the form
\begin{equation}\label{mm_st2}
W(\kappa_{1},\kappa_{2},K,\boldsymbol{L}_{1},\boldsymbol{L}_{2}, r):=\max_{\boldsymbol{c} \in \cD(K)}\min\left\{\cV(\hat{\boldsymbol{c}};\kappa_{1},\boldsymbol{L}_{1}), r \; \cV(\check{\boldsymbol{c}};\kappa_{2},\boldsymbol{L}_{2})\right\},
\end{equation}
where the function  $\cV$ is defined in \eqref{mn_m}, and 
 $\boldsymbol{c}:=(\hat{\boldsymbol{c}},\check{\boldsymbol{c}},\boldsymbol{0})$, the size of $\boldsymbol{0}$ being $M-|\boldsymbol{L}_{1}|-|\boldsymbol{L}_{2}|$, and 
\begin{align}
\begin{split}
\hat{\boldsymbol{c}}&:=(\hat{c}_{1},\ldots,\hat{c}_{|\boldsymbol{L}_{1}|}),\\ \check{\boldsymbol{c}}&:=(\check{c}_{1},\ldots,\check{c}_{|\boldsymbol{L}_{2}|}).
\end{split}
\end{align}

As we show in the following lemma, the value $W(\kappa_{1},\kappa_{2},K,\boldsymbol{L}_{1},\boldsymbol{L}_{2}, r)$ of the max-min optimization problem \eqref{mm_st2} is equal to the value of the following optimization problem,
\begin{equation}\label{rt}
\max_{(K_1, K_2)} V(\kappa_1,K_1,\boldsymbol{L}_1),
\end{equation}
such that the following two constraints hold 
\begin{equation}\label{ck}
\begin{aligned}
K_1 + K_2 &\leq K,\\
V(\kappa_1,K_1,\boldsymbol{L}_1) &= r \, V(\kappa_2,K_2,\boldsymbol{L}_2).
\end{aligned}
\end{equation}

\begin{definition}
We denote by $(K^*_1, K^*_2)$ the maximizer of the constrained optimization problem \eqref{rt} with the minimum $\cL^{1}$ norm, i.e., the minimum sum $K_1 + K_2$, among all maximizers. Based on Lemma \ref{mm_g}, we denote by $x_1, y_1, u_1, v_1$ the parameters such that
\begin{equation*}
V(\kappa_1,K^{*}_1,\boldsymbol{L}_1) = 
\begin{cases}
(\kappa_1 - u_1)\,\frac{K^{*}_1}{|\boldsymbol{L}_{1}|-u_1}\, \widetilde{L}_{1,u_1 +1}, \quad & \mbox{if}\quad v_1 =0,\\
x_1 L_{1,v_1 -1} + \sum_{i=v_1}^{u_1} L_{1,i} + (\kappa_1 - u_1)\, y_1 \, \widetilde{L}_{1,u_1 +1}, \quad & \mbox{if}\quad v_1 \geq 1,
\end{cases}
\end{equation*}
and by $x_2, y_2, u_2, v_2$ the parameters such that
\begin{equation*}
V(\kappa_2,K^{*}_2,\boldsymbol{L}_2) = 
\begin{cases}
(\kappa_2 - u_2)\,\frac{K^{*}_2}{|\boldsymbol{L}_{2}|-u_2}\, \widetilde{L}_{2,u_2 +1}, \quad & \mbox{if}\quad v_2 =0,\\
x_2 L_{2,v_2 -1} + \sum_{i=v_2}^{u_2} L_{2,i} + (\kappa_2 - u_2)\, y_2 \, \widetilde{L}_{2,u_2 +1}, \quad & \mbox{if}\quad v_2 \geq 1.
\end{cases}
\end{equation*}
The parameters $x_1, y_1, u_1, v_1$, and $x_2, y_2, u_2, v_2$ can be computed by applying Algorithm \ref{alg} for each case.
\end{definition}
The parameters $x_1, y_1, u_1, v_1$, and $x_2, y_2, u_2, v_2$ are used in the expression of the maximizers of \eqref{mm_st2} with the minimum $\cL^1$ norm.

\begin{lemma}\label{equiv_opt}
The value of the max-min optimization problem \eqref{mm_st2} is equal to 
\begin{equation}\label{wf}
\begin{aligned}
W(\kappa_{1},\kappa_{2},K,\boldsymbol{L}_{1},\boldsymbol{L}_{2},r) &= V(\kappa_1,K^{*}_1,\boldsymbol{L}_1) \\
 &=r\,V(\kappa_2,K^{*}_2,\boldsymbol{L}_2).
\end{aligned}
\end{equation}
The maximizer of \eqref{mm_st2} with the minimum $\cL^1$ norm is given by
\begin{equation}\label{c_pr2}
\boldsymbol{c}'(\kappa_{1},\kappa_{2},K,\boldsymbol{L}_{1},\boldsymbol{L}_{2}, r):=\left(\hat{c}'_{1},\ldots,\hat{c}'_{|\boldsymbol{L}_{1}|},\check{c}'_{1},\ldots,\check{c}'_{|\boldsymbol{L}_{2}|},0,\ldots,0 \right),
\end{equation}
where
\begin{itemize}
	\item for all $i \in \{u_1 +1,\ldots,|\boldsymbol{L}_1|\}$ we have
	\begin{equation}\label{ts}
	\hat{c}'_{i} \, L_{1,i} =
	\begin{cases}
	y_1 \, \widetilde{L}_{1,u_1 +1},\, \qquad \text{if} \quad u_1 < \kappa_1, \\
	L_{1,\kappa_1}, \qquad  \qquad \, \text{if}  \quad  u_1=\kappa_1,
	\end{cases}
	\end{equation}
	
	\item if $v_1 =0$ then
	\begin{equation}
	\hat{c}'_{i} = 0, \qquad \mbox{for all} \quad i \, \in \, \{1,\ldots,u_1\},
	\end{equation}
	
	\item if $v_1 \geq 1$ then 
	\begin{equation}
	\hat{c}'_{i} = 1, \qquad \mbox{for all} \quad i \, \in \, \{v_1,\ldots,u_1\},
	\end{equation}
	
	\item if $v_1 \geq 2$ then
	\begin{equation}
	\hat{c}'_{v_1 -1} =
	\begin{cases}
	x_1, \;\, \quad \text{if} \quad x_1 >0,\\
	0, \qquad \text{if} \quad x_1 =0,
	\end{cases}
	\end{equation}
	
	\item if $v_1 \geq 3$ then
	\begin{equation}
	\hat{c}'_{i} = 0,  \qquad \mbox{for all} \quad i \leq v_1 -2,
	\end{equation}
\end{itemize}
and
\begin{itemize}
	\item for all $i \in \{u_2 +1,\ldots,|\boldsymbol{L}_2|\}$ we have
	\begin{equation} 
	\check{c}'_{i} \, L_{2,i} =
	\begin{cases}
	y_2 \, \widetilde{L}_{2,u_2 +1},\, \qquad \text{if} \quad u_2 < \kappa_2, \\
	L_{2,\kappa_2}, \qquad  \qquad \, \text{if}  \quad  u_2=\kappa_2,
	\end{cases}
	\end{equation}
	
	\item if $v_2 =0$ then
	\begin{equation}
	\check{c}'_{i} = 0, \qquad \mbox{for all} \quad i \, \in \, \{1,\ldots,u_2\},
	\end{equation}
	
	\item if $v_2 \geq 1$ then 
	\begin{equation}
	\check{c}'_{i} = 1, \qquad \mbox{for all} \quad i \, \in \, \{v_2,\ldots,u_2\},
	\end{equation}
	
	\item if $v_2 \geq 2$ then
	\begin{equation}
	\check{c}'_{v_2 -1} =
	\begin{cases}
	x_2, \;\, \quad \text{if} \quad x_2 >0,\\
	0, \qquad \text{if} \quad x_2 =0,
	\end{cases}
	\end{equation}
	
	\item if $v_2 \geq 3$ then
	\begin{equation}\label{te}
	\check{c}'_{i} = 0,  \qquad \mbox{for all} \quad i \leq v_2 -2.
	\end{equation}
\end{itemize}   	      
\end{lemma}

\begin{IEEEproof}
	Appendix \ref{miscl}.
\end{IEEEproof}

\section{Universal asymptotic lower bounds}\label{sec:lb_mis}
In this section, we fix an arbitrary $A \subseteq [M]$, and  $k,k_1,k_2,K$ as in Section \ref{sec:formul}, and we establish a universal asymptotic lower bound for \eqref{optimal_unknown} and \eqref{optimal_both_kin}, as $\alpha \to 0$, and $\alpha, \beta \to 0$ respectively, under the moment assumption \eqref{momr}.

\subsection{The case of generalized misclassification error rate}\label{mmiscl}
To state the asymptotic lower bound for $\cJ_A(\alpha; k, K)$ as $\alpha \to 0$, we need to introduce some additional notation. Thus, we denote by 
\begin{equation}\label{F(A)}
\boldsymbol{F}(A):=\{ F_{i}(A) \,:\, i \in [M] \}
\end{equation}
the ordered set that consists of the Kullback-Leibler numbers in $\{ I_{i},\, J_{j}\, : \, i \in A,\, j \notin A \}$. In particular, for each $i \in [M]$, $F_{i}(A)$ is the $i^{th}$ smallest element in $\boldsymbol{F}(A)$, and can be interpreted  as a measure of the  difficulty of the $i^{th}$ most difficult testing problem. 

The overall difficulty of the testing problem is determined by the quantity $V(k,K,\boldsymbol{F}(A))$, as described in the following theorem.

\begin{theorem}\label{th:lower_bound1} 
Suppose \eqref{momr} holds. As $\alpha \to 0$, we have	
	\begin{equation}
	\cJ_A(\alpha; k, K) \gtrsim \frac{|\log \alpha|}
	{V(k,K,\boldsymbol{F}(A))},
	\end{equation}
	where $V$ is defined in  \eqref{mm_st}.
\end{theorem}

\begin{IEEEproof}
	Appendix \ref{LBs}.
\end{IEEEproof}

\begin{remark}
Consider the homogeneous and symmetric setup where the difficulty is the same across all testing problems,  in the sense that  
$$I_{i}=J_{j}=I, \quad \forall \; i,j \in [M].$$ 
Then, for every $A \subseteq [M]$  we have
$$F_i(A)=I, \quad \forall \; i \in [M].$$
Consequently, by Remark \ref{rem} we can see that
\begin{equation*}
V(k,K,\boldsymbol{F}(A))=k \, (K/M) \, I.
\end{equation*}
\end{remark}

\subsection{The case of generalized familywise error rates}\label{gener_fam}
We next establish an asymptotic lower bound for $\mathcal{J}_{A}(\alpha,\beta;k_1,k_2,K)$ as $\alpha \to 0$ and / or $\beta \to 0$. For this, we need to introduce the following notation.
\begin{enumerate}[(i)]
\item If $A \neq \emptyset$,  for each $i \in [|A|]$  we denote
\begin{itemize}
\item by  $I_{i}(A)$ the $i^{th}$ smallest element in $\{I_{j} \,:\, j \in A\}$, 
\item  by
 $\widetilde{I}_{i}(A)$ the harmonic mean of the $|A|-i+1$ largest elements in $\{I_{j} \,:\, j \in A\}$.
 \end{itemize}
Moreover,  we denote by $\boldsymbol{I}(A)$ the ordered set that consists of the Kullback-Leibler numbers in $\{I_{j} \,:\, j \in A\}$, i.e., 
$$\boldsymbol{I}(A):= \left\{I_{i}(A) \,:\, i \in [|A|] \right\},$$	
and for each $l \in \{0,\ldots,|A|-1\}$, we denote by $\boldsymbol{I_{l}}(A)$ the  set that consists of the $|A|-l$ largest elements in $\{I_{j} \,:\, j \in A\}$, i.e., 
$$\boldsymbol{I_{l}}(A):=\{ I_{i}(A) \,:\, l<i\leq|A|\}.$$ 
	
	\item If $A \neq [M]$,    for each $i \in [|A^c|]$ we denote
	 \begin{itemize}
	 \item  by $J_{i}(A)$ the $i^{th}$ smallest element in $\{J_{j} \,:\, j \in A^c\}$, 
	 \item  by $\widetilde{J}_{i}(A)$  the harmonic mean of the largest $|A^c|-i+1$ elements in $\{J_{j} \,:\, j \in A^c\}$.
	 \end{itemize}
	 Moreover, we denote by $\boldsymbol{J}(A)$ the ordered set that consists of the Kullback-Leibler numbers in $\{J_{j} \,:\, j \in A^c\}$, i.e., 
$$\boldsymbol{J}(A):= \left\{J_{i}(A) \,:\, i \in [|A^c|] \right\},$$
	and for each $l \in \{0,\ldots,|A^c|-1\}$, we denote by $\boldsymbol{J_{l}}(A)$ the ordered set that consists of the $|A^c|-l$ largest elements in $\{J_{j} \,:\, j \in A^c\}$, i.e., 
 $$\boldsymbol{J_{l}}(A):=\{J_{i}(A) \,:\, l<i\leq |A^c|\}.$$
\end{enumerate}

We state first  the asymptotic lower bound for $\mathcal{J}_{A}(\alpha,\beta;k_1,k_2,K)$ when $A=\emptyset$ and  when $A=[M]$, as $\beta \to 0$ and $\alpha \to 0$, respectively. 

\begin{theorem}\label{tf1}
Suppose \eqref{momr} holds. 
\begin{itemize}
	\item Let $A=\emptyset$. For any given $\alpha \in (0,1)$, as $\beta \to 0$ we have
	      \begin{equation}
	      \mathcal{J}(\alpha,\beta;k_1,k_2,K) \gtrsim \frac{|\log \beta|}{V(k_2,K,\boldsymbol{J_{k_1-1}}(A))}.
	      \end{equation}
	
	\item Let $A=[M]$. For any given $\beta \in (0,1)$, as $\alpha \to 0$ we have
	      \begin{equation}
	      \mathcal{J}(\alpha,\beta;k_1,k_2,K) \gtrsim \frac{|\log \alpha|}{V(k_1,K,\boldsymbol{I_{k_2-1}}(A))}.
	      \end{equation}
\end{itemize}	
where $V$ is defined in \eqref{mm_st}.
\end{theorem}

\begin{IEEEproof}
	Appendix \ref{LBs}.\\
\end{IEEEproof}
 
\begin{remark}
If $A=\emptyset$ the $k_1 -1$ sources with the smallest KL numbers in $\boldsymbol{J}(A)$ are not considered in the evaluation of the difficulty of the testing problem. This is because we can intentionally misclassify these $k_1 -1$ sources as anomalous without exceeding the tolerance level of $k_1$ false alarm errors in order to reduce the expected stopping time. The respective remark holds for the case $A=[M]$.  
\end{remark} 

We continue with the asymptotic lower bound  when  $0<|A|<M$ as $\alpha,\, \beta \to 0$ so that
\begin{equation}\label{r}
	|\log \alpha|\sim r |\log \beta|, \quad \mbox{for some }\; r \in (0, \infty).
\end{equation}
For this, we need to introduce the following definition. 
	
\begin{definition}\label{def_v_A}
We denote by $v_{A}(k_{1},k_{2},K,r)$ the maximum of the following quantities. Each quantity is included to the overall maximum given that the respective condition is satisfied. We recall that $k_1 \leq |A|$ or $k_2 \leq |A^c|$, because otherwise we would have $k_1 + k_2 > |A| + |A^c|=M$ which contradicts the initial assumption $k_1 +k_2 \leq M$.
\begin{itemize}
	\item If $k_2 \leq |A^c|$ we include the
	      \begin{equation}\label{v1}
	      \max\{W(k_1 -l, k_2, K,\boldsymbol{I}(A),\boldsymbol{J_{l}}(A), r)\}
	      \end{equation}
	      over all $l \in \{(k_1 -|A|)^+, \ldots, (k_1 -1) \wedge (|A^c|-k_2)\}$.
	
	\item  If $k_1 \leq |A|$ we include the
	       \begin{equation}\label{v2}
	       \max\{W(k_1,k_2-l,K,\boldsymbol{I_{l}}(A),\boldsymbol{J}(A),r)\}
	       \end{equation}
	       over all $l \in \{(k_2-|A^c|)^{+},\ldots,(k_2 -1)\wedge(|A|-k_1)\}$.
	
	\item If $k_1 -1 \geq |A^c|-k_2 +1$ we include the
	      \begin{equation}\label{v3}
	      V(k_1-(|A^c|-k_2 +1)^+,K,\boldsymbol{I}(A)).
	      \end{equation}
	      
	\item If $k_2 -1 \geq |A|-k_1 +1$ we include the
	      \begin{equation}\label{v4}
	      r\, V(k_2-(|A|-k_1 +1)^+,K,\boldsymbol{J}(A)).
	      \end{equation}
\end{itemize}
where $W$ is defined in \eqref{mm_st2}, and $V$ in \eqref{mm_st}.
\end{definition}
	
\begin{remark}
We note that for any $A \subseteq [M]$,  it holds
\begin{align}
(k_1 -|A|)^+ &\leq (k_1 -1) \wedge (|A^c|-k_2),\label{aa1}   \\
(k_2 -|A^c|)^+ &\leq (k_2 -1) \wedge (|A|-k_1). \label{aa2}
\end{align}
Indeed for \eqref{aa1} we observe that if $k_1 \leq |A|$ then $(k_1 -|A|)^+ =0$ and the result is evident, whereas if $k_1 > |A|$ it holds
\begin{align}
k_1 - |A| &\leq k_1 -1, \label{ab1} \\
k_1 - |A| &\leq |A^c| - k_2, \label{ab2}
\end{align}
where \eqref{ab1} follows by the fact that $0 < |A| < M$, and \eqref{ab2} holds because we have assumed $k_1 + k_2 \leq M$. Similarly, we can verify that \eqref{aa2} holds.
\end{remark}	
	
The difficulty of the testing problem is determined by the quantity $v_{A}(k_{1},k_{2},K,r)$ as described in the following theorem.	

\begin{theorem}\label{lower_bound_bk}
Suppose \eqref{momr} holds, and let $0<|A|<M$. As $\alpha,\, \beta \to 0$ so that \eqref{r} holds, 
	we have 
	\begin{equation} \label{LB}
	\mathcal{J}_{A}(\alpha,\beta;k_1,k_2,K) \gtrsim \frac{|\log \alpha|}{v_{A}(k_{1},k_{2},K,r)}.
	\end{equation}
where $v_{A}(k_{1},k_{2},K,r)$ is given by Definition \ref{def_v_A}.
\end{theorem}

\begin{IEEEproof}
Appendix \ref{LBs}.
\end{IEEEproof}

We denote by $l_A$ the value of the parameter $l$ which corresponds to the maximum of the quantities in Definition \ref{def_v_A}, as it is used in the formulation of the following results.  

 \begin{definition}\label{def_l_A}			
  The quantity $l_A$ is defined as follows.
  
  \begin{itemize} 	
   \item If $v_{A}(k_{1},k_{2},K,r)$ is equal to \eqref{v1}, then  $l_A$ is the number such that
	 \begin{equation} \label{l_1}
	  v_{A}(k_{1},k_{2},K,r)= W(k_1-l_{A},k_2,K,\boldsymbol{I}(A),\boldsymbol{J_{l_A}}(A), r).
	 \end{equation}

	\item If $v_{A}(k_{1},k_{2},K,r)$  is equal to $\eqref{v2}$, then $l_A$ is the number such that  
	\begin{equation} \label{l_2}
	v_{A}(k_{1},k_{2},K,r)=W(k_1,k_2-l_A,K,\boldsymbol{I_{l_A}}(A),\boldsymbol{J}(A), r).
	\end{equation}

	\item If $v_{A}(k_{1},k_{2},K,r)$ is equal to $\eqref{v3}$, then $l_A=(|A^c|-k_2 +1)^+$.
 
	\item If $v_{A}(k_{1},k_{2},K,r)$ is equal to $\eqref{v4}$, then $l_A=(|A|-k_1 +1)^+$.
	\end{itemize}
\end{definition}
	
\begin{remark}
If $v_{A}(k_{1},k_{2},K,r)$ is equal to \eqref{v1}, the $l_{A}$ sources with the smallest KL numbers in $\boldsymbol{J}(A)$ are not considered in the evaluation of the difficulty of the testing problem. This is because we can intentionally misclassify these $l_A$ sources as anomalous without exceeding the tolerance level of $k_1$ false alarms in order to reduce the expected stopping time. The respective remark holds if $v_{A}(k_{1},k_{2},K,r)$ is equal to \eqref{v2}.

If $k_1 -1 \geq |A^c|-k_2 +1$ then we can intentionally misclassify as anomalous the $|A^c|-k_2 +1$ sources with the smallest KL numbers in $\boldsymbol{J}(A)$ without exceeding the tolerance level of $k_1$ false alarms in order to reduce the expected stopping time. The remaining $k_2 -1$ sources in $A^c$ are already less than the tolerance level of the $k_2$ missed detection errors and this is why the difficulty of the testing problem is determined only by the KL numbers in $\boldsymbol{I}(A)$ in \eqref{v3}. The respective remark holds if $k_2 -1 \geq |A|-k_1 +1$. 
\end{remark}	

\begin{corollary}
Suppose \eqref{momr} holds, and let $0<|A|<M$. As $\alpha,\, \beta \to 0$ so that \eqref{r} holds, we can distinguish the following cases.

\begin{itemize}
		\item If $v_{A}(k_{1},k_{2},K,r)$ is equal to \eqref{v1}, 
		\begin{equation*}
		\mathcal{J}_{A}(\alpha,\beta;k_1,k_2,K)  \gtrsim \frac{|\log \alpha|}{V(k_1 -l_A, K_1^{*},\boldsymbol{I}(A))} \sim \frac{|\log \beta|}{V(k_2, K_2^{*},\boldsymbol{J_{l_A}}(A))}.
		\end{equation*}
		
		\item If $v_{A}(k_{1},k_{2},K,r)$ is equal to \eqref{v2},
		\begin{equation*}
		\mathcal{J}_{A}(\alpha,\beta;k_1,k_2,K)  \gtrsim \frac{|\log \alpha|}{V(k_1, K_1^{*},\boldsymbol{I_{l_A}}(A))} \sim \frac{|\log \beta|}{V(k_2 -l_A, K_2^{*},\boldsymbol{J}(A))}.
		\end{equation*}

		\item If $v_{A}(k_{1},k_{2},K,r)$ is equal to \eqref{v3}, then			
		\begin{equation*}
		\mathcal{J}_{A}(\alpha,\beta;k_1,k_2,K)  \gtrsim \frac{|\log \alpha|}{V(k_1 -l_A, K,\boldsymbol{I}(A))}.
	    \end{equation*}
	    
	    \item If $v_{A}(k_{1},k_{2},K,r)$ is equal to \eqref{v4}, then
	    \begin{equation*}
	    \mathcal{J}_{A}(\alpha,\beta;k_1,k_2,K)  \gtrsim \frac{|\log \beta|}{V(k_2 -l_A, K,\boldsymbol{J}(A))}.
	    \end{equation*}
	\end{itemize}
In the first two cases $(K_1^{*},K_{2}^{*})$ is the maximizer of \eqref{rt} for the respective case.
\end{corollary}

\section{A criterion for asymptotic optimality}\label{sec:famil}
In this section, for each of the two problems  under consideration, we first introduce a  stopping and a decision rule so that  the corresponding error constraint is satisfied for any choice of sampling rule. For this, we  adopt  the approach in \cite{song_fel_gener}, where the  full sampling case was considered. Subsequently, we establish the second main result of this paper, which is a criterion on the sampling rule for the resulting policy to achieve the corresponding universal asymptotic lower bound in the previous section.

In what follows, for any sampling rule $R$ and  for each source  $i \in [M]$, we denote by  $\Lambda^R_{i}(n)$  the \textit{local} log-likelihood ratio (LLR) of source $i$ based on the  observations from it in the first $n$ time instants, i.e., 
\begin{align}\label{LLR}
\Lambda^R_{i}(n) &:= \sum_{m=1}^n   \log \left( \frac{f_{1i} (X_i(m)) }{ f_{0i}(X_i(m))} \right) \, R_i(m), \quad n \in \bN.
\end{align}

\subsection{The case of generalized misclassification error}
For any sampling rule $R$, we denote by  $T^{R}_{si}$ the first time that the sum of the $k$ smallest LLRs in absolute value is larger than some threshold $d>0$,  and by $\Delta_{si}^R$ the subset of data sources with positive LLRs upon stopping, i.e., 
\begin{align}
T^{R}_{si} &:= \inf\left\{ n\geq 1 \,:\, \sum_{i=1}^{k} \bar{\Lambda}^{R}_{i}(n) \geq d \right\},\label{st_r} \\
\Delta^{R}_{si}&:=\left\{ i \in [M] \,:\, \Lambda^R_{i}(T^{R}_{si}) > 0 \right\}, \label{st_d}
\end{align}
where, for each $i \in [M]$ and $n \in \bN$, $\bar{\Lambda}^{R}_{i}(n)$ denotes the \textit{$i^{th}$ smallest LLR in absolute value at time $n$}, i.e., the $i^{th}$ smallest element in $\{ |{\Lambda}^{R}_{j}(n)|\,:\, j \in [M] \}$.  In the full-sampling case, 
where $R(n)=[M]$ for every $n \in \bN$, $(T^{R}_{si},\Delta^{R}_{si})$ coincides with the \textit{sum-intersection rule}, introduced in \cite{song_fel_gener}. Similarly to  \cite[Theorem 3.1]{song_fel_gener}, it can be shown that if the threshold $d$ in \eqref{st_r} is selected as 
\begin{equation}\label{unk_thr}
d:=|\log \alpha|+\log {M \choose k},
\end{equation}
then the policy $(R,T^{R}_{si},\Delta^{R}_{si})$ satisfies the error constraint \eqref{error_const} for any sampling rule $R$. Next, we show that this policy, with this choice of threshold,  also achieves the asymptotic lower bound in  Theorem  \ref{th:lower_bound1} when, for each $i \in [M]$, the long-run frequency of the source that corresponds to $F_{i}(A)$ is not smaller than the quantity $c'_i(k,K,\boldsymbol{F}(A))$, defined  according to  \eqref{c_pr} of Lemma \ref{mm_g}. To be specific, we need the following definition.

\begin{definition}\label{def_c_star1}
For each $A \subseteq [M]$, we denote by 
$\boldsymbol{c}^{*}(A):=(c_1^*(A), \ldots, c_M^*(A))$ the 
permutation of $\boldsymbol{c}'(k,K,\boldsymbol{F}(A))$, defined in \eqref{c_pr}, so that
	\begin{align}\label{cv01}
	c^*_{(i)}(A):= c'_i(k,K,\boldsymbol{F}(A)), \quad i \in [M],
	\end{align}
	where $(i)$ denotes the source with the $i^{th}$ smallest number in the set  $\boldsymbol{F}(A)$, defined in \eqref{F(A)}. 
\end{definition}

We are now ready to state the first main result of this section. 

\begin{theorem}\label{asy_opt1}
	Consider a policy of the form $(R, T^R_{si}, \Delta^R_{si})$, where the threshold $d$ in \eqref{st_r} is selected according to \eqref{unk_thr} and the sampling constraint \eqref{samp_cons} is satisfied. Fix $A \subseteq [M]$, and suppose that, for all $i \in [M]$, the sampling rule $R$ satisfies 
    \begin{equation} \label{condi_misclas}	
	\sum_{n=1}^\infty \Pro_{A} \left(\pi_{i}^{R}(n)< c^*_{i}(A) -\epsilon \right)< \infty, \qquad \forall \; \epsilon>0,
	\end{equation}
	where $\boldsymbol{c}^{*}(A):=(c_1^*(A), \ldots, c_M^*(A))$ is defined according to Definition
\ref{def_c_star1}. Then, as $\alpha \to 0$, we have 
	\begin{equation}\label{AO1}
	\Exp_{A}\left[T^{R}_{si}\right] \sim \mathcal{J}_{A}(\alpha;k,K) \sim \frac{|\log \alpha |}{V(k,K,\boldsymbol{F}(A))},
	\end{equation}
	where $V$ is defined in \eqref{mm_st}.
\end{theorem}

\begin{IEEEproof}
	Appendix \ref{AO}.
\end{IEEEproof}

\subsection{The case of generalized familywise error metric}
For any sampling rule $R$, we denote by $p^{R}(n)$ the number of non-negative LLRs at time $n$, by $\hat{w}^{R}_1(n), \ldots, \hat{w}^{R}_{p^{R}(n)}(n)$ the indices of the increasingly ordered non-negative LLRs at time $n$, and by
$\check{w}^{R}_1(n), \ldots, \check{w}^{R}_{M-p^{R}(n)}(n)$ the indices of the decreasingly ordered negative LLRs at time $n$, i.e.,
\begin{align}
\begin{split}
0 &\leq \Lambda^{R}_{\hat{w}^{R}_{1}(n)}(n) \leq  \ldots \leq \Lambda^{R}_{\hat{w}^{R}_{p^{R}(n)}}(n), \\
0 &> \Lambda^{R}_{\check{w}^{R}_{1}(n)}(n) \geq \ldots \geq \Lambda^{R}_{\check{w}^{R}_{M-p^{R}(n)}}(n).
\end{split}
\end{align}
We also set
\begin{equation}
\begin{aligned}
\hat{\Lambda}^{R}_{i}(n)&:=
\begin{cases}
\Lambda^{R}_{\hat{w}^{R}_i(n)}(n), &\qquad i \leq p^{R}(n),\\
+\infty, &\qquad i> p^{R}(n), 
\end{cases}\\
\check{\Lambda}^{R}_{i}(n)&:=
\begin{cases}
-\Lambda^{R}_{\check{w}^{R}_i(n)}(n), &\quad i \leq M-p^{R}(n),\\
+\infty, &\quad i> M-p^{R}(n). 
\end{cases}
\end{aligned}
\end{equation}
For any integer $l$ such that $0 \leq l < k_{1}$, we set 
\begin{equation}\label{omh}
\hat{\tau}(l):= \inf \left\{n\geq 1 : \sum_{i=1}^{k_{1}-l}\hat{\Lambda}^{R}_{i}(n) \geq b,\,\, \sum_{i=1+l}^{k_{2}+l}\check{\Lambda}^{R}_{i}(n) \geq a\right\},
\end{equation}
Also, for any integer $l$ with $1\leq l < k_{2}$, we set
\begin{equation}\label{omc}
\check{\tau}(l):= \inf \left\{n \geq 1 : \sum_{i=1+l}^{k_{1}+l}\hat{\Lambda}^{R}_{i}(n) \geq b,\,\, \sum_{i=1}^{k_{2}-l}\check{\Lambda}^{R}_{i}(n) \geq a \right\}.
\end{equation}

We denote by $T^{R}_{leap}$ the minimum of the stopping times  in \eqref{omh} and \eqref{omc}, i.e.,
\begin{align}\label{tr_leap}
T^{R}_{leap}:=\min \left\{ \min\limits_{0\leq l <k_{1}}\hat{\tau}(l) , \min\limits_{1\leq l <k_{2}}\check{\tau}(l)  \right\},
\end{align}
and depending on whether the minimum is attained by a $\hat{\tau}(l)$ for some $l \in [0,k_1)$, or by a $\check{\tau}(l)$ for some $l \in [1, k_2)$, and we  set
\begin{align}\label{dr_leap} 
\Delta^{R}_{leap}:=
\begin{cases}
\{\hat{w}_{1}(\hat{\tau}(l)),..,\hat{w}_{p(\hat{\tau}(l))}(\hat{\tau}(l))\}\bigcup\{\check{w}_{1}(\hat{\tau}(l)),..,\check{w}_{l\wedge (M-p^{R}(n))}(\hat{\tau}(l))\}, & \mbox{if}
\quad  T^{R}_{leap}=\hat{\tau}(l),\\
\{\hat{w}_{l+1}(\check{\tau}(l)),..,\hat{w}_{p(\check{\tau}(l))}(\check{\tau}(l))\},  &  \mbox{if}  \quad T^{R}_{leap}=\check{\tau}(l). 
\end{cases} 
\end{align}
\\
In the full-sampling case, where $R(n)=[M]$ for every $n \in \bN$, $(T^{R}_{si},\Delta^{R}_{si})$ coincides with the \textit{leap rule},  introduced in \cite{song_fel_gener}. Similarly to \cite[Theorem 4.1]{song_fel_gener}, it can be shown that if the thresholds $a$ and $b$ in \eqref{omh}-\eqref{omc} are selected as
\begin{align}\label{ab_thrs}
\begin{split}
a&:=|\log \beta|+\log\left(2^{k_{2}} {M \choose k_{2}} \right), \\
b&:=|\log \alpha|+\log\left(2^{k_{1}} {M \choose k_{1}} \right),
\end{split}
\end{align}
then the policy $(R,T^{R}_{leap},\Delta^{R}_{leap})$ satisfies the error constraint \eqref{b_error_const} for any sampling rule $R$.  We next show that this policy, with this choice of threshold,  achieves the asymptotic lower bound in Theorem \ref{lower_bound_bk}, as long as the long-run sampling frequency of each source is sufficiently large.  To be specific, we introduce the following definition, for which we recall the definition of the quantity $l_A$ in Definition \ref{def_l_A}. 

\begin{definition} \label{def_c_star2}
For each $A \subseteq [M]$, we define the  vector  $\boldsymbol{c}^{*}(A):=(c_1^*(A), \ldots, c_M^*(A))$  as follows:
\begin{enumerate}[(i)]
  \item If $A=\emptyset$, then 
  $\boldsymbol{c}^{*}(A)$ 
is a permutation of the vector $\boldsymbol{c}'(k_2,K,\boldsymbol{J_{k_1 -1}}(A))$, defined in \eqref{c_pr}, such that 
  \begin{equation}
   c_{\{i\}}^*(A):=c'_{i}(k_2,K,\boldsymbol{J_{k_1 -1}}(A)), \quad i \in [M-k_1 +1],
  \end{equation}
  where   $\{i\}$ denotes  the source with the $i^{th}$ smallest element in $\boldsymbol{J_{k_1 -1}}(A)$.
       
  \item If  $A=[M]$, then 
  $\boldsymbol{c}^{*}(A)$ 
is a permutation of  the vector
   $\boldsymbol{c}'(k_1,K,\boldsymbol{I_{k_2 -1}}(A))$, defined in \eqref{c_pr}, such that 
  \begin{equation}
  c_{<i>}^*(A):=c_{i}'(k_1,K,\boldsymbol{I_{k_2 -1}}(A)), \quad i \in [M-k_2 +1],
  \end{equation}
where    $<i>$ denotes the 
 source with the $i^{th}$ smallest element in $\boldsymbol{I_{k_2 -1}}(A)$.
  
  \item If $0<|A|<M$ and $v_{A}(k_1, k_2, K,r)$ is equal to \eqref{v3}, then $\boldsymbol{c}^{*}(A)$ 
is a permutation of  the vector
$\boldsymbol{c}'(k_1 -l_A,K,\boldsymbol{I}(A)),$ defined in \eqref{c_pr},   such that 
  \begin{equation}
  c_{<i>}^*(A):=c_{i}'(k_1 -l_A, K,\boldsymbol{I}(A)), \quad i \in [|A|],
  \end{equation}
where  $<i>$ denotes the  source in $A$ with the $i^{th}$ smallest element in $\boldsymbol{I}(A)$.
  
  \item If  $0<|A|<M$ and $v_{A}(k_1, k_2, K,r)$ is equal to \eqref{v4},      then $\boldsymbol{c}^{*}(A)$ 
is  a permutation of  the vector
 $\boldsymbol{c}'(k_2 -l_A,K,\boldsymbol{J}(A)),$ defined according to  \eqref{c_pr}, such that
  \begin{equation}
  c_{\{i\}}^*(A):=c'_{i}(k_2 -l_A,K,\boldsymbol{J}(A)), \quad i \in [|A^c|],
  \end{equation}
  where  $\{i\}$  denotes the source  in $A^c$ with the $i^{th}$ smallest element in $\boldsymbol{J}(A)$.
  
  \item If  $0<|A|<M$ and  $v_{A}(k_1, k_2, K,r)$ is equal to \eqref{v1},     then $\boldsymbol{c}^{*}(A)$  is  a permutation of  the vector
$\boldsymbol{c}'(k_1 -l_{A},k_2,K,\boldsymbol{I}(A),\boldsymbol{J_{l_A}}(A))$, defined according to \eqref{c_pr2},  such that
  \begin{align}
   c^{*}_{<i>}(A) &:= \hat{c}'_{i}(k_1 -l_{A},k_2,K,\boldsymbol{I}(A),\boldsymbol{J_{l_A}}(A)),  \quad i \in [|A|],     \label{cv02} \\
   c^{*}_{\{j\}}(A) &:= \check{c}'_{j}(k_1 -l_{A},k_2,K,\boldsymbol{I}(A),\boldsymbol{J_{l_A}}(A)), \quad j \in [|A^c|-l_{A}], \label{cv03}
  \end{align}  
where   $<i>$ is the source  in $A$ with the $i^{th}$ smallest element in $\boldsymbol{I}(A)$,  and $\{j\}$ the source in $A^c$ with the  $j^{th}$ smallest element in $\boldsymbol{J_{l_A}}(A)$.
  
  \item If $0<|A|<M$ and $v_{A}(k_1, k_2, K,r)$ is equal to \eqref{v2},
then $\boldsymbol{c}^{*}(A)$  is  a permutation of  the vector
$\boldsymbol{c}'(k_1,k_2-l_{A},K,\boldsymbol{I_{l_{A}}}(A),\boldsymbol{J}(A))$, defined according to  \eqref{c_pr2}, such that
  \begin{align}
  c^{*}_{<i>}(A) &:= \hat{c}'_{i}(k_1,k_2-l_{A},K,\boldsymbol{I_{l_A}}(A),\boldsymbol{J}(A)),  \quad i \in [|A|-l_A],     \label{cv04} \\
  c^{*}_{\{j\}}(A) &:= \check{c}'_{j}(k_1,k_2-l_{A},K,\boldsymbol{I_{l_A}}(A),\boldsymbol{J}(A)), \quad j \in [|A^c|], \label{cv05}
  \end{align}  
 where  $<i>$ denotes the source in $A$ with the $i^{th}$ smallest element in $\boldsymbol{I_{l_{A}}}(A)$, and 
 $\{j\}$   the source in $A^c$ with the $j^{th}$ smallest element in $\boldsymbol{J}(A)$.
\end{enumerate}
\end{definition}

We are now ready to state the second main result of this section. 

\begin{theorem}\label{asy_opt}
Consider a policy of the form $(R,T^{R}_{leap},\Delta^{R}_{leap})$,	where the thresholds $a$ and $b$  in \eqref{omh}-\eqref{omc} are selected according to \eqref{ab_thrs} and the sampling constraint \eqref{samp_cons} is satisfied. Fix $A \subseteq [M]$,  and suppose that, for all $i \in [M]$, the sampling rule $R$ satisfies
\begin{equation} \label{condi_generalized}
\sum_{n=1}^\infty \Pro_{A} \left(\pi_i^{R}(n)< c^*_i(A)- \epsilon \right)< \infty, \quad \forall \; \epsilon > 0,
\end{equation} 
where  the vector $\boldsymbol{c}^{*}(A):=(c_1^*(A), \ldots, c_M^*(A))$
is defined according to Definition \ref{def_c_star2}.  
\begin{itemize}
\item If $A=\emptyset$, then, for any $\alpha \in (0,1)$,
as $\beta \to 0$ we have
  \begin{align} \label{AO21}
  \begin{split}
	\Exp_{A}\left[ T^{R}_{leap} \right] &\sim \mathcal{J}_{A}(\alpha,\beta; k_1,k_2,K) \sim 
	\frac{|\log \beta|}{V(k_2,K, \boldsymbol{J_{k_1 -1}}(A))}.
	\end{split}
	\end{align} 
	
	\item If $A=[M]$, then, for any $\beta \in (0,1)$,
as $\alpha \to 0$ we have
	  \begin{align} \label{AO22}
  \begin{split}
	\Exp_{A}\left[ T^{R}_{leap} \right] &\sim \mathcal{J}_{A}(\alpha,\beta; k_1,k_2,K)
	\sim  
	\frac{|\log \alpha|}{V(k_1,K, \boldsymbol{I_{k_2 -1}}(A))}.
	\end{split}
	\end{align} 

 \item If $0 < |A|< M$, then,
as $\alpha, \beta \to 0$ so that \eqref{r} holds we have
  \begin{align} \label{AO23}
  \begin{split}
	\Exp_{A}\left[ T^{R}_{leap} \right] &\sim \mathcal{J}_{A}(\alpha,\beta; k_1,k_2,K) \sim  
	\frac{|\log \alpha |}{v_{A}(k_1, k_2, K,r)},
	\end{split}
	\end{align} 
	\end{itemize}
where $V$ is defined in \eqref{mm_st}, and $v_A$ in Definition \ref{def_v_A}.
\end{theorem}

\begin{IEEEproof}
	Appendix \ref{AO}.
\end{IEEEproof}

\section{Asymptotically optimal probabilistic sampling rules}\label{sec:ao}
In this section, we design sampling rules that satisfy the criteria for asymptotic optimality  established in Section \ref{sec:famil} \text{simultaneously for every possible subset of anomalies}. For this, we first introduce a notion of consistency, which applies to   an arbitrary sampling rule. Then, we define  a family of \textit{probabilistic} sampling rules, and finally, we show how to design a probabilistic sampling rule  so that  condition \eqref{condi_misclas}  (resp. \eqref{condi_generalized}) is satisfied for every $i \in [M]$ and  $A \subseteq [M]$, and subsequently so that  the first-order asymptotic optimality property \eqref{AO1}  (resp.  \eqref{AO21}-\eqref{AO23}) holds for every $A \subseteq [M]$.


\subsection{Consistency}
We say that a sampling rule $R$ is \textit{consistent} if the subset of sources with non-negative LLRs at time $n$, i.e., 
\begin{equation}\label{dfr}
\mathfrak{D}_n^{R} := \{ i \in [M] \,:\, \Lambda^R_{i}(n) \geq 0\},
\end{equation}
converges \textit{quickly}  to the true subset of anomalous sources $A$.

\begin{definition}   
    Fix $A \subseteq [M]$.  We say that a sampling rule $R$ is \textit{consistent under} $\Pro_A$ if	
	\begin{equation}\label{es}
	\Exp_{A} \left[\sigma^R_{A} \right] < \infty,
	\end{equation}
	where $\sigma_A^R$ is the random time starting from which the sources in $A$ are the only ones with non-negative LLR, i.e.,
	\begin{equation}\label{sigmaA}
	\sigma^R_{A} :=\inf\left\{n \in \bN: \mathfrak{D}^{R}_m =A  \quad \text{for all} \; \; m \geq n\right\}.
	\end{equation}  
	\end{definition}
	
	\begin{definition}   
	We say that a sampling rule $R$  is \textit{consistent} if it is consistent under $\Pro_A$, for every $A \subseteq [M]$.
\end{definition}

From \cite[Theorem 3.1]{a_prob} we know that if there is a $\rho > 0$ such that, for each $i \in [M]$,  the sequence $\{ \Pro_{A}\left(\pi^{R}_{i}(n) < \rho \right) \, :\, n \in \bN\}$ is exponentially decaying, then the sequence $\{\Pro_{A}\left(\sigma^R_{A} > n\right) \, :\, n \in \bN \}$ is exponentially decaying, and, as a result, the sampling rule $R$ is consistent under $\Pro_A$. Next, we state a less restrictive criterion for consistency, according to which a sampling rule is consistent even if the long-run sampling frequency of all sources is equal to $0$, as long as the decay is not fast. 

\begin{theorem}\label{cons}
Suppose that condition \eqref{momr0} holds for some $\mathfrak{p}>4$, and let $\delta \in \left(0, \frac{1}{2}-\frac{2}{\mathfrak{p}}\right)$, and $C>0$. Fix $A \subseteq [M]$.  If   $R$ is an arbitrary sampling rule for which
\begin{equation}\label{spar}
\sum_{n=1}^{\infty} n \, \Pro_{A}\left(\pi^{R}_{i}(n) < C \, n^{-\delta}\right) < \infty, \quad \forall \, i\, \in [M],
\end{equation} 
then $R$ is consistent under $\Pro_{A}$.
\end{theorem}

\begin{IEEEproof}
Appendix \ref{AO}.
\end{IEEEproof}

\subsection{Probabilistic sampling rules}
We say that a sampling rule  $R$ is \textit{probabilistic} if there exists a function
$$q^R : 2^{[M]} \times \bN \times 2^{[M]} \to [0,1]$$
such that, for every $n  \in \bN$, $D \subseteq [M]$, and $B \subseteq [M]$, $q^R\left(B;n,D\right)$ is the probability that $B$ is the subset of sampled sources at time $n$ when $D$ is the subset of sources with non-negative LLRs at time $n-1$, i.e.,
\begin{align}\label{q}
q^R\left(B;n,D\right) := \Pro \left( R(n)=B \, | \, \cF^R_{n-1} , \mathfrak{D}^{R}_{n-1}=D\right).
\end{align} 
For such a sampling rule,  for each source $i \in [M]$ we denote by $c^R_{i}\left(n,D\right)$ the probability with which source $i$ is sampled at time $n$ when $D$ is the subset of sources with non-negative LLRs at time $n-1$, i.e.,
\begin{align} \label{c}
c^R_{i} \left(n,D \right) &:=\Pro \left( R_i(n)=1 \, |  \, \cF^R_{n-1}, \mathfrak{D}^{R}_{n-1}=D \right) 
\end{align}
thus, 
\begin{align} \label{chr}
c^R_{i} \left(n,D \right) = \sum_{B \subseteq [M]: \, i \in B}  q^R\left(B;n,D \right).
\end{align}

In the following theorem, we state a condition under which a \textit{consistent} probabilistic sampling rule satisfies condition \eqref{condi_misclas} or \eqref{condi_generalized} simultaneously for every $A \subseteq [M]$. 

\begin{theorem}\label{c_st}
	Let $R$ be a consistent probabilistic sampling rule, and  fix $A \subseteq [M]$.
	\begin{itemize}
	\item[(i)]  If for all $i \in [M]$ 
	we have 
	\begin{equation}\label{cr_geq}
	\liminf_{n \to \infty} c_i^{R}(n,A) \geq c^*_{i}(A),
	\end{equation}  	
	where $(c^*_{1}(A), \ldots, c^*_{M}(A))$ is given by Definition \ref{def_c_star1}  (resp. Definition \ref{def_c_star2}), then condition   \eqref{condi_misclas} (resp. \eqref{condi_generalized}) is satisfied for every $i \in [M]$.
	\item[(ii)]   If, also, the sampling constraint \eqref{samp_cons} is satisfied, then the first-order asymptotic optimality property \eqref{AO1} (resp. \eqref{AO21}-\eqref{AO23}) holds.
	\end{itemize}
\end{theorem}

\begin{IEEEproof}
Part (i) is proven in Appendix \ref{AO}, and part (ii) follows by Theorem \ref{asy_opt1} (resp. Theorem \ref{asy_opt}).
\end{IEEEproof}

Condition \eqref{cr_geq} is clearly satisfied if source $i$ is sampled at each instant with probability $c^*_{i}(D)$ when the estimated anomalous subset at the previous time instant is $D$, i.e.,
\begin{equation}\label{cr_eq} 
c_i^{R}(n, D) = c^*_{i}(D), \quad \forall \; n \in \bN, \; D \subseteq [M], \; i \in [M]. 
\end{equation}
From \cite[Theorem 4.1]{a_prob} it follows that 
this choice  implies that the sampling rule is consistent 
 when 
\begin{equation}\label{tt} 
c^*_{i}(D)>0, \quad \forall \;  D \subseteq [M], \; i \in [M]. 
\end{equation}

In case condition \eqref{tt} does not hold, selection  \eqref{cr_eq} no longer guarantees the consistency of the sampling rule, and needs to be modified. Indeed, for a  source  $i \in [M]$ for which $c^{*}_{i}(D)=0$, then the sampling frequency of source $i$ should converge to $0$ slowly enough when the true  subset of anomalous sources is $D$, so that Theorem 
\ref{c_st} is applicable.  To be more specific,  let $\{b_n : n \in \bN\}$ be a sequence of positive reals that converges to 0, which we will specify later, and for each $n \in \bN$, $i \in [M]$, and $ D \subseteq [M]$  set
\begin{equation}\label{cr_pr}
c_i^{R}(n, D) =
\begin{cases}
b_n , \quad & \text{if} \quad c^*_{i}(D)=0,\\
c^*_{i}(D)- ( \tilde{l}_D/ (M- \tilde{l}_D) )\,  b_n, \quad   & \text{if} \quad c^*_{i}(D)>0,
\end{cases}
\end{equation}
where  $\tilde{l}_D$ is the number of zero entries in the vector $\mathbf{c}^{*}(D)$.

As we show in the following proposition, a suitable selection of the sequence $\{b_n : n \in \bN\}$ guarantees the consistency of the sampling rule. 

\begin{proposition}\label{cst_rp}
\begin{itemize}
\item[(i)]
If  \eqref{tt} holds,  and the  probabilistic sampling rule $R$ satisfies \eqref{cr_eq}, then condition \eqref{condi_misclas}  (resp. \eqref{condi_generalized}) is satisfied for every $i \in [M]$.

\item[(ii)] Suppose condition \eqref{momr0} holds for some $\mathfrak{p}>4$. If $R$ is a probabilistic sampling rule that satisfies  \eqref{cr_pr} for every $i \in [M], D \subseteq [M]$, $n \in \bN$, 
where  $b_n = C_p \, n^{-\delta}$ for some  $\delta \in \left(0,\frac{1}{2}-\frac{2}{\mathfrak{p}}\right)$ and some  $C_{p}>0$  small  enough such that 
\begin{equation}\label{c_gd}
c_i^{R}(n, D) \geq C_p \, n^{-\delta}
\end{equation}
holds for all $n \in \bN$,  $D \subseteq [M]$, and  $i \in [M]$, then condition \eqref{condi_misclas}  (resp. \eqref{condi_generalized}) is satisfied for every $i \in [M]$.
\end{itemize}
\end{proposition}	

\begin{IEEEproof}
Since  \eqref{cr_eq} and  \eqref{cr_pr} both   imply  \eqref{cr_geq} for any $A \subseteq [M]$, by Theorem \ref{c_st}  it follows that it suffices to show that $R$ is consistent. As discussed earlier,  for (i)  this follows from  \cite[Theorem 4.1]{a_prob}. For (ii),  the proof of consistency is presented in Appendix \ref{AO}.	\\
\end{IEEEproof}

The previous proposition provides concrete selections of $c^R(n,D)$ that guarantee the consistency of the sampling rule and conditions  \eqref{condi_misclas} or
  \eqref{condi_generalized}. To achieve the corresponding asymptotic optimality property, the sampling rule should satisfy the sampling constraint  \eqref{samp_cons}. This  is clearly the case if at most $ \lfloor K \rfloor$ sources are sampled at each instant, i.e., 
\begin{align}\label{chernoff0}
q^R\left(B; n, D \right) =0 \quad   \text{for all}  \quad n \in \bN, \; D \subseteq [M], \quad \text{and }\; B \subseteq [M] \; \text{ such that } \;   |B| >  \lfloor K \rfloor. 
\end{align}
Combining this  observation with the previous proposition we can now state the theorem that summarizes the main result of this section.  

\begin{theorem}\label{c_st_cor}
Consider an integer $K$.
\begin{itemize}
\item[(i)] Suppose \eqref{tt} holds. If $R$ is a probabilistic sampling rule that satisfies \eqref{cr_eq} and \eqref{chernoff0}, then the first-order asymptotic optimality property \eqref{AO1} (resp. \eqref{AO21}-\eqref{AO23}) holds for every $A \subseteq [M]$.
\item[(ii)] If $R$ is a probabilistic sampling rule that satisfies \eqref{cr_pr} and \eqref{chernoff0} for every $i \in [M], D \subseteq [M]$, $n \in \bN$, 
where  $b_n = C_p \, n^{-\delta}$ for some  $\delta \in \left(0,\frac{1}{2}-\frac{2}{\mathfrak{p}}\right)$ and some  $C_{p}>0$  small  enough such that  \eqref{c_gd} 
holds for all $n \in \bN$,  $D \subseteq [M]$, and  $i \in [M]$,  then the first-order asymptotic optimality property \eqref{AO1} (resp. \eqref{AO21}-\eqref{AO23}) holds for every $A \subseteq [M]$.
\end{itemize}
\end{theorem}	
\begin{IEEEproof}
The claim follows by Proposition \ref{cst_rp}, and Theorems \ref{asy_opt1}, \ref{asy_opt}, respectively.
\end{IEEEproof}

\begin{remark}
Theorem \ref{c_st_cor} remains valid even if $K$ is not an integer as long as 
\begin{equation*}
\sum_{i=1}^{M} c^{*}_{i}(D) \leq \lfloor K \rfloor, \quad \forall \, D \subseteq [M].
\end{equation*}	
\end{remark}

\begin{remark}
Proposition \ref{cst_rp} and Theorem \ref{c_st_cor} remain valid even if $c_i^{R}(n, D)$ are chosen greater than or equal to \eqref{cr_eq} (resp. \eqref{cr_pr}), as long as the sampling constraint \eqref{chernoff0} is satisfied.	
\end{remark}

\section{Simulation study}\label{sec:sim}
In this section, we present the results of various simulations by which we illustrate the asymptotic theory that was developed in the previous sections. In all simulations, there are $M=10$ sources and, for each source $i \in [M]$, the observations are normally distributed with variance $1$ and mean $0$ if source $i$ is not anomalous, whereas the mean is $\mu_i$ if it is anomalous, i.e., $f_{0i}=\mathcal{N}(0,1)$ and  $f_{1i}=\mathcal{N}(\mu_i,1)$, and thus  $I_{i}=J_{i}=(\mu_{i})^2/2$. For our simulations, we consider a non-homogeneous setup, where 
\begin{equation*}
\mu_{i}=
\begin{cases}
0.5, \quad &1\leq i \leq 3,\\
0.7, \quad &4\leq i \leq 7,\\
1, \quad   &8 \leq i \leq 10.
\end{cases}
\end{equation*} 
We apply the probabilistic sampling rule (see Section \ref{sec:ao}), which observes the values of exactly $K=5$ sources per sampling instant. When controlling the generalized misclassification error, we apply the sum-intersection rule \eqref{st_r}-\eqref{st_d}, whereas when  controlling the generalized familywise errors, we apply the leap rule \eqref{tr_leap}-\eqref{dr_leap}. For the computation of each expected time for stopping, we apply $10^4$ simulation runs, and the standard error in each expected time for stopping is $1$, whereas the standard error for each ratio of expected times for stopping is $10^{-2}$, in all cases.

In all simulations, we fix the true, but unknown, anomalous set to be $A=\{1, \ldots, 5\}$. The numbers $\{F_{i}(A) \,:\, i \in [M]\}$ defined in Subsection \ref{mmiscl} for $I_{i}=J_{i}=(\mu_{i})^2/2$, are equal to
\begin{align*}
F_{i}(A) &=
\begin{cases}
0.125, \quad & 1 \leq i \leq 3, \\
0.245, \quad & 4 \leq i \leq 7, \\
0.5, \quad & 8 \leq i \leq 10,
\end{cases}
\end{align*}
and for $\{I_{i}(A) \,:\, i \in [|A|]\}$, $\{J_{i}(A) \,:\, i \in [|A^c|]\}$ defined in Subsection \ref{gener_fam}, we have
\begin{align*}
I_{i}(A)&=F_{i}(A), \quad 1 \leq i \leq  5, \\
J_{i-5}(A)&=F_{i}(A), \quad  6 \leq i \leq 10.
\end{align*}

\subsection{Controlling the generalized misclassification error}

\begin{figure}[htb!]
	\centering
	\centerline{\includegraphics[width=8.5cm]{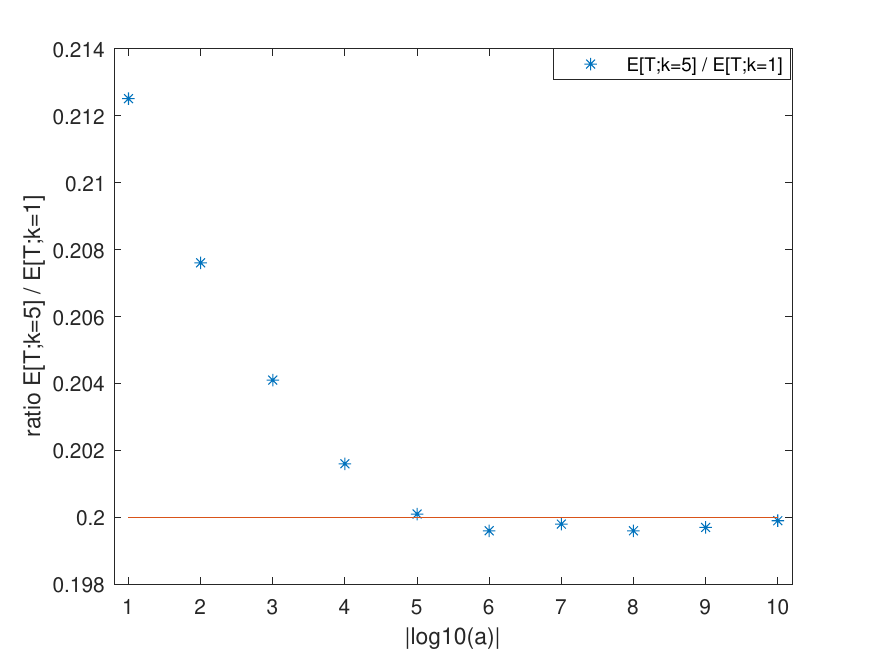}}
	\caption{Ratio of expected stopping times $\Exp[T;k=5] / \Exp[T;k=1]$ versus $|\log_{10}(\alpha)|$.}
	\label{fig3}
\end{figure}

In the case of generalized \textit{misclassification} error rate, by Theorem \ref{asy_opt1} as  $\alpha \to 0$
\begin{equation}
\mathcal{J}_{A}(\alpha;k,K) \sim \frac{|\log \alpha |}{V(k,K,\boldsymbol{F}(A))},
\end{equation}
where for $K=5$ we have
\begin{equation}\label{kkf}
V(k,K,\boldsymbol{F}(A))=
\begin{cases}
k \, \frac{K}{M} \, \widetilde{F}_{1}(A), \quad &\mbox{if }\; k \in \{1,\ldots,5\},\\
(k-3)\, \frac{K}{M-3} \, \widetilde{F}_{4}(A), \quad &\mbox{if }\; k \in \{6,\ldots,8\},\\
\sum_{i=6}^{9} F_{i}(A), \quad &\mbox{if }\; k=9,\\
\sum_{i=6}^{10} F_{i}(A), \quad &\mbox{if }\; k=10.
\end{cases}
\end{equation}
where for the computation of $V(k,K,\boldsymbol{F}(A))$ we used Algorithm \ref{alg}. We can easily verify that $k(K/M)\widetilde{F}_{1}(A)$ is less than or equal to the respective value of $V(k,K,\boldsymbol{F}(A))$ for all $k \in [M]$, and as a result
\begin{equation}\label{ast}
\frac{\mathcal{J}_{A}(\alpha;k,K)}{ \mathcal{J}_{A}(\alpha;1,K)} \lesssim \frac{1}{k}, \quad \forall \, k \in [M].
\end{equation}
In Figure \ref{fig3}, we plot the ratio of the expected time for stopping for $k=5$ over that for $k=1$, against $|\log_{10}(\alpha)|$ for $\alpha \in \{  10^{-1}, \ldots, 10^{-10} \}$. For each value of $\alpha$, the thresholds are selected according to \eqref{unk_thr}. As expected by the form of $V(k,K,\boldsymbol{F}(A))$ for $k=5$ in \eqref{kkf}, the ratio converges to $1/5$, and this is also depicted in Figure \ref{fig3}, where by $\Exp[T;k=1]$, $\Exp[T;k=5]$ we denote the expected stopping time for $k=1$, $k=5$, respectively.

\begin{figure}[htb!]
	\centering
	\centerline{\includegraphics[width=8.5cm]{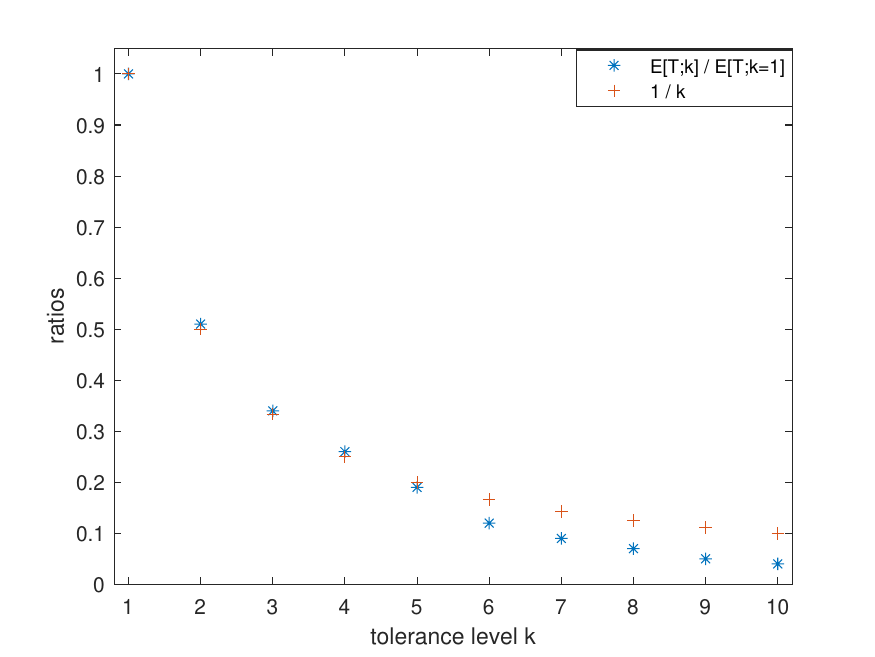}}
	\caption{Ratio of expected stopping times $\Exp[T;k] / \Exp[T;k=1]$ versus the tolerance level k.}
	\label{fig1}
\end{figure}

In order to verify the limit \eqref{ast} as an approximation in the finite regime, we fix $\alpha=10^{-3}$ and select the thresholds of the sum-intersection rule in \eqref{st_r}, via Monte Carlo simulation, so that the probability of at least $k$ errors, of any kind, is equal to $10^{-3}$. In Figure \ref{fig1}, we plot the ratio of the expected time for stopping for $k$, over that for $k=1$, against $k \in [M]$. In Figure \ref{fig1}, we can see that the expected time for stopping when $k=1$ reduces by a factor approximately equal to $1/k$ for $1 \leq k \leq 5$, and clearly less than $1/k$ for $5 \leq k \leq 10$. In Figure \ref{fig1}, we denote by $\Exp[T;k]$ the expected stopping time for the respective value of $k$.

\subsection{Controlling the generalized familywise errors}

In the case of generalized \textit{familywise} error rates, for $r=1$ and since $0<|A|<M$, by Theorem \ref{asy_opt} as $\alpha, \, \beta \to 0$
\begin{equation}
\mathcal{J}_A(\alpha, \beta; k_1, k_2, K) \sim \frac{|\log(\alpha)|}{v_{A}(k_{1},k_{2},K,r)},
\end{equation}
where for $K=5$, $v_{A}(k_{1},k_{2},K,r)$ has the following form
\begin{equation}\label{vak}
v_{A}(k_{1},k_{2},K,r)=
\begin{cases}
k_1 \, y_1\, \widetilde{I}_{1}(A) = k_2 \, y_2 \, \widetilde{J}_{1}(A), \quad & k_1 = k_2 \in \{1,2\}, \\
k_1 \, y_1\, \widetilde{I}_{2}(A) = (k_2 -1)\, y_2 \, \widetilde{J}_{1}(A), \quad & k_1 = k_2=3,\\
x_1 \, I_{3}(A) + \sum_{i=4}^{5} I_{i}(A) = (k_2 -1)\, y_2 \, \widetilde{J}_{1}(A), \quad & k_1 = k_2=4,\\
x_1 \, I_{2}(A) + \sum_{i=3}^{5} I_{i}(A) = x_2 \, J_{4}(A) + J_{5}(A), \quad & k_1 = k_2=5.
\end{cases}
\end{equation}
and $x_1$, $x_2$, $y_1$, $y_2$, and $l_{A}$ are given in Table \ref{ttb}.

\begin{table}[h!]
	\begin{center}
		\begin{tabular}{ |c||c|c|c|c|c| } 
			\hline
			$k_1=k_2 $ & 1 & 2 & 3 & 4 & 5\\ 
			\hline
			\hline
			$l_A$ & 0 & 0 & 1 & 1 & 0\\ 
			\hline
			$x_{1}$ & - & - & - & 0.61 & 0.61 \\ 
			\hline
			$x_{2}$ & - & - & - & - & 0.39 \\ 
			\hline
			$y_{1}$ & 0.69 & 0.69 & 0.66 & - & - \\ 
			\hline
			$y_{2}$ & 0.30 & 0.30 & 0.46 & 0.53 & - \\ 
			\hline
		\end{tabular}
	\end{center}
	\caption{The $l_A$, $x_{1}$, $x_{2}$, $y_{1}$, $y_{2}$ for each $k_1=k_2 \in [M/2]$.}
	\label{ttb}
\end{table}
For the computation of $v_{A}(k_{1},k_{2},K,r)$ we first use Algorithm \ref{alg2}, and then for the computation of $x_{1}$, $x_{2}$, $y_{1}$, $y_{2}$ we used Algorithm \ref{alg}. We can easily verify that $k_1 \, y_1\, \widetilde{I}_{1}(A)$ is less than or equal to $v_{A}(k_{1},k_{2},K,r)$ for all $k_1 = k_2 \in [M/2]$, and as a result
\begin{equation}\label{ast2}
\frac{\mathcal{J}_A(\alpha, \beta; k_1, k_1, K)}{ \mathcal{J}_A(\alpha, \beta; 1, 1, K)} \lesssim \frac{1}{k_1}. 
\end{equation}

\begin{figure}[htb!]
	\centering
	\centerline{\includegraphics[width=8.5cm]{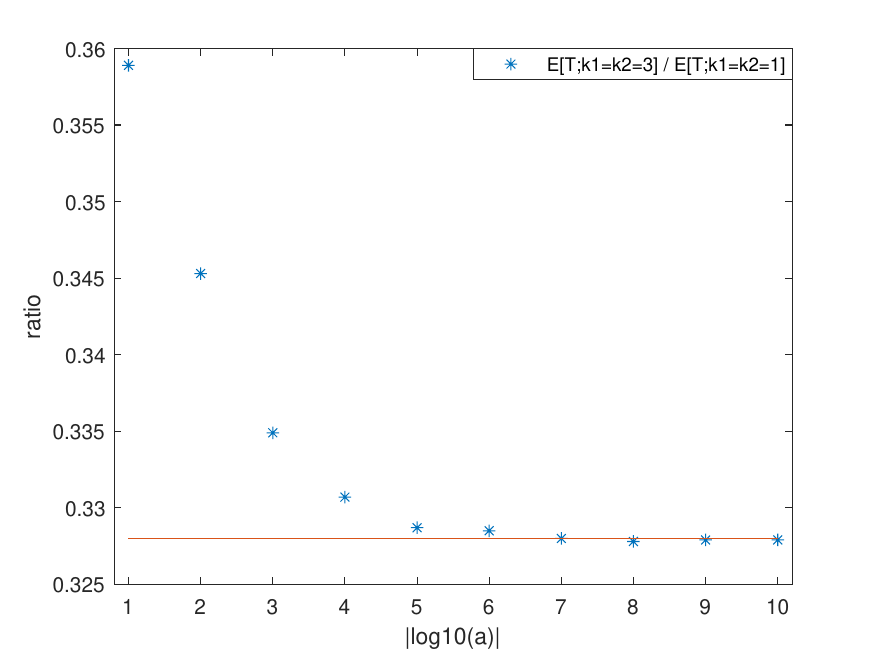}}
	\caption{Ratio of expected stopping times $\Exp[T;k_1=k_2=3] / \Exp[T;k_1=k_2=1]$ versus $|\log_{10}(\alpha)|$.}
	\label{fig4}
\end{figure}

In Figure \ref{fig4}, we plot the ratio of the expected stopping time for $k_1=k_2=3$ over that for $k_1=k_2=1$, against $|\log_{10}(\alpha)|$ for  $\alpha \in \{  10^{-1}, \ldots, 10^{-10}\}$. For each value of $\alpha$, and $\beta=\alpha$, the thresholds are selected according to \eqref{ab_thrs}. As expected by the form of $v_A$ in \eqref{vak}, it follows that this ratio converges to $0.328$, which is also depicted in Figure \ref{fig4}, where by $\Exp[T;k_1=k_2=3]$, $\Exp[T;k_1=k_2=1]$ we denote the expected stopping time for $k_1=k_2=1$, $k_1=k_2=3$, respectively.

\begin{figure}[htb!]
	\centering
	\centerline{\includegraphics[width=8.5cm]{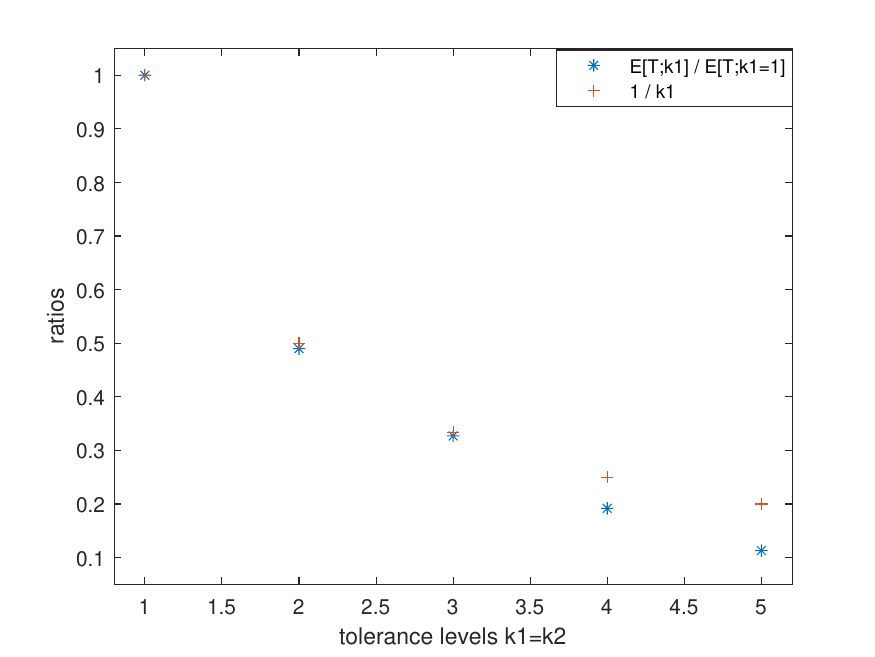}}
	\caption{Ratio of expected stopping times $\Exp[T;k_1] / \Exp[T;k_1=1]$ versus the tolerance level $k_{1}=k_{2}$.}
	\label{fig2}
\end{figure}

In order to verify the limit \eqref{ast2} as an approximation in the finite regime, we fix $\alpha=\beta=10^{-3}$ and we select the thresholds of the leap rule \eqref{tr_leap}, via Monte Carlo sumulation, so that the probabilities of at least $k_1$ false positives and at least $k_2$ false negatives are both equal to $\alpha=\beta=10^{-3}$. In Figure \ref{fig2}, we depict the ratio of the expected time for stopping with $k_1 = k_2$, over that for $k_1 = k_2 = 1$, against $k_1 = k_2 \in [M/2]$. In Figure \ref{fig2}, we observe that the expected stopping time for $k_1=k_2=1$, reduces by a factor of approximately $1/k_1$, as $k_1$ increases. 

\section{Conclusion}\label{sec:con}
In this work,  we study the sequential anomaly identification problem in the presence of a sampling constraint under two different formulations that involve  generalized error metrics. For each of them we establish  a universal asymptotic lower bound, and we show that  it is attained by a policy that combines the stopping and decision rule  that was proposed in the case of full-sampling in \cite{song_fel_gener} and a probabilistic sampling rule which is designed to achieve specific long-run sampling frequencies. 
The optimal performance is characterized, and the impact of the sampling constraint and tolerance to errors is  assessed, 
both to a first-order asymptotic approximation as the error probabilities go to 0. 
These theoretical asymptotic results are also illustrated via simulation studies.

  
Directions for further research involve (i) the incorporation of prior information on the number of anomalous sources, as  in \cite{a_prob} in the case of classical familywise error control $(k_1 =k_2 = 1)$, (ii) consideration of composite hypotheses for the testing problem in each source, as in \cite{song_fel_gener} in the full sampling case, and in \cite{hemo_g} for the case we know a priori that there is only one anomalous source and we can observe one source at a time, (iii) varying sampling or switching cost per source, as in \cite{gur19} and \cite{lambez2021anomaly}, respectively. Other directions include the framework where the acquired observations are not conditionally independent of the past \cite{veer15}, as well as the consideration of a dependence structure within the observations from different sources \cite{hey16}.

\section*{Acknowledgments}
This work was supported by the University of Illinois at Urbana–Champaign research support award RB21036.

\bibliographystyle{IEEEtran}
\bibliography{biblio_g}

\appendices
\renewcommand{\thetheorem}{ \Alph{section}.\arabic{theorem}}
\renewcommand{\thelemma}{ \Alph{section}.\arabic{lemma}}

\section{}\label{miscl}
In Appendix \ref{miscl}, we prove Lemma \ref{mm_g} and Lemma \ref{equiv_opt}, which provide the solutions to the max-min problems \eqref{mm_st} and \eqref{mm_st2}, respectively. We also provide Algorithm \ref{alg} for the computation of $x,y \in [0,1)$ and $u,v \in [0,\kappa]$, and Algorithm \ref{alg2} for the computation of $(K^{*}_{1},K^{*}_{2})$.

\begin{IEEEproof}[Proof of Lemma \ref{mm_g}]
Since $\cD(K)$ is compact, and $\{U \subseteq [|\boldsymbol{L}|]\, : \; |U|=\kappa \}$ is finite, the max-min problem \eqref{mm_st} has a solution. The max-min structure of \eqref{mm_st} implies that a maximizer $(\tilde{c}_1,\ldots,\tilde{c}_{|\boldsymbol{L}|})$ must satisfy the following two conditions:
\begin{enumerate}[(i)]
	\item for each $i \in \{1,\ldots,\kappa\}$, and each $j \in \{\kappa+1,\ldots,|\boldsymbol{L}|\}$, it holds
          \begin{equation}\label{cnd1}
	      \tilde{c}_{i} L_{i} \leq \tilde{c}_j L_j,
 	      \end{equation}
	           
	\item the $\{\tilde{c}_i L_i \,:\, i \in \{1,\ldots,\kappa\} \}$ follow an ascending order, i.e.,
	\begin{equation}\label{cnd2}
	\tilde{c}_1 L_1 \leq \ldots \leq \tilde{c}_i L_i \leq \ldots \leq \tilde{c}_{\kappa} L_{\kappa},
	\end{equation}
\end{enumerate}
By definition of $\cV(\boldsymbol{c};\kappa,\boldsymbol{L})$ in \eqref{mn_m}, the first condition \eqref{cnd1} implies that
\begin{equation}\label{str}
V(\kappa,K,\boldsymbol{L}) = \sum_{i=1}^{\kappa} \tilde{c}_i L_i.
\end{equation}

To prove why condition \eqref{cnd1} must hold, we apply an argument by contradiction. Let us assume that there are $i \in \{1,\ldots,\kappa\}$, and $j \in \{\kappa+1,\ldots,|\boldsymbol{L}|\}$ such that $\tilde{c}_{i} L_{i} > \tilde{c}_j L_j$, and let us denote by $j^*$ the index which corresponds to smallest such $\tilde{c}_j L_j$, in case the assumed inequality holds for more than one $j \in \{\kappa+1,\ldots,|\boldsymbol{L}|\}$, and by $i^*$ the index which corresponds to largest such $\tilde{c}_i L_i$ in case the assumed inequality holds for more than one $i \in \{1,\ldots,\kappa\}$. Clearly, $\tilde{c}_{i^*} L_{i^*} > \tilde{c}_{j^*} L_{j^*}$, and the $\tilde{c}_{j^*} L_{j^*}$ is included in the sum of $V(\kappa,K,\boldsymbol{L})$, as it is one of the $\kappa$ smallest $\{\tilde{c}_u L_u : u \in [|\boldsymbol{L}|] \}$, whereas  $\tilde{c}_{i^*} L_{i^*}$ is not included. Since $\tilde{c}_{i^*} L_{i^*} > \tilde{c}_{j^*} L_{j^*}$, and by definition $L_{i^*} \leq L_{j^*}$, then $\tilde{c}_{j^*} < \tilde{c}_{i^*}$. Thus, simply by swapping the values of $\tilde{c}_{j^*}$, $\tilde{c}_{i^*}$ we could increase the value of $V(\kappa,K,\boldsymbol{L})$ by a size of $\tilde{c}_{i^*}L_{j^*} - \tilde{c}_{j^*} L_{j^*}$, which is a contradiction because $(\tilde{c}_1,\ldots,\tilde{c}_{|\boldsymbol{L}|})$ is a maximizer. 

Assuming that \eqref{cnd1} holds which implies \eqref{str}, in order to prove why condition \eqref{cnd2} must hold, we apply again an argument by contradiction. Let us assume that there are $i < j$, both in $\{1,\ldots,\kappa\}$, such that $\tilde{c}_i L_i > \tilde{c}_j L_j$. In such a case, it is clear that $\tilde{c}_i > \tilde{c}_j$ because by definition $L_i \leq L_j$. If we swap the values of $\tilde{c}_i$, $\tilde{c}_j$, then we can increase the value of $V(\kappa,K,\boldsymbol{L})$ since
\begin{equation}
\tilde{c}_i L_i + \tilde{c}_j L_j < \tilde{c}_j L_i + \tilde{c}_i L_j \Leftrightarrow (\tilde{c}_i - \tilde{c}_j)L_i < (\tilde{c}_i - \tilde{c}_j)L_j,
\end{equation}
which is a contradiction because $(\tilde{c}_1,\ldots,\tilde{c}_{|\boldsymbol{L}|})$ is a maximizer.

The maximizer with the minimum $\cL^1$ norm is the one that satisfies \eqref{cnd2}, and \eqref{cnd1} by the following equality,
\begin{equation}\label{cnd12}
\tilde{c}_{\kappa} L_{\kappa} = \tilde{c}_j L_j, \quad \forall \, j \in \{\kappa+1,\ldots,|\boldsymbol{L}|\}.
\end{equation}

If $K$ is large enough so that
\begin{equation}\label{crs}
K \geq \kappa + L_{\kappa} \sum_{j=\kappa+1}^{|\boldsymbol{L}|} 1/L_{j},
\end{equation}
then the maximizer with the minimum $\cL^1$ norm is
\begin{equation}
\begin{aligned}
&\tilde{c}_1=\ldots=\tilde{c}_{\kappa}=1,\\
\tilde{c}_j =& L_{\kappa}/L_{j}, \quad \forall \; j \in \{\kappa+1,\ldots,|\boldsymbol{L}|\},
\end{aligned}
\end{equation}
which implies that
\begin{equation*}
V(\kappa,K,\boldsymbol{L}) = \sum_{i=1}^{\kappa} L_i,
\end{equation*}
and thus $V(\kappa,K,\boldsymbol{L})$ equals to the form \eqref{st0mb} with $x=y=0$, $v=1$, $u=\kappa$. Hence, it suffices to prove that $V(\kappa,K,\boldsymbol{L})$ equals to the form \eqref{st0mb} when
\begin{equation}\label{kk}
K < \kappa +  L_{\kappa} \sum_{j=\kappa+1}^{|\boldsymbol{L}|}1/L_{j}.
\end{equation}

When \eqref{kk} holds, the maximizer $(\tilde{c}_1,\ldots,\tilde{c}_{|\boldsymbol{L}|}) \in \cD(K)$ satisfies
\begin{equation}\label{cndk}
\sum_{i=1}^{|\boldsymbol{L}|} \tilde{c}_{i} = K,
\end{equation}
by which we get that
\begin{equation*}
K-\sum_{i=1}^{\kappa-1} \tilde{c}_{i} = \sum_{i=\kappa}^{|\boldsymbol{L}|} \tilde{c}_{i} = \tilde{c}_{\kappa} L_{\kappa} \sum_{i=\kappa}^{|\boldsymbol{L}|} 1/L_j
\end{equation*}
where the second equality follows by \eqref{cnd12}. Therefore, 
\begin{equation}\label{cnd123}
\tilde{c}_j L_j = \tilde{c}_{\kappa} L_{\kappa} = \frac{K-\sum_{i=1}^{\kappa-1} \tilde{c}_{i}}{|\boldsymbol{L}|-(\kappa-1)}\widetilde{L}_{\kappa}, \quad \forall \, j \in \{\kappa+1,\ldots,|\boldsymbol{L}|\}.
\end{equation}
where $\widetilde{L}_{\kappa}$ defined in \eqref{L2} is the harmonic mean of the $|\boldsymbol{L}|-(\kappa-1)$ largest elements in $\boldsymbol{L}$, and it stands for an average of them. It holds
\begin{equation}
0 \leq \frac{K-\sum_{i=1}^{\kappa-1} \tilde{c}_{i}}{|\boldsymbol{L}|-(\kappa-1)} \leq 1
\end{equation}
because
\begin{equation*}
K-\sum_{i=1}^{\kappa-1}\tilde{c}_{i} = \sum_{i=\kappa}^{|\boldsymbol{L}|} \tilde{c}_{i} \leq |\boldsymbol{L}|-(\kappa-1),
\end{equation*}
since $\tilde{c}_{i} \in [0,1]$ for all $i \in [|\boldsymbol{L}|]$. By replacing $\tilde{c}_{\kappa} L_{\kappa}$ in \eqref{str} we get
\begin{equation}\label{cfr}
V(\kappa,K,\boldsymbol{L}) = \tilde{c}_1 L_1 + \ldots + \tilde{c}_{\kappa-1} L_{\kappa-1} + \frac{K-\sum_{i=1}^{\kappa-1} \tilde{c}_{i}}{|\boldsymbol{L}|-(\kappa-1)}\; \widetilde{L}_{\kappa}.
\end{equation}
Therefore, the values of $\tilde{c}_1, \ldots,\tilde{c}_{\kappa-1}$ are deduced by the solution of the following maximization problem
\begin{equation}\label{mmx}
\mbox{maximize}\left\{ c_1 L_1 + \ldots + c_{\kappa-1} L_{\kappa-1} + \frac{z}{|\boldsymbol{L}|-(\kappa-1)} \widetilde{L}_{\kappa}     \right\}
\end{equation}
with respect to $c_1,\ldots,c_{\kappa-1} \in [0,1]$ and $z$ subject to 
\begin{enumerate}[(i)]
	\item 
	\begin{equation*}
	0 \leq z \leq L_{\kappa} \sum_{i=\kappa}^{|\boldsymbol{L}|} 1/L_i,
	\end{equation*}
	      
	\item       
	\begin{equation*}
	\sum_{i=1}^{\kappa-1} c_{i} + z = K,
	\end{equation*}
	
	\item 
	\begin{equation*}
	c_1 L_1 \leq \ldots \leq c_{\kappa-1} L_{\kappa-1} \leq \frac{z}{|\boldsymbol{L}|-(\kappa-1)} \widetilde{L}_{\kappa}.
	\end{equation*}
\end{enumerate}
	
In the special case $\widetilde{L}_{\kappa}/(|\boldsymbol{L}|-(\kappa-1)) \geq L_{\kappa-1}$, and since by definition $L_{\kappa-1} \geq \ldots \geq L_1$, the maximization method is to distribute $K$ so that $z, c_{\kappa-1}, \ldots, c_1$ respectively takes its' largest possible value, independently, one by one with respect to the priority list they are written. On the other hand, if $\widetilde{L}_{\kappa}/(|\boldsymbol{L}|-(\kappa-1)) < L_{\kappa-1}$ then $c_{\kappa-1}$ would become first in the priority list. Due to constraint (iii) for any value of $c_{\kappa-1}$, the value of $z$ must be large enough so that the inequality is satisfied. Thus, the maximization method for the former case does not work, as the value of $z$ depends on the value of $c_{\kappa-1}$. In order to resolve this issue, we consider 
\begin{equation*}
u^* :=\max\left\{u \in \{0,\ldots,\kappa-1\} \,:\, \frac{\kappa - u}{|\boldsymbol{L}|-u}\;\widetilde{L}_{u+1} \geq L_{u}\right\},
\end{equation*}
thus $V(\kappa,K,\boldsymbol{L})$ reduces to 
\begin{equation}
V(\kappa,K,\boldsymbol{L}) = \tilde{c}_1 L_1 + \ldots + \tilde{c}_{u^*} L_{u^*} + (\kappa -u^*)\frac{z}{|\boldsymbol{L}|-u^*} \widetilde{L}_{u^* +1},
\end{equation}
and it holds
\begin{equation*}
\frac{\kappa -u^*}{|\boldsymbol{L}|-u^*}\widetilde{L}_{u^* +1} \geq L_{u^*} \geq \ldots \geq L_{1}.
\end{equation*}

The values of $\tilde{c}_1, \ldots, \tilde{c}_{|\boldsymbol{L}|}$ are deduced by the solution of the following maximization problem,
\begin{equation}\label{mp}
\mbox{maximize}\left\{ c_1 L_1 + \ldots + c_{u^*} L_{u^*} + (\kappa -u^*)\frac{z}{|\boldsymbol{L}|-u^*} \widetilde{L}_{u^* +1}\right\}
\end{equation}
with respect to $c_1,\ldots,c_{u^*} \in [0,1]$ and $z$ subject to

\begin{enumerate}[(i)]
	\item 
	\begin{equation*}
	0 \leq z \leq L_{u^* +1} \sum_{i=u^* +1}^{|\boldsymbol{L}|} 1/L_i,
	\end{equation*}
	
	\item 
	\begin{equation*}
	\sum_{i=1}^{u^*} c_{i} + z = K,
	\end{equation*}	
	
	\item 
	\begin{equation*}
	c_1 L_1 \leq \ldots \leq c_{u^*} L_{u^*} \leq \frac{z}{|\boldsymbol{L}|-u^*} \widetilde{L}_{u^* +1}.
	\end{equation*}
\end{enumerate}

If
\begin{equation}\label{crs2}
K < L_{u^* +1}\sum_{i=u^* +1}^{|\boldsymbol{L}|} 1/L_i,
\end{equation}
then $c_1= \ldots =c_{u^*}=0$,
\begin{equation}\label{str_c}
V(\kappa,K,\boldsymbol{L})= (\kappa -u^*) \frac{K}{|\boldsymbol{L}|-u^*}\;\widetilde{L}_{u^* +1},
\end{equation}
and thus $V(\kappa,K,\boldsymbol{L})$ equals to the form \eqref{st0mb} with $x=0$, $v=0$, $u=u^*$, $y=K/(|\boldsymbol{L}|-u^*)$. Therefore, in view of \eqref{kk}, it suffices to prove that $V(\kappa,K,\boldsymbol{L})$ equals to the form \eqref{st0mb} when
\begin{equation}\label{rk}
L_{u^* +1}\sum_{i=u^* +1}^{|\boldsymbol{L}|} 1/L_i \leq K < \kappa +  L_{\kappa} \sum_{j=\kappa+1}^{|\boldsymbol{L}|}1/L_{j}.
\end{equation}

If $K$ satisfies \eqref{rk} then $c_{u^* +1}=1$, the objective of \eqref{mp} becomes
\begin{equation*}
c_1 L_1 + \ldots + c_{u^*} L_{u^*} + L_{u^* +1} + (\kappa -(u^* +1))\frac{z}{|\boldsymbol{L}|-(u^* +1)} \widetilde{L}_{u^* +2}.
\end{equation*}
and if we set $u := u^* +1$, $v:=u^* +1$ the maximization problem \eqref{mp} takes the more general form
\begin{equation}\label{mp2}
\mbox{maximize}\left\{ c_1 L_1 + \ldots + c_{v-1} L_{v-1} + \sum_{i=v}^{u} L_i + (\kappa -u)\frac{z}{|\boldsymbol{L}|-u} \widetilde{L}_{u +1}\right\}
\end{equation}
with respect to $c_1,\ldots,c_{v-1} \in [0,1]$ and $z$ subject to
\begin{enumerate}[(i)]
	\item 
	\begin{equation*}
	L_{u} \sum_{i=u +1}^{|\boldsymbol{L}|} 1/L_i \leq z \leq L_{u +1} \sum_{i=u +1}^{|\boldsymbol{L}|} 1/L_i,
	\end{equation*}
	
	\item 
	\begin{equation*}
	\sum_{i=1}^{v-1} c_{i} + z = K -(u-v-1),
	\end{equation*}	
	
	\item 
	\begin{equation*}
	c_1 L_1 \leq \ldots \leq c_{v-1} L_{v-1}\leq L_{v} \leq \ldots \leq L_{u} \leq \frac{z}{|\boldsymbol{L}|-u} \widetilde{L}_{u +1}.
	\end{equation*}
\end{enumerate}

For any $u,v \in [1,\kappa]$ with $v \leq u$, if
\begin{equation}\label{ll}
L_{v-1} \geq \frac{\kappa -u}{|\boldsymbol{L}|-u} \widetilde{L}_{u +1},
\end{equation}
then our priority is to make $c_{v-1}$ as large as possible. If $c_{v-1}=1$ then we decrease variable $v$ by $1$. We observe that the value of $c_{v-1}$ does not affect $z$ as the values $L_{v},\ldots,L_{u}$ interfere. If \eqref{ll} does not hold, our priority is to make $z$ as large as possible. If $z=L_{u+1} \sum_{i=u+1}^{|\boldsymbol{L}|} 1/L_i$ then we set $z$ equal to $L_{u+1} \sum_{i=u+2}^{|\boldsymbol{L}|} 1/L_i$, and we increase $u$ by $1$. Depending on the size of $K$ this procedure continues for a number of times, and it ends up with a form
\begin{equation}\label{vvn}
V(\kappa,K,\boldsymbol{L})=x\, L_{v-1} + \sum_{i=v}^{u} L_{i} + (\kappa-u) \frac{z}{|\boldsymbol{L}|-u} \widetilde{L}_{u+1},
\end{equation}
which for $y:=z/(|\boldsymbol{L}|-u)$ equals to the form \eqref{st0mb} in Lemma \ref{mm_g}.

In order to prove that the maximizer with the minimum $\cL^1$ norm satisfies \eqref{c1}-\eqref{c4}, we observe that \eqref{c2}-\eqref{c4} follow directly by the form \eqref{vvn} of $V(\kappa,K,\boldsymbol{L})$ by matching the $c'_{i}$ with the factor of the respective $L_i$. In order to prove \eqref{c1} we distinguish the following two cases.

\begin{enumerate}[(i)]
	\item If $u=\kappa$ then $\tilde{c}_{\kappa}=1$ and the result follows by \eqref{cnd12}.
	
	\item If $u < \kappa$ then 
	      \begin{equation*}
	      \tilde{c}_{j} L_{j} = \tilde{c}_{\kappa} L_{\kappa} = \frac{z}{|\boldsymbol{L}|-u} \widetilde{L}_{u+1}, \quad \forall \, j\, \in \, \{\kappa+1,\ldots, |\boldsymbol{L}|\},
	      \end{equation*}
	      which proves the claim since $y:=z/(|\boldsymbol{L}|-u)$.
\end{enumerate}
\end{IEEEproof}

In Algorithm \ref{alg}, we describe explicitly how $x,y,u,v$ are computed. The algorithm primarily focuses on the case that $K$ satisfies \eqref{rk}, as the other two cases are covered in \eqref{crs} and \eqref{crs2}. Algorithm \ref{alg} solves the constrained optimization problem \eqref{mp2} by implementing the steps described in the paragraph that follows \eqref{ll}.

\begin{algorithm}
	\caption{Computation of $x,y \in [0,1)$ and $u,v \in [0,\kappa]$.}
	\begin{multicols}{2}
	\begin{enumerate}[(1)]
		\item Input: $\kappa$, $K$, $\boldsymbol{L}$.
		
		\item Compute
		      \begin{equation*}
		      u^* \leftarrow \max\left\{u \in \{0,..,\kappa-1\} \,:\, \frac{\kappa -u}{|\boldsymbol{L}|-u} \widetilde{L}_{u+1} \geq L_{u}\right\}.
		      \end{equation*}
		      
	    \item If $K \geq \kappa + L_{\kappa} \sum_{j=\kappa+1}^{|\boldsymbol{L}|} 1/L_{j}$ then
	    
	          $x \leftarrow 0$, $y \leftarrow 0$, $v \leftarrow 1$, $u \leftarrow \kappa$,
	    
	          go to output
	          
	          end-if.
	          
	   \item If $K < L_{u^* +1}\sum_{i=u^* +1}^{|\boldsymbol{L}|} 1/L_i$ then
	   
	         $x \leftarrow 0$, $v \leftarrow 0$, $u \leftarrow u^*$, $y \leftarrow K/(|\boldsymbol{L}|-u^*)$,
	          
	         go to output
	         
	         end-if.
		
	   \item Initialize: $x \leftarrow 0$, $u \leftarrow u^*+1$, $v \leftarrow u^*+1$, 
		
	   \hspace{1.5 cm} $z \leftarrow L_{u} \sum_{i=u +1}^{|\boldsymbol{L}|} 1/L_i$,
	   
	   \hspace{1.5 cm} $K \leftarrow K - L_{u^* +1} \sum_{i=u^* +1}^{|\boldsymbol{L}|} 1/L_i$.
		
		\item \textbf{While} $(K > 0)$ 
		       \begin{enumerate}[(i)] 
		       	\item If $\left( L_{v-1} \geq \frac{\kappa -u}{|\boldsymbol{L}|-u} \widetilde{L}_{u+1} \right)$ or $(u=\kappa)$ then
		       	
		       	\begin{enumerate}[]
		       		\item If $(K \geq 1)$ then
		       		
		       		$v \leftarrow v-1$,
		       		
		       		$K \leftarrow K-1$,
		       		
		       		\item else
		       		
		       		$x \leftarrow K$,
		       		
		       		$K \leftarrow 0$,
		       		
		       		end-if
		       	\end{enumerate}
		       	
		       	end-if
		       	
		       	\item If $\left( L_{v-1} < \frac{\kappa -u}{|\boldsymbol{L}|-u} \widetilde{L}_{u+1} \right)$ or $(v=1)$ then
		       	
		       	$w \leftarrow L_{u+1} \sum_{i=u +1}^{|\boldsymbol{L}|} 1/L_i  -z$,
		       	
		       	$z \leftarrow z+ \min\{ K, w \}$,
		       	
		       	$K \leftarrow K - \min\{ K, w \}$.
		       	
		       	\begin{enumerate}[]
		       		\item If $\left( z \geq L_{u+1} \sum_{i=u+1}^{|\boldsymbol{L}|} 1/L_i \right)$ then
		       		
		       		$z \leftarrow L_{u+1} \sum_{i=u+2}^{|\boldsymbol{L}|} 1/L_i$,
		       		
		       		$u \leftarrow u+1$,
		       		
		       		end-if
		       	\end{enumerate}
		         	
		       	end-if
		       \end{enumerate}
	       
		end-while.
		  
		  $y \leftarrow z/(|\boldsymbol{L}|-u)$.
		    
		\item Output: $x$, $u$, $v$, $y$.
	\end{enumerate}
\end{multicols}
\label{alg}
\end{algorithm} 

\begin{IEEEproof}[Proof of Lemma \ref{equiv_opt}]
	For the proof of Lemma \ref{equiv_opt}, it suffices to show that
	\begin{equation*}
	W(\kappa_{1},\kappa_{2},K,\boldsymbol{L}_{1},\boldsymbol{L}_{2}, r) = V(\kappa_{1},K_{1}^{*},\boldsymbol{L}_{1}),
	\end{equation*}
	where $(K^{*}_{1},K^{*}_{2})$ is the maximizer with the minimum $\cL^1$ norm of the problem \eqref{rt}, then the form \eqref{wf} of $W(\kappa_{1},\kappa_{2},K,\boldsymbol{L}_{1},\boldsymbol{L}_{2}, r)$, as well as the form \eqref{ts}-\eqref{te} of the maximizer of \eqref{mm_st2} with the minimum $\cL^1$ norm, follow by Lemma \ref{mm_g}.
	
	Since $\cD(K)$ is compact and $\{U \subseteq [|\boldsymbol{L}_{1}|] \,:\, |U|=\kappa_1\}$, $\{U \subseteq [|\boldsymbol{L}_{2}|] \, :\, |U|=\kappa_2\}$ are finite sets, the max-min problem \eqref{mm_st2} has a solution. We denote by $(\hat{\boldsymbol{c}}',\check{\boldsymbol{c}}',\boldsymbol{0})$ the maximizer of \eqref{mm_st2} with the minimum $\cL^1$ norm, where  $\hat{\boldsymbol{c}}':=(\hat{c}'_{1},\ldots,\hat{c}'_{|\boldsymbol{L}_{1}|})$, $\check{\boldsymbol{c}}':=(\check{c}'_{1},\ldots,\check{c}'_{|\boldsymbol{L}_{2}|})$, and we consider 
	\begin{equation}\label{KK}
	\hat{K} := \sum_{i=1}^{|\boldsymbol{L}_{1}|} \hat{c}'_{i}, \qquad \check{K} := \sum_{i=1}^{|\boldsymbol{L}_{2}|}\check{c}'_{i}.
	\end{equation}
	The max-min structure of \eqref{mm_st2} implies that for $\hat{\boldsymbol{c}}'$, $\check{\boldsymbol{c}}'$, we have
	\begin{equation*}
	W(\kappa_{1},\kappa_{2},K,\boldsymbol{L}_{1},\boldsymbol{L}_{2}, r) = \cV(\hat{\boldsymbol{c}}';\kappa_1,\boldsymbol{L}_1) = r\, \cV(\check{\boldsymbol{c}}';\kappa_2,\boldsymbol{L}_2) 
	\end{equation*}
	and
	\begin{equation*}
	\begin{aligned}
	\cV(\hat{\boldsymbol{c}}';\kappa_1,\boldsymbol{L}_1) &= \max_{\hat{\boldsymbol{c}} \in \cD(\hat{K})} \cV(\hat{\boldsymbol{c}};\kappa_1,\boldsymbol{L}_1) = V(\kappa_1,\hat{K},\boldsymbol{L}_1),\\
	\cV(\check{\boldsymbol{c}}';\kappa_2,\boldsymbol{L}_2) &=\max_{\check{\boldsymbol{c}} \in \cD(\check{K})} \cV(\check{\boldsymbol{c}};\kappa_2,\boldsymbol{L}_2) = V(\kappa_2,\check{K},\boldsymbol{L}_2).
	\end{aligned}
	\end{equation*}
	Hence,
	\begin{equation}
	W(\kappa_{1},\kappa_{2},K,\boldsymbol{L}_{1},\boldsymbol{L}_{2}, r) = V(\kappa_1,\hat{K},\boldsymbol{L}_1) = r\, V(\kappa_2,\check{K},\boldsymbol{L}_2),
	\end{equation}
	and since $(\hat{\boldsymbol{c}}',\check{\boldsymbol{c}}',\boldsymbol{0}) \in \cD(K)$ we have $\hat{K} + \check{K} \leq K$.
	
	Therefore, it suffices to prove that $(\hat{K},\check{K})$ is the maximizer of \eqref{rt} with the minimum $\cL^1$ norm. To prove this, we apply an argument by contradiction. Let us assume that $(\hat{K},\check{K})$ is not a maximizer of \eqref{rt}, then there is $(K_{1},K_{2})$ such that $K_{1} + K_{2} \leq K$, and
	\begin{equation*}
	V(\kappa_1,K_1,\boldsymbol{L}_1) = r \, V(\kappa_2,K_2,\boldsymbol{L}_2) > V(\kappa_1,\hat{K},\boldsymbol{L}_1) = r\, V(\kappa_2,\check{K},\boldsymbol{L}_2).
	\end{equation*}
	Let us denote by $\hat{\boldsymbol{c}}'':=(\hat{c}''_{1},\ldots,\hat{c}''_{|\boldsymbol{L}_{1}|})$ the maximizer of $V(\kappa_1,K_1,\boldsymbol{L}_1)$, and by $\check{\boldsymbol{c}}'':=(\check{c}''_{1},\ldots,\check{c}''_{|\boldsymbol{L}_{2}|})$ the maximizer of $V(\kappa_2,K_2,\boldsymbol{L}_2)$. Then,
	\begin{equation*}
	\min\left\{ \cV(\hat{\boldsymbol{c}}'';\kappa_1,\boldsymbol{L}_1) , r\, \cV(\check{\boldsymbol{c}}'';\kappa_2,\boldsymbol{L}_2) \right\} > \min\left\{ \cV(\hat{\boldsymbol{c}}';\kappa_1,\boldsymbol{L}_1) , r\, \cV(\check{\boldsymbol{c}}';\kappa_2,\boldsymbol{L}_2) \right\} 
	\end{equation*} 
	which is a contradiction because we assumed that $(\hat{\boldsymbol{c}}',\check{\boldsymbol{c}}',\boldsymbol{0})$ is a maximizer of \eqref{mm_st2}. Also, if we assume that $(\hat{K},\check{K})$ is not the maximizer \eqref{rt} with the minimum $\cL^1$ norm, by \eqref{KK} we deduce that $(\hat{\boldsymbol{c}}',\check{\boldsymbol{c}}',\boldsymbol{0})$ is not the maximizer of \eqref{mm_st2} with the minimum $\cL^1$ norm, which is a contradiction. 
\end{IEEEproof}
\vspace{1 cm}
Next, we provide Algorithm \ref{alg2} for the computation of the maximizer $(K^{*}_{1},K^{*}_{2})$ of the constrained optimization problem \eqref{rt} with the minimum $\cL^1$ norm.

\begin{algorithm}
	\caption{Computation of $(K^{*}_{1},K^{*}_{2})$.}
		\begin{enumerate}[(1)]
			\item Input: $\kappa_1$, $\kappa_2$, $K$, $\boldsymbol{L}_1$, $\boldsymbol{L}_2$, $r$.\\
			
			\item Compute
			      
			      $\tilde{K}_1 \leftarrow \kappa_{1} + L_{1,\kappa_1} \sum_{i=\kappa_{1}+1}^{|\boldsymbol{L}_{1}|} 1/L_{1,i},$
			      
			      $\tilde{K}_2 \leftarrow \kappa_{2} + L_{2,\kappa_2} \sum_{i=\kappa_{2}+1}^{|\boldsymbol{L}_{2}|} 1/L_{2,i}.$

		   \item If $\sum_{i=1}^{\kappa_1} L_{1,i} \leq r\, \sum_{i=1}^{\kappa_2} L_{2,i}$
		         
		         $K^{r}_{2} \leftarrow$ root of the equation $r\, V(\kappa_2,K_2,\boldsymbol{L}_2) = V(\kappa_1,\tilde{K}_1 ,\boldsymbol{L}_1)$ with respect to $K_2 \in [0, \tilde{K}_2]$.
		         
		         $\tilde{K} \leftarrow \min\{\tilde{K}_1 + K^{r}_{2}, K\}$.
		              
			     $K^{*}_{1} \leftarrow$ root of the equation $V(\kappa_1,K_1 ,\boldsymbol{L}_1) = r\, V(\kappa_2,\tilde{K}- K_1,\boldsymbol{L}_2)$ with respect to $K_1 \in [0, \tilde{K}]$.
			     
			     $K^{*}_{2} \leftarrow \tilde{K} - K^{*}_{1}.$
			     
			     end-if
			      
		  \item If $\sum_{i=1}^{\kappa_1} L_{1,i} > r\, \sum_{i=1}^{\kappa_2} L_{2,i}$	      
			      
			    $K^{r}_{1} \leftarrow$ root of the equation $V(\kappa_1,K_1 ,\boldsymbol{L}_1) = r\, V(\kappa_2,\tilde{K}_2,\boldsymbol{L}_2)$ with respect to $K_1 \in [0, \tilde{K}_1]$.
			    
			    $\tilde{K} \leftarrow \min\{ K^{r}_1 + \tilde{K}_{2}, K\}$.
			    
			    $K^{*}_{2} \leftarrow$ root of the equation $V(\kappa_1,\tilde{K}-K_2,\boldsymbol{L}_1) = r\, V(\kappa_2,K_2,\boldsymbol{L}_2)$ with respect to $K_2 \in [0, \tilde{K}]$.
			    
			    $K^{*}_{1} \leftarrow \tilde{K} - K^{*}_{2}.$
			    
			    end-if 
			
			\item Output: $K^{*}_{1}, K^{*}_{2}$.
		\end{enumerate}
	\label{alg2}
\end{algorithm} 

We note that $K^{r}_{2}$ (resp. $K^{r}_{1}$) and $K^{*}_{1}$ (resp. $K^{*}_{2}$) are unique can be computed using the bisection method. Without loss of generality, we assume that
\begin{equation*}
\sum_{i=1}^{\kappa_1} L_{1,i} \leq r\, \sum_{i=1}^{\kappa_2} L_{2,i},
\end{equation*}
then for the computation of $K^{r}_{2}$, we consider the function
\begin{equation*}
g(K_2) := r\, V(\kappa_2,K_2,\boldsymbol{L}_2) - V(\kappa_1,\tilde{K}_1 ,\boldsymbol{L}_1), \quad K_2 \in [0, \tilde{K}_2].
\end{equation*} 
which is continuous with
\begin{equation*}
\begin{aligned}
g(0) &= -V(\kappa_1,\tilde{K}_1 ,\boldsymbol{L}_1) < 0,\\
g(\tilde{K}_2) &= r\, V(\kappa_2,\tilde{K}_2,\boldsymbol{L}_2) - V(\kappa_1,\tilde{K}_1 ,\boldsymbol{L}_1)\\
       &= r\, \sum_{i=1}^{\kappa_{2}}L_{2,i} - \sum_{i=1}^{\kappa_1} L_{1,i} \geq 0. 
\end{aligned}
\end{equation*}
If $g(\tilde{K}_2)=0$ then $K^{r}_{2} = \tilde{K}_2$, otherwise $g(\tilde{K}_2)>0$ and the solution follows by the bisection method. Since $V(\kappa_2,K_2,\boldsymbol{L}_2)$ is increasing over $[0,\tilde{K}_2]$, the function $g(K_2)$ is also increasing over $[0,\tilde{K}_2]$ and thus $K^{r}_{2}$ is unique.

For the computation of $K^{*}_1$, we consider the function
\begin{equation}
h(K_1) = V(\kappa_1,K_1 ,\boldsymbol{L}_1) - r\, V(\kappa_2,\tilde{K}-K_1,\boldsymbol{L}_2), \quad K_1 \in [0,\tilde{K}]
\end{equation}
which is continuous with
\begin{equation*}
\begin{aligned}
h(0) &= -r\, V(\kappa_2,\tilde{K},\boldsymbol{L}_2) < 0, \\
h(\tilde{K}) &= V(\kappa_1,\tilde{K} ,\boldsymbol{L}_1) >0.
\end{aligned}
\end{equation*}
Thus, the solution follows by the bisection method and it is less than or equal to $\tilde{K}_1 \wedge \tilde{K}$. Since $V(\kappa_1,K_1 ,\boldsymbol{L}_1)$ is increasing over $[0, \tilde{K}_1]$, and $V(\kappa_2,\tilde{K}-K_1,\boldsymbol{L}_2)$ is non-increasing over $[0, \tilde{K}]$ the function $h(K_1)$ is increasing over $[0, \tilde{K}_1 \wedge \tilde{K}]$ and thus $K^{*}_{1}$ is unique. 

By restricting the value of $K$ to $\tilde{K}:=(\tilde{K}_1 + K^{r}_{2})\wedge K$, we do not reduce the maximum possible value that the objective of \eqref{rt}, i.e. $V(\kappa_1,K_1,\boldsymbol{L}_1)$, can take. In the same time, we restrict the number of possible maximizers to $1$, by keeping the one with the minimum $\cL^1$ norm, which is given by the unique root of
\begin{equation*}
V(\kappa_1,K_1 ,\boldsymbol{L}_1) - r\, V(\kappa_2,\tilde{K}-K_1,\boldsymbol{L}_2), \quad K_1 \in [0,\tilde{K}].
\end{equation*}

\section{}\label{LBs}
In Appendix \ref{LBs}, we prove Theorems \ref{th:lower_bound1}, \ref{tf1}, \ref{lower_bound_bk}, which establish a universal asymptotic lower bound on the expected stopping time when controlling the misclassification error rate \eqref{error_const}, and the familywise error rates \eqref{b_error_const}, respectively.

\subsection{Misclassification error rate}

As a first step towards the proof of Theorem \ref{th:lower_bound1} we provide the following auxiliary lemma. We fix $A \subseteq [M]$, and we denote by $Z(A;k)$ the family of subsets of $[M]$ whose Hamming distance from $A$ is at least $k$, i.e.,
\begin{equation}\label{Z_unkn}
Z(A;k):=\{ C \subseteq [M]: |C \triangle A| \geq k \}.
\end{equation}

\begin{lemma}\label{a1}
	For any $B \notin Z(A;k)$, and $(c_{1},\ldots,c_{M}) \in \cD(K)$, we have the following inequality,
	\begin{equation}
	\min_{G \in Z(B;k)}\left(\sum_{i \in A \setminus G} c_{i} I_{i}  +\sum_{j \in G \setminus A} c_{j} J_{j} \right) \leq V(k,K,\boldsymbol{F}(A)).
	\end{equation}
\end{lemma}

\begin{IEEEproof}
	We fix $(c_{1},\ldots,c_{M}) \in \cD(K)$, and we denote by $C^{*}$ the set minimizer of
	\begin{equation}
	\min_{ C \in Z(A;k)}\left(\sum_{i \in A \setminus C}  c_{i} I_{i} +\sum_{j \in C \setminus A} c_{j} J_{j}  \right),
	\end{equation}
	where in place of $Z(B;k)$ we have $Z(A;k)$. By \cite[Lemma B.2]{song_fel_gener}, there exists a set $G^{*}$ such that for the sets $A,\, A \triangle C^{*},\, B$ it holds
	\begin{equation}
	A \triangle G^{*} \subseteq  A \triangle C^{*} \subseteq B \triangle G^{*}.
	\end{equation}
	By the right inclusion we have
	\begin{equation}\label{as}
	|B \triangle G^{*}| \geq |A \triangle C^{*}| \geq k,
	\end{equation}
	which means that $G^{*} \in Z(B;k)$, whereas by the left inclusion we have
	\begin{equation}\label{bs}
	A \setminus G^{*} \subseteq A \setminus C^{*}, \quad G^{*} \setminus A \subseteq C^{*} \setminus A.
	\end{equation}
	As a result,
	\begin{equation}\label{oot}
	\begin{aligned}
	\min_{G \in Z(B;k)}\left(\sum_{i \in A \setminus G} c_{i} I_{i} +\sum_{j \in G \setminus A} c_{j} J_{j} \right) \leq& \sum_{i \in A \setminus G^{*}} c_{i} I_{i}  +\sum_{j \in G^{*} \setminus A} c_{j} J_{j} \\
	\leq& \sum_{i \in A \setminus C^{*}} c_{i} I_{i} +\sum_{j \in C^{*} \setminus A} c_{j} J_{j} \\
	&=\min_{ C \in Z(A;k) }\left(\sum_{i \in A \setminus C} c_{i} I_{i} +\sum_{j \in C \setminus A} c_{j} J_{j} \right)\\
	&\leq \max_{(c_1,\ldots, c_M) \in \cD(K)} \min_{ C \in Z(A;k) }\left(\sum_{i \in A \setminus C} c_{i} I_{i} +\sum_{j \in C \setminus A} c_{j} J_{j} \right),
	\end{aligned}
	\end{equation}
	where the first inequality follows by the fact that $G^{*} \in Z(B;k)$, and the second one by \eqref{bs}. The set minimizer $C^{*} \in Z(A;k)$ has Hamming distance equal to $k$, for any $(c_{1},\ldots,c_{M}) \in \cD(K)$, as the addition of extra terms exceeds the minimum, which implies
	\begin{equation}
	\begin{aligned}
	\max_{(c_1,\ldots, c_M) \in \cD(K)}& \min_{ C \in Z(A;k) }\left(\sum_{i \in A \setminus C} c_{i} I_{i} +\sum_{j \in C \setminus A} c_{j} J_{j} \right) \\
	&= \max_{(c_{1}, \ldots, c_{M}) \in \cD(K)}  \; \min_{U \subseteq [M]: \, |U|= k} \; \sum_{i \in U} c_i \, F_{i}(A)\\
	&=V(k,K,\boldsymbol{F}(A)).
	\end{aligned}
	\end{equation}
	where the last equality follows by definition of $V(k,K,\boldsymbol{F}(A))$ according to \eqref{mm_st}. 
\end{IEEEproof}

For the proof of Theorem \ref{th:lower_bound1}, we introduce the log-likelihood ratio of $\Pro_{A}$ versus $\Pro_{C}$, for any sampling rule $R$ and any subsets $A,\,C \subseteq [M]$, based on the observations from all sources in the first $n$ sampling instants, i.e.,
\begin{align}
\Lambda^R_{A,C} (n) :=  \log  \frac{d\Pro_{A}}{d\Pro_{{C}}}  \left(\cF^R_{n} \right), \quad n \in \bN.
\end{align}
Since  $R(n)$ is $\cF^R_{n-1}$-measurable, $X_i(n)$ is  independent of $\cF^R_{n-1}$, and $Z(n)$ is  independent of $\cF^R_{n-1}$ and of $\{X_i(n)\,:\, i \in [M]\}$, we have 
\begin{align}\label{LLR_global}
\begin{split}
\Lambda^R_{A,C} (n) :=  \Lambda^R_{A,C} (n-1) 
&+ 
\sum_{i \in A\setminus C} g_i(X_i(n) )   \, R_i(n)  \\
&- \sum_{j \in C\setminus A} g_j(X_j(n) )  \, R_j(n), 
\end{split}
\end{align}
where we recall that $g_i := \log \left( f_{1i} / f_{0i}  \right)$, and we set $\Lambda^R_{A,C} (0):=0$. Comparing with \eqref{LLR}, it is clear that 
\begin{align}\label{LLR_global_repre}
\Lambda^R_{A,C} (n) = \sum_{i \in A\setminus C} \Lambda^R_{i}(n) - \sum_{j \in C\setminus A}  \Lambda^R_{j}(n), \quad n \in \bN.
\end{align}
We also set
\begin{align} \label{tilde}
\begin{split}
\tilde{\Lambda}^R_i(n) &:= \tilde{\Lambda}^R_i(n-1) + \Bigl( g_{i}(X_i(n))   - \Exp_A[ g_i(X_i(n))] \Bigr) \, R_i(n),\\
\tilde{\Lambda}^R_i(0) &:= 0,
\end{split}
\end{align} 
which implies that
\begin{align} \label{decompose}
\tilde{\Lambda}^R_i(n) &= 
\begin{cases}
\Lambda^R_i(n) - I_i \, n\pi^R_i(n), \quad i \in A, \\ 
\Lambda^R_i(n) + J_i \, n\pi^R_i(n), \quad i \notin A. 
\end{cases}
\end{align}
\\
\begin{IEEEproof}[Proof of Theorem \ref{th:lower_bound1}]
	We have to show that 
	\begin{equation}
	\mathcal{J}_{A}(\alpha;k,K) \geq \frac{|\log(\alpha)|}{V(k,K,\boldsymbol{F}(A))}(1 + o(1)),
	\end{equation}
	where $o(1)$ is a quantity that tends to zero as $\alpha \to 0$. We set
	\begin{equation}\label{ll_1}
	f(\alpha):=\frac{|\log(\alpha)|}{V(k,K,\boldsymbol{F}(A))}, \quad \alpha \in (0,1).
	\end{equation}
	By Markov's inequality, for any stopping time $T$ and $q \in (0,1)$ we have
	\begin{equation}
	\Exp_{A}[T] \geq q\, f(\alpha) \Pro_{A}(T\geq q \, f(\alpha)).
	\end{equation}
	Thus, it suffices to show that for every $q \in (0,1)$ we have
	\begin{equation}\label{suffic}
	\liminf_{\alpha \to 0} \inf_{(R,T,\Delta) \in \mathcal{C}(\alpha;k,K) }\Pro_{A}(T \geq q\, f(\alpha))\geq 1,
	\end{equation}
	as this will imply that 
	\begin{equation}
	\liminf\limits_{\alpha \to 0} \frac{\mathcal{J}_{A}(\alpha;k,K)}{ |\log(\alpha)|} \geq \frac{q}{V(k,K,\boldsymbol{F}(A))},
	\end{equation}
	and the desired result will follow by letting $q \to 1$. 
	
	In the rest of the proof we fix some arbitrary $q \in (0,1)$. Then, for any $\alpha \in \, (0,1)$, $(R,T,\Delta) \in \mathcal{C}(\alpha;k,K)$, and $B \notin Z(A;k)$, where $Z(A;k)$ is defined in \eqref{Z_unkn}, we have
	\begin{equation}\label{terms}
	\begin{aligned}
	\Pro_{A}(\Delta=B) \leq & \Pro_{A}\left(\min_{  G \in Z(B;k) }\Lambda^{R}_{A,G}(T)<\log\left(\frac{\eta}{\alpha}\right),\Delta=B\right)\\
	&+\Pro_{A}\left( T \leq q f(\alpha)  ,\, \min_{  G \in Z(B;k)}\Lambda^{R}_{A,G}(T) \geq \log\left(\frac{\eta}{\alpha}\right),\Delta=B\right)\\
	&+\Pro_{A}\left(T \geq q f(\alpha),\Delta=B\right),
	\end{aligned}
	\end{equation}
	where $\eta$ is an arbitrary constant in $(0,1)$. By summing up \eqref{terms} with respect to all $B \notin Z(A;k)$, we have
	\begin{equation}
	\begin{aligned}
	\Pro_{A}(\Delta \notin Z(A;k))\leq & \sum_{B \notin Z(A;k)} \Pro_{A}\left(\min_{G \in Z(B;k)}\Lambda^{R}_{A,G}(T)<\log\left(\frac{\eta}{\alpha}\right),\Delta=B\right)\\
	&+ \sum_{B \notin Z(A;k)} \Pro_{A}\left(T \leq q f(\alpha),\min_{G \in Z(B;k)}\Lambda^{R}_{A,G}(T) \geq \log\left(\frac{\eta}{\alpha}\right),\Delta=B \right)\\
	&+\Pro_{A}\left(T \geq q f(\alpha),\Delta \notin Z(A;k)\right).
	\end{aligned}
	\end{equation}
	In view of the fact that
	\begin{equation}
	1-\alpha \leq \Pro_{A}(\Delta \notin Z(A;k)),
	\end{equation}
	we obtain
	\begin{equation}
	\begin{aligned}
	\Pro_{A}&\left(T \geq q f(\alpha)\right)\\ 
	&\geq \Pro_{A}\left(T \geq q f(\alpha),\Delta \notin Z(A;k)\right)\\ 
	&\geq 1-\alpha -\sum_{B \notin Z(A;k)} \Pro_{A}\left(\min_{G \in Z(B;k)}\Lambda^{R}_{A,G}(T)<\log\left(\frac{\eta}{\alpha}\right),\Delta=B\right)\\
	&\qquad \qquad -\sum_{B \notin Z(A;k)}\Pro_{A}\left(T \leq q f(\alpha),\min_{G \in Z(B;k)}\Lambda^{R}_{A,G}(T) \geq \log\left(\frac{\eta}{\alpha}\right),\Delta=B \right).
	\end{aligned}
	\end{equation}
	Thus, in order to show  \eqref{suffic}, it suffices to show that for any $B \notin Z(A;k)$
	\begin{equation}\label{t_cl1}
	\lim_{\alpha \to 0}\sup_{(R,T,\Delta) \in \mathcal{C}(\alpha;k,K) } \Pro_{A}\left(\min_{G \in Z(B;k)}\Lambda^{R}_{A,G}(T)<\log\left(\frac{\eta}{\alpha}\right),\Delta=B\right)=0,
	\end{equation}
	and 
	\begin{equation}\label{t_cl2}
	\lim_{\alpha \to 0}\sup_{(R,T,\Delta) \in \mathcal{C}(\alpha;k,K)} \Pro_{A}\left(T \leq q f(\alpha),\min_{G \in Z(B;k)}\Lambda^{R}_{A,G}(T) \geq \log\left(\frac{\eta}{\alpha}\right),\Delta=B\right) =0.
	\end{equation}
	In order to show \eqref{t_cl1}, we fix $\alpha \in \, (0,1)$ and $(R,T,\Delta) \in \mathcal{C}(\alpha;k,K)$. By application of Boole's inequality we have
	\begin{equation}\label{star1}
	\Pro_{A}\left(\min_{G \in Z(B;k)} \Lambda^{R}_{A,G}(T)<\log\left(\frac{\eta}{\alpha}\right),\Delta=B \right) \leq \sum_{G \in Z(B;k)}\Pro_{A}\left(\Lambda^{R}_{A,G}(T)<\log\left(\frac{\eta}{\alpha}\right),\Delta=B \right).
	\end{equation}
	For all $G \in Z(B;k)$, we apply the change of measure $\Pro_{A} \to \Pro_{G}$ and by Wald's likelihood ratio identity  we obtain
	\begin{equation}
	\begin{aligned}
	\Pro_{A}\left(\Lambda^{R}_{A,G}(T)<\log\left(\frac{\eta}{\alpha} \right),\Delta=B \right) &= \Exp_{G} \left[\exp\{\Lambda^{R}_{A,G}(T)\};{\Lambda^{R}_{A,G}}(T)<\log\left(\frac{\eta}{\alpha}\right),\Delta=B \right]\\
	&\leq \frac{\eta}{\alpha} \, \Pro_{G}(\Delta=B) \leq \eta,
	\end{aligned}
	\end{equation}
	where the last inequality is deduced by the fact that for any $G \in Z(B;k)$, it holds $\Pro_{G}(\Delta=B) \leq \alpha$. In view of \eqref{star1}, we obtain
	\begin{equation}\label{terms1}
	\Pro_{A}\left(\min_{G \in Z(B;k)}{\Lambda^{R}_{A,G}}(T)<\log\left(\frac{\eta}{\alpha}\right),\Delta=B \right) \leq \big{|} Z(B;k)\big{|}\eta.
	\end{equation}
	Since  $\eta \in (0,1)$ is arbitrary,  letting $\eta \to 0$, we prove \eqref{t_cl1}.\\
	
	In order to show \eqref{t_cl2}, we observe that by decomposition \eqref{decompose} we have
	\begin{align*}
	\frac{1}{T}\Lambda^{R}_{A,G}(T) &= \frac{1}{T} \left(\sum_{i \in A \setminus G} \tilde{\Lambda}^{R}_{i}(T) + \sum_{j \in G \setminus A} -\tilde{\Lambda}^{R}_{j}(T)\right)+ \sum_{i \in A \setminus G}I_{i} \, \pi^{R}_{i}(T) + \sum_{j \in G \setminus A} J_{j} \, \pi^{R}_{j}(T)  \\
	&\leq \frac{1}{T} \sum_{i \in [M]} |\tilde{\Lambda}^{R}_{i}(T) | + \sum_{i \in A \setminus G}I_{i} \, \pi^{R}_{i}(T) + \sum_{j \in G \setminus A} J_{j}  \, \pi^{R}_{i}(T) .
	\end{align*} 
	Considering the minimum over all $G \in Z(B;k)$ on both sides of the above equality, we have
	\begin{equation}\label{xii_ineq}
	\begin{aligned}
	\frac{1}{T}	\min_{G \in Z(B;k)} \Lambda^{R}_{A,G}(T) 
	& \leq  \frac{1}{T} \sum_{i \in [M]} |\tilde{\Lambda}^{R}_{i}(T) |+\min_{G \in Z(B;k)}\left(\sum_{i \in A \setminus G}I_{i} \, \pi^{R}_{i}(T) +\sum_{j \in G \setminus A}J_{j} \, \pi^{R}_{j}(T) \right) \\
	& \leq \frac{1}{T} \sum_{i \in [M]} |\tilde{\Lambda}^{R}_{i}(T) |+ V(k,K,\boldsymbol{F}(A)),
	\end{aligned}
	\end{equation}
	where the second inequality follows by Lemma \ref{a1} and the fact that $(R,T,\Delta)$ belongs to $\cC(K)$, which implies 
	$(\pi^{R}_{1}(T),\ldots, \pi^{R}_{M}(T)) \in \cD(K)$. Therefore,
	
	\begin{equation}\label{oth}
	\begin{aligned}
	\Pro_{A}&\left(T \leq q f(\alpha),\min_{G \in Z(B;k)}\Lambda^{R}_{A,G}(T) \geq \log\left(\frac{\eta}{\alpha}\right),\Delta=B \right) \\	
	&\leq \Pro_{A}\left(T \leq q f(\alpha) , \, 
	\xi(T) \geq |\log \alpha|+\log \eta \right),
	\end{aligned}
	\end{equation}
	where
	\begin{equation}\label{W_A_B}
	\begin{aligned}
	\xi(n):=\Bigg( \sum_{i \in [M]} \left| \frac{\tilde{\Lambda}^{R}_{i}(n)}{n} \right| + V(k,K,\boldsymbol{F}(A))\Bigg )n,\quad \forall n \in \mathbb{N}.
	\end{aligned}
	\end{equation}
	By the moment assumption \eqref{momr} and \cite[Theorem 2.19]{hall14}, we have
	\begin{equation}
	\lim_{n \to \infty} \frac{\tilde{\Lambda}^{R}_{i}(n)}{n}  = 0, \quad \mbox{a.s.} \quad \forall \,i \in \, [M],
	\end{equation}
	and as a result,
	\begin{equation}
	\lim_{n \to \infty} \frac{\xi(n)}{n} = V(k,K,\boldsymbol{F}(A)), \quad \mbox{a.s.}  
	\end{equation}
	Hence, by \cite[Lemma F.1]{song_fel_gener} we have
	\begin{equation}
	\lim_{\alpha \to 0}\sup_{(R,T,\Delta) \in \mathcal{C}(\alpha;k,K)} \Pro_{A}\left(T \leq q f(\alpha) , \xi(T) \geq |\log(\alpha)|+\log(\eta) \right) =0,
	\end{equation}
	From this and \eqref{oth}, we conclude that \eqref{t_cl2} holds.
\end{IEEEproof}

\subsection{Familywise error rates}
As a first step towards the proof of Theorems \ref{tf1}, \ref{lower_bound_bk}, we provide the following auxiliary lemma. To this end, for fixed $A \subseteq [M]$, we consider the set
\begin{equation}\label{er_bot_set}
H_{k_{1},k_{2}}(A):=\{ B\subseteq [M]: |B \setminus A| <k_{1}, \, |A\setminus B|<k_{2} \},
\end{equation}
and for any $B \in H_{k_{1},k_{2}}(A) $, we also consider the sets
\begin{equation}
\begin{aligned}
U_{k_{1}}(B)&:=\{ C\subseteq A: |B\setminus C|\geq k_{1} \},\\ 
Y_{k_{2}}(B)&:=\{ C\supseteq A: |C\setminus B|\geq k_{2} \}.
\end{aligned}
\end{equation}

\begin{lemma}\label{a2}
	For any $B \in H_{k_{1},k_{2}}(A)$ and $(c_1,\ldots,c_M) \in \cD(K)$, we have the following inequalities:
	
	\begin{enumerate}[(i)]
		\item If $A=\emptyset$,
		      \begin{equation}\label{sst1}
		      \min_{G \in Y_{k_2}(B)}\sum_{j \in G} c_{j} J_{j} \leq V(k_2,K,\boldsymbol{J_{k_1-1}}(A)).
		      \end{equation}
		   
		\item If $A=[M]$,
		      \begin{equation}\label{sst2}
		      \min_{G \in U_{k_1}(B)}\sum_{i \in [M]\setminus G} c_{i} \, I_{i} \leq V(k_1,K,\boldsymbol{I_{k_2-1}}(A)).
		      \end{equation}
		
		\item If $0<|A|<M$,
		      \begin{equation}\label{sst3}
		      \min\Bigg\{\min_{G \in U_{k_1}(B)}\sum_{i \in A\setminus G} c_{i} \, I_{i}, \,r\min_{G \in Y_{k_2}(B)}\sum_{j \in G \setminus A} c_{j} \, J_{j}\Bigg\} \leq v_{A}(k_1,k_2,K,r),
		      \end{equation}
		      where $v_{A}(k_1,k_2,K,r)$ is defined in Theorem \ref{lower_bound_bk}.
	\end{enumerate}
\end{lemma}

\begin{IEEEproof}
We prove (i), (iii). Case (ii) is symmetric to (i).
\\
For $A=\emptyset$, without loss of generality, we consider $B \in H_{k_{1},k_{2}}(\emptyset)$ such that $Y_{k_2}(B) \neq \emptyset$, or equivalently $|B|<k_1$ and $|B^c| \geq k_2$, otherwise \eqref{sst1} holds trivially. Since $|B^c| \geq k_2$, there exists $\Gamma_{2} \subseteq B^c$ such that $|\Gamma_{2}|=k_2$, and $\Gamma_{2}$ contains the sources that correspond to the $k_{2}$ smallest terms in $\{c_{j} J_j\,:\, j \in B^c \}$. Clearly, $\Gamma_{2}\in Y_{k_2}(B)$, and
\begin{equation}
[M] = B \cup B^c \supseteq B \cup \Gamma_{2}.
\end{equation}	
Since $|B|<k_1$, even if $B$ contains the $k_1 -1$ smallest elements in $\{c_{j} J_j\,:\, j\in [M] \}$, the set $\Gamma_{2}$ contains the following $k_2$ smallest elements in $[M]\setminus B = B^c$. Thus, let us denote by $\{j\}$ the identity of the source with the $j^{th}$ smallest element in $\{c_{i} J_i\,:\, i \in [M] \}$, then we have
\begin{equation}
\begin{aligned}
\min_{G \in Y_{k_2}(B)} \sum_{j \in G} c_{j} \, J_j \leq & \sum_{j \in \Gamma_{2}} c_{j} \, J_j \\
\leq & \sum_{j=k_1}^{k_2 +(k_1 -1)} c_{\{j\}} \, J_{\{j\}} \leq V(k_2,K,\boldsymbol{J_{k_1-1}}(A)).
\end{aligned}
\end{equation}	
\\
For $0<|A|<M$, without loss of generality, we consider $B \in H_{k_{1},k_{2}}(A)$ such that $U_{k_{1}}(B) \neq \emptyset$ and $Y_{k_2}(B) \neq \emptyset$, or equivalently $|B| \geq k_1$ and $|B^c|\geq k_2$, otherwise \eqref{sst3} holds trivially. We consider the quantities 
\begin{equation}
l_1:=|B \setminus A|, \quad  l_{2}:=|A \setminus B|.
\end{equation}
\begin{enumerate}
	\item The fact that $|B|\geq k_{1}$ implies that $|A\cap B|\geq k_{1}-|B\setminus A|=k_{1}-l_{1}$. Thus, there exists $\Gamma_{1} \subseteq A\cap B$ such that $|\Gamma_{1}|=k_{1}-l_{1}$, and $\Gamma_{1}$ contains the sources that correspond to the $k_{1}-l_{1}$ smallest elements in $\{ c_i I_i\,:\, i\in A\cap B \}$. We set $B^{*}_{1}:=A \setminus \Gamma_{1}$. Since $\Gamma_{1} \subseteq A\cap B$, it holds $A\setminus B^{*}_{1}=\Gamma_{1}$, and
	\begin{equation}
	B\setminus B^{*}_{1}=B \cap (A\setminus \Gamma_{1})^c = \Gamma_{1} \cup \left( B \setminus A \right).
	\end{equation}
	Therefore,
	\begin{equation*}
	|B\setminus B^{*}_{1}|=|\Gamma_{1}| + |B \setminus A| = k_1 -l_1 + l_1 = k_1,
	\end{equation*}
	which implies $B^{*}_{1} \in U_{k_{1}}(B)$, and
	\begin{equation}
	\min_{G \in U_{k_1}(B)}\sum_{i \in A\setminus G} c_{i} \, I_{i} \leq \sum_{i \in \Gamma_{1}} c_{i} \, I_{i}.
	\end{equation}
	Even if $A\setminus B$ contained the $l_2$ smallest elements in $\{ c_{i} I_{i}\,:\, i\in A \}$, the set $\Gamma_{1} \subseteq A \cap B$ would contain the following $k_1 -l_1$ smallest elements in the same set. Thus, for $<j>$ to denote the identity of the $j^{th}$ smallest element in $\{ c_{i} I_{i}\,:\, i\in A \}$, we always have 
	\begin{equation}\label{s11}
	\sum_{i \in \Gamma_{1}} c_{i} \, I_{i} \leq \sum_{i=1+l_2}^{k_1 -l_1 +l_2} c_{<i>} \,  I_{<i>}.
	\end{equation}
	
	\item The fact that $|B^c| \geq k_{2}$ implies that $|A^c \cap B^c|\geq k_{2} - |A\cap B^c|=k_{2}-l_{2}$. Thus, there exists $\Gamma_{2} \subseteq A^c \cap B^c$ such that $|\Gamma_{2}|=k_{2}-l_{2}$ and $\Gamma_{2}$ contains the sources that correspond to the $k_{2}-l_{2}$  smallest elements in $\{c_{j} J_j\,:\, j\in A^c \cap B^c \}$. We set $B^{*}_{2}:=A\cup \Gamma_{2}$. Since $\Gamma_{2} \subseteq A^c \cap B^c$, it holds $B^{*}_{2} \setminus A=\Gamma_{2}$, and
	\begin{equation}
	B^{*}_{2} \setminus B = (A\cup \Gamma_{2}) \cap B^c = \Gamma_{2} \cup \left( A \setminus B \right).
	\end{equation}
	Therefore, 
	\begin{equation*}
	|B^{*}_{2} \setminus B| = |\Gamma_{2}| + |A\setminus B| = k_2 - l_2 + l_2 = k_2,
	\end{equation*}
	which implies $B^{*}_{2} \in Y_{k_{2}}(B)$, and
	\begin{equation}
	\min_{G \in Y_{k_2}(B)}\sum_{j \in G \setminus A} c_{j} \, J_{j} \leq \sum_{j \in \Gamma_{2}}  c_{j} \, J_{j}.
	\end{equation}
	Even if $B\setminus A$ contained the $l_1$ smallest elements in $\{c_{j} J_j\,:\, j\in A^c \}$, the set $\Gamma_{2}\subseteq A^c \cap B^c$ would contain the following $k_2 - l_2$ smallest elements in the same set. Thus, for $\{j\}$ to denote the identity of the source with the $j^{th}$ smallest element in $\{c_{j} J_j\,:\, j\in A^c\}$, we always have
	\begin{equation}\label{s22}
	\sum_{j \in \Gamma_{2}}  c_{j} \, J_{j} \leq \sum_{j=1+l_{1}}^{k_{2}-l_{2}+l_{1}} c_{\{j\}} \, J_{\{j\}}.
	\end{equation}	
\end{enumerate}

Let us assume, without loss of generality, that $l_1 \geq l_2$. We set $l:=l_1 - l_2$, and for \eqref{s11}, \eqref{s22} we further have,
\begin{equation}
\begin{aligned}
\sum_{i=1+l_{2}}^{k_{1}-l_{1}+l_{2}} c_{<i>} \, I_{<i>} &\leq \sum_{i=1+l_{2}}^{k_{1}-l} c_{<i>} \, I_{<i>} \leq \sum_{i=1}^{k_{1}-l} c_{<i>} \, I_{<i>},\\
\sum_{j=1+l_{1}}^{k_{2}-l_{2}+l_{1}}c_{\{j\}} \,J_{\{j\}} &\leq \sum_{j=1+l+l_{2}}^{k_{2}+l}c_{\{j\}}\,J_{\{j\}} \leq \sum_{j=1+l}^{k_{2}+l}c_{\{j\}}\,J_{\{j\}}.
\end{aligned}
\end{equation}
As a result,
\begin{equation*}
\begin{aligned}
\min&\left\{ \min_{G \in U_{k_1}(B)} \sum_{i \in A\setminus G} c_{i} I_{i}, \, r \,\min_{G \in Y_{k_2}(B)}\sum_{j \in G \setminus A} c_{j} \, J_{j}  \right\}\\ 
& \leq \min\left\{ \sum_{i \in \Gamma_{1}} c_{i} \, I_{i}, r \, \sum_{j \in \Gamma_{2}}  c_{j} \, J_{j} \right\}\\
& \leq \min\left\{ \sum_{i=1}^{k_{1}-l} c_{<i>} \, I_{<i>},r\, \sum_{j=1+l}^{k_{2}+l} c_{\{j\}} \, J_{\{j\}} \right\}\\
& \leq W(k_1-l,k_2,K,\boldsymbol{I}(A) ,\boldsymbol{J_{l}}(A), r ) \\
& \leq v_{A}(k_1,k_2,K,r),
\end{aligned}
\end{equation*}
which concludes the claim.
\end{IEEEproof}

We proceed to the proof of Theorem \ref{lower_bound_bk} and by which we deduce the proof of Theorem \ref{tf1}, as explained in the last part of the following proof. 

\begin{IEEEproof}[Proof of Theorem \ref{lower_bound_bk}]
    We have to show that
    \begin{equation}
    \mathcal{J}_{A}(\alpha,\beta;k_{1},k_{2},K) \geq \frac{|\log(\alpha)|}{v_{A}(k_{1},k_{2},K,r)}(1+o(1)),
    \end{equation}
    where $v_{A}(k_{1},k_{2},K,r)$ is defined in Theorem \ref{lower_bound_bk} and $o(1)$ is a quantity that tends to zero as $\alpha \to 0$. We define
    \begin{equation}
    f(\alpha):=\frac{|\log(\alpha)|}{v_{A}(k_{1},k_{2},K,r)}, \quad \alpha \in (0,1).
    \end{equation}
    By Markov's inequality, for any stopping time $T$ and $q, \alpha \in (0,1)$,
    \begin{equation}
    \Exp_{A}[T] \geq q \, f(\alpha) \Pro_{A}(T\geq q\, f(\alpha)).
    \end{equation}
    Thus, it suffices to show that for every $q \in (0,1)$, we have
    \begin{equation}\label{suffic0}
    \liminf_{\alpha \to 0} \inf_{(R,T,\Delta) \in \mathcal{C}(\alpha,\beta;k_1,k_2,K) }\Pro_{A}(T \geq q\, f(\alpha))\geq 1,
    \end{equation}
    as this will imply that 
    \begin{equation}
    \liminf\limits_{\alpha \to 0}\mathcal{J}_{A}(\alpha,\beta;k_{1},k_{2},K)/ |\log(\alpha)| \geq q/v_{A}(k_{1},k_{2},K,r),
    \end{equation}
    and the desired result will follow by letting $q \to 1$. 
    
    In the rest of the proof, we fix some arbitrary $q \in (0,1)$. Moreover, we note that for any $B \subseteq [M]$ we have either 
    $|B| \geq k_{1}$ or $|B^c| \geq k_{2}$, because otherwise $M=|B|+|B^c|<k_{1}+k_{2}$, which contradicts the assumption $k_{1}+k_{2} \leq M$. In what follows, we let $B \in H_{k_{1},k_{2}}(A)$ and we focus on the case $|B|\geq k_{1}$ and $|B^c| \geq k_{2}$. The other two cases $|B|\geq k_{1},\, |B^c|<k_{2}$, and $|B|<k_{1},\, |B^c| \geq k_{2}$ are simpler and they can be treated in the same way as described in the last part of this proof.   
    
    Then, for any $\alpha \in\, (0,1)$ and  $(R,T,\Delta) \in \mathcal{C}(\alpha,\beta;k_1,k_2,K)$ we have
    \begin{equation}\label{terms0}
    \begin{aligned}
    \Pro_{A}\left( \Delta=B \right) \leq& \Pro_{A}\left(\min_{G \in U_{k_1}(B)} \Lambda^{R}_{A,G}(T)<\log\left(\frac{\eta}{\alpha}\right),\Delta=B\right)\\
    +&\Pro_{A}\left(\min_{G \in Y_{k_2}(B)} \Lambda^{R}_{A,G}(T)<\log\left(\frac{\eta}{\beta}\right),\Delta=B\right) \\
    +&\Pro_{A}\left(\min_{G \in U_{k_1}(B)} \Lambda^{R}_{A,G}(T)\geq \log\left(\frac{\eta}{\alpha}\right),\min_{G \in Y_{k_2}(B)} \Lambda^{R}_{A,G}(T)\geq \log\left(\frac{\eta}{\beta}\right),\Delta=B\right),
    \end{aligned}
    \end{equation} 
    where $\Lambda^{R}_{A,G}(T)$ is defined in \eqref{LLR_global}, and $\eta$ is an arbitrary constant in $(0,1)$. The third term in \eqref{terms0} can be equivalently written as
    \begin{equation}
    \Pro_{A}\left( \min\Big\{ \min_{G \in U_{k_1}(B)} \Lambda^{R}_{A,G}(T), \rho(\alpha,\eta) \min_{G \in Y_{k_2}(B)} \Lambda^{R}_{A,G}(T) \Big\} \geq \log\left(\frac{\eta}{\alpha}\right), \Delta=B \right),
    \end{equation}
    where $\rho(\alpha,\eta):=\log(\eta / \alpha)/ \log(\eta / \beta)$.
    
    For simplicity in what follows we just write $\rho$ instead of $\rho(\alpha,\eta)$. We upper bound the third term in \eqref{terms0} by
    \begin{equation}\label{l2}
    \begin{aligned}
    \Pro_{A}&\left( \min\Big\{ \min_{G \in U_{k_1}(B)} \Lambda^{R}_{A,G}(T),\rho \min_{G \in Y_{k_2}(B)} \Lambda^{R}_{A,G}(T) \Big\} \geq \log\left(\frac{\eta}{\alpha}\right), \Delta=B \right)\\
    \leq& \Pro_{A}\left( T\leq q f(\alpha) ,\min\Big\{ \min_{G \in U_{k_1}(B)} \Lambda^{R}_{A,G}(T),\rho \min_{G \in Y_{k_2}(B)} \Lambda^{R}_{A,G}(T) \Big\} \geq \log\left(\frac{\eta}{\alpha}\right), \Delta=B \right)\\
    &+\Pro_{A}(T \geq q f(\alpha),\Delta=B).
    \end{aligned}
    \end{equation}
    By summing up \eqref{terms0} over all $B \in H_{k_{1},k_{2}}(A)$, we have 
    \begin{equation*}
    \begin{aligned}
    \Pro_{A}&\left( \Delta \in H_{k_{1},k_{2}}(A) \right)\\ 
    &\leq \sum_{B \in H_{k_{1},k_{2}}(A)} \Pro_{A}\left( \min_{G \in U_{k_1}(B)} \Lambda^{R}_{A,G}(T)< \log\left(\frac{\eta}{\alpha}\right), \Delta=B \right)\\
    &+ \sum_{B \in H_{k_{1},k_{2}}(A)} \Pro_{A}\left( \min_{G \in Y_{k_2}(B)} \Lambda^{R}_{A,G}(T)< \log\left(\frac{\eta}{\beta}\right), \Delta=B \right)\\
    &+ \sum_{B \in H_{k_{1},k_{2}}(A)} \Pro_{A}\left( T\leq q f(\alpha) ,\min\Big\{ \min_{G \in U_{k_1}(B)} \Lambda^{R}_{A,G}(T),\rho \min_{G \in Y_{k_2}(B)} \Lambda^{R}_{A,G}(T) \Big\} \geq \log\left(\frac{\eta}{\alpha}\right), \Delta=B \right)\\
    &+ \Pro_{A}(T \geq q f(\alpha),\Delta \in H_{k_{1},k_{2}}(A)).
    \end{aligned}
    \end{equation*}
    In view of the fact that
    \begin{equation}
    1-(\alpha+\beta) \leq \Pro_{A}\left( \Delta \in H_{k_{1},k_{2}}(A) \right),
    \end{equation}
    we obtain
    \begin{equation*}
    \begin{aligned}
    \Pro_{A}&\left(T \geq q f(\alpha) \right)\\ 
    & \geq  \Pro_{A}\left(T \geq q f(\alpha),\Delta \in H_{k_{1},k_{2}}(A) \right) \\
 &    \geq  1-(\alpha+\beta)\\
    &- \sum_{B \in H_{k_{1},k_{2}}(A)} \Pro_{A}\left( \min_{G \in U_{k_1}(B)} \Lambda^{R}_{A,G}(T)< \log\left(\frac{\eta}{\alpha}\right), \Delta=B \right)\\
    &- \sum_{B \in H_{k_{1},k_{2}}(A)} \Pro_{A}\left( \min_{G \in Y_{k_2}(B)} \Lambda^{R}_{A,G}(T)< \log\left(\frac{\eta}{\beta}\right), \Delta=B \right)\\
    &- \sum_{B \in H_{k_{1},k_{2}}(A)} \Pro_{A}\left( T\leq q f(\alpha) ,\min\Big\{ \min_{G \in U_{k_1}(B)} \Lambda^{R}_{A,G}(T),\rho \min_{G \in Y_{k_2}(B)} \Lambda^{R}_{A,G}(T) \Big\} \geq \log\left(\frac{\eta}{\alpha}\right), \Delta=B \right).
    \end{aligned}
    \end{equation*}
    Thus, in order to show \eqref{suffic0} it suffices to show that for all $B \in H_{k_{1},k_{2}}(A)$
    \begin{equation}\label{sh1}
    \begin{aligned}
    \lim_{\alpha \to 0} \sup_{(R,T,\Delta) \in \mathcal{C}(\alpha,\beta;k_1,k_2,K)} \Pro_{A}\left( \min_{G \in U_{k_1}(B)} \Lambda^{R}_{A,G}(T)< \log\left(\frac{\eta}{\alpha}\right), \Delta=B \right) = 0,\\
    \lim_{\beta \to 0} \sup_{(R,T,\Delta) \in \mathcal{C}(\alpha,\beta;k_1,k_2,K)} \Pro_{A}\left( \min_{G \in Y_{k_2}(B)} \Lambda^{R}_{A,G}(T)< \log\left(\frac{\eta}{\beta}\right), \Delta=B \right) = 0,
    \end{aligned}
    \end{equation}
    and 
    \begin{equation}\label{sh2}
    \begin{aligned}	
    \lim_{\alpha  \to 0} \sup_{(R,T,\Delta)} &\Pro_{A}\left( T{\leq} q f(\alpha) ,\min\{ \min_{G \in U_{k_1}(B)} \Lambda^{R}_{A,G}(T),\rho \min_{G \in Y_{k_2}(B)} \Lambda^{R}_{A,G}(T) \} {\geq} \log\left(\frac{\eta}{\alpha}\right), \Delta{=}B \right)\\
    &=0,
    \end{aligned}	
    \end{equation}
    where the supremum is evaluated over $(R,T,\Delta) \in \mathcal{C}(\alpha,\beta;k_1,k_2,K)$. In order to show \eqref{sh1}, we apply Boole's inequality and we have:
    \begin{equation}
    \begin{aligned}
    \Pro_{A}\left(\min_{G \in U_{k_1}(B)}\Lambda^{R}_{A,G}(T)<\log\left(\frac{\eta}{\alpha}\right),\Delta=B\right) &\leq \sum_{G \in U_{k_1}(B)}\Pro_{A}\left( \Lambda^{R}_{A,G}(T)<\log\left(\frac{\eta}{\alpha}\right),\Delta=B\right), \\
    \Pro_{A}\left(\min_{G \in Y_{k_2}(B)} \Lambda^{R}_{A,G}(T)<\log\left(\frac{\eta}{\beta}\right),\Delta=B\right) &\leq \sum_{G \in Y_{k_2}(B)}  \Pro_{A}\left(\Lambda^{R}_{A,G}(T)<\log\left(\frac{\eta}{\beta}\right),\Delta=B\right).
    \end{aligned}
    \end{equation}
    For all $G \in U_{k_1}(B)$, we apply the change measure $\Pro_{A} \to \Pro_{G}$, and we have
    \begin{equation}
    \begin{aligned}
    \Pro_{A}\left(\Lambda^{R}_{A,G}(T)<\log\left(\frac{\eta}{\alpha}\right),\Delta=B\right)&=\Exp_{G}\left[\exp\{\Lambda^{R}_{A,G}(T) \}; \Lambda^{R}_{A,G}(T)<\log\left(\frac{\eta}{\alpha}\right),\Delta=B \right]\\
    &\leq \frac{\eta}{\alpha}\Pro_{G}(\Delta =B) \leq \frac{\eta}{\alpha}\Pro_{G}\left( |B\setminus G|\geq k_{1} \right)=\eta,
    \end{aligned}
    \end{equation}
    where the last inequality is deduced by the fact that for any $G \in U_{k_1}(B)$, the error constraint \eqref{b_error_const} implies $\Pro_{G}\left( |B\setminus G|\geq k_{1} \right) \leq \alpha$. Therefore, for every $\eta\, \in (0,1)$,
    \begin{equation}
    \limsup_{\alpha \to 0} \sup_{(R,T,\Delta)\in \mathcal{C}(\alpha,\beta;k_1,k_2,K)} \Pro_{A}\left(\min_{G \in U_{k_1}(B)}\Lambda^{R}_{A,G}(T)<\log\left(\frac{\eta}{\alpha}\right),\Delta=B\right) \leq |U_{k_1}(B)|\eta. 
    \end{equation}
    In the same way, we show that
    \begin{equation}
    \limsup_{\beta \to 0} \sup_{(R,T,\Delta)\in \mathcal{C}(\alpha,\beta;k_1,k_2,K)} \Pro_{A}\left(\min_{G \in Y_{k_2}(B)} \Lambda^{R}_{A,G}(T)<\log\left(\frac{\eta}{\beta}\right),\Delta=B\right)\leq |Y_{k_2}(B)| \eta.
    \end{equation}
    Letting $\eta \to 0$, we prove \eqref{sh1}. 
    
    In order to show \eqref{sh2}, we observe that by decomposition  \eqref{decompose} we have
    \begin{equation}
    \begin{aligned}
    \frac{1}{T} \, \Lambda^{R}_{A,G}(T)&= \sum_{i \in A\setminus G} \frac{\tilde{\Lambda}^{R}_{i}(T)}{T}+I_{i}\, \pi^{R}_{i}(T), &\quad \forall \,G \in U_{k_1}(B),\\
    \frac{1}{T}\Lambda^{R}_{A,G}(T)&= \sum_{j \in G\setminus A} \frac{-\tilde{\Lambda}^{R}_{j}(T)}{T}+J_{j}\, \pi^{R}_{j}(T), &\quad \forall \, G \in Y_{k_2}(B),
    \end{aligned}
    \end{equation}
    which further implies that
    \begin{equation}
    \begin{aligned}
    \frac{1}{T}\min&\left\{\min_{G \in U_{k_1}(B)}\Lambda^{R}_{A,G}(T),\rho \, \min_{G \in Y_{k_2}(B)}\Lambda^{R}_{A,G}(T)\right\}\\
    \leq&  \max\{1,\rho \}\sum_{i \in [M]} \frac{|\tilde{\Lambda}^{R}_{i}(T)|}{T}\\
    &+\min\Bigg\{\min_{G \in U_{k_1}(B)}\sum_{i \in A\setminus G} I_{i} \, \pi^{R}_{i}(T),\rho \,  \min_{G \in Y_{k_2}(B)}\sum_{j \in G \setminus A}J_{j} \, \pi^{R}_{j}(T)\Bigg\}.
    \end{aligned}
    \end{equation}
    We note that for any fixed $\eta > 0$, 
    \begin{equation}
    \rho(\alpha,\eta)=\frac{(\log(\eta)/\log(\alpha) -1)\log(\alpha)}{(\log(\eta)/\log(\beta) -1)\log(\beta)} \to r, \quad \mbox{as }\, \alpha \to 0,
    \end{equation}
    where $r$ is defined in \eqref{r}. Also, by Lemma \ref{a2} we have
    \begin{equation}
    \min\Bigg\{\min_{G \in U_{k_1}(B)}\sum_{i \in A\setminus G}I_{i}\pi^{R}_{i}(T),\,r  \min_{G \in Y_{k_2}(B)}\sum_{j \in G \setminus A}J_{j}\pi^{R}_{j}(T)\Bigg\} \leq v_{A}(k_1,k_2,K,r).
    \end{equation}
    Therefore, 
    \begin{equation}
    \frac{1}{T}\min\left\{\min_{G \in U_{k_1}(B)}\Lambda^{R}_{A,G}(T),r \min_{G \in Y_{k_2}(B)}\Lambda^{R}_{A,G}(T)\right\}
    \leq \max\{1,r \}\sum_{i \in [M]} \frac{|\tilde{\Lambda}^{R}_{i}(T)|}{T} + v_{A}(k_1,k_2,K,r).
    \end{equation}
    As a result, the  probability in \eqref{sh2} is bounded above by 
    \begin{equation}
    \Pro_{A}\left(T \leq q f(\alpha) ,\, \xi(T) \geq |\log(\alpha)|+\log(\eta) \right), 
    \end{equation}
    where
    \begin{equation}\label{e_a_b}
    \xi(n):=\left( \max\{1,r\}\sum_{i \in [M]} \frac{\big|\tilde{\Lambda}^{R}_{i}(n)\big|}{n} + v_{A}(k_1,k_2,K,r) \right) n, \quad n \in \mathbb{N},
    \end{equation}
    and it suffices to show that 
    \begin{equation}\label{ohn}
    \limsup_{\alpha  \to 0} \sup_{(R,T,\Delta)\in \mathcal{C}(\alpha,\beta;k_1,k_2,K)} \Pro_{A}\left(T \leq q f(\alpha) ,\, \xi(T) \geq |\log(\alpha)|+\log(\eta) \right)=0.
    \end{equation}
    In order to show \eqref{ohn}, it suffices to prove that 
    \begin{equation}\label{va12}
    \lim_{n\to \infty} \frac{\xi(n)}{n} = v_{A}(k_1,k_2,K,r), \quad \mbox{a.s.}
    \end{equation}
    and then the claim follows by \cite[Lemma F.1]{song_fel_gener}. Indeed, by \cite[Theorem 2.19]{hall14} and moment assumption \eqref{momr} we have
    \begin{equation}
    \lim_{n\to \infty} \frac{\tilde{\Lambda}^{R}_{i}(n)}{n}  = 0, \quad \mbox{a.s.} \quad \forall \, i \in [M],
    \end{equation}
    which combined with \eqref{e_a_b}, implies \eqref{va12}.\\
    
    \textit{The other cases}: As for the cases $|B|\geq k_{1},\, |B^c|<k_{2}$, and $|B|<k_{1},\, |B^c|\geq k_{2}$ we note the following. 
    \begin{itemize}
    	\item In case $|B| \geq k_{1}$, $|B^c|<k_{2}$, i.e., $Y_{k_2}(B) = \emptyset$ and for any $B \in H_{k_1,k_2}(A)$, we have
    	\begin{equation}
    	\begin{aligned}
    	\Pro_{A}(\Delta=B) \leq \,& \Pro_{A}\left(\min_{G \in U_{k_1}(B)}{\Lambda^{R}_{A,G}}(T)<\log\left(\frac{\eta}{\alpha}\right),\Delta=B\right) \\
    	&+\Pro_{A}\left(\min_{G \in U_{k_1}(B)}{\Lambda^{R}_{A,G}}(T)\geq \log\left(\frac{\eta}{\alpha}\right),\Delta=B\right).
    	\end{aligned}
    	\end{equation}
    	This also applies for the case $A=[M]$.
    	
    	\item In case $|B^c|\geq k_{2}$, $|B|<k_{1}$, i.e., $U_{k_1}(B)= \emptyset$ and for any $B \in H_{k_1,k_2}(A)$, we have
    	\begin{equation}
    	\begin{aligned}
    	\Pro_{A}(\Delta=B) \leq \, & \Pro_{A}\left(\min_{G \in Y_{k_2}(B)} \Lambda^{R}_{A,G}(T)<\log\left(\frac{\eta}{\beta}\right),\Delta=B\right) \\
    	&+\Pro_{A}\left(\min_{G \in Y_{k_2}(B)} \Lambda^{R}_{A,G}(T)\geq \log\left(\frac{\eta}{\beta}\right),\Delta=B\right).
    	\end{aligned}
    	\end{equation}
    	This also applies for the case $A=\emptyset$.
    \end{itemize}
    In both cases, we deduce the claim following the same reasoning as for the main case ($|B|\geq k_{1},\, |B^c|\geq k_{2}$).
\end{IEEEproof}

\section{}\label{AO}
In Appendix \ref{AO}, we prove Theorems \ref{asy_opt1}, \ref{asy_opt}, Theorems \ref{cons}, \ref{c_st}, and Proposition \ref{cst_rp}. Throughout Appendix \ref{AO}, we fix $A \subseteq [M]$.

\begin{IEEEproof}[Proof of Theorem \ref{asy_opt1}]
	In view of the asymptotic lower bound in Theorem \ref{th:lower_bound1}, it suffices to show that 
	\begin{equation}\label{tsh1}
	\Exp_{A}[T^{R}] \lesssim \frac{|\log(\alpha)|}{V(k,K,\boldsymbol{F}(A))}, \quad \mbox{as }\; \alpha \to 0,
	\end{equation}
	where $T^R$ is the stopping rule defined in \eqref{st_r}. For $d>0$ selected according to \eqref{unk_thr}, it suffices to show that for an arbitrarily small $\epsilon > 0$, we have
	\begin{equation} \label{ubk}
	\Exp_{A}[T^{R}] \lesssim  d \Big{/} \left(\sum_{i=1}^{k} ( c^{*}_{(i)}(A)-\epsilon)F_{i}(A)-\epsilon \right), \quad \mbox{as }\; d \to \infty, 
	\end{equation}
	where $c^{*}_{(i)}(A)$ are defined in \eqref{cv01}, and by letting $\epsilon \to 0$ we obtain \eqref{tsh1}. To this end, we fix $\epsilon > 0$ small enough, and we set
	\begin{equation}
	L_\epsilon(d):=\frac{d}{ \sum_{i=1}^{k} ( c^{*}_{(i)}(A)-\epsilon)F_{i}(A)-\epsilon},
	\end{equation}
	for which it holds
	\begin{equation}\label{ttn}
	\Exp_A[T^R] \leq L_\epsilon(d)+ \sum_{n>L_\epsilon(d)} \Pro_A(T^R>n).
	\end{equation}
	On the event $\{T^R>n\}$, there exists $U \subseteq [M]$ with $|U|=k$ such that 
	\begin{equation}
	\sum_{i \in A \cap U} \Lambda^R_i(n) - \sum_{j \in A^c \cap U} \Lambda^R_j(n) \leq \sum_{i \in U} |\Lambda^R_i(n)| < d.
	\end{equation}
	Moreover, for every $n > L_\epsilon(d)$ we have
	\begin{equation*}
	d < n \, \left(\sum_{i=1}^{k} ( c^{*}_{(i)}(A)-\epsilon) F_{i}(A)-\epsilon \right).
	\end{equation*}
	By Lemma \ref{mm_g}, for every set $U \subseteq [M]$ with $|U|=k$, we have
	\begin{equation*}
	V(k,K,\boldsymbol{F}(A)) = \sum_{i=1}^{k} c^{*}_{(i)}(A) \, F_{i}(A) \leq \sum_{i \in A \cap U } c^{*}_{i}(A) \, I_i  +  \sum_{j \in A^c \cap U } c^{*}_{j}(A) \, J_j.
	\end{equation*}
	which further implies that there is an $\epsilon' >0$ sufficiently smaller than $\epsilon$ such that
	\begin{equation*}
	d < n \, \left( \sum_{i \in A \cap U } (c^{*}_{i}(A)-\epsilon') I_i  +  \sum_{j \in A^c \cap U } (c^{*}_{j}(A)-\epsilon') J_j        -\epsilon' \right).
	\end{equation*}
	Therefore, for every $n > L_\epsilon(d)$ the event $\{T^R>n\}$ is included in the event
	\begin{equation*}
	\left\{\sum_{i \in A \cap U} \Lambda^R_i(n) - \sum_{j \in A^c \cap U} \Lambda^R_j(n) <  n  \left( \sum_{i \in A \cap U } (c^{*}_{i}(A)-\epsilon') I_i  +  \sum_{j \in A^c \cap U } (c^{*}_{j}(A)-\epsilon') J_j -\epsilon'\right) \right\}.
	\end{equation*}
	In view of \eqref{ttn} and by application of Boole's inequality over all $U\subseteq [M]$ such that $|U|=k$, we deduce that 
	\begin{equation}\label{series_bound}
	\Exp_A[T^R] \leq  L_\epsilon(d) + \sum_{ \{U\subseteq [M]: |U|=k\} } \sum_{n=1}^{\infty} \; S(n;U),
	\end{equation}
	where
	\begin{equation}
	\begin{aligned}
	S(n;U) := &\sum_{i\in A \cap U } \Pro_{A}\left( \frac{\Lambda^R_i(n)}{n}< (c^{*}_{i}(A)-\epsilon')I_{i}-\epsilon'/k \right) \\
	&+ \sum_{j\in A^c \cap U } \Pro_{A}\left( -\frac{\Lambda^R_j(n) }{n} < (c^{*}_{j}(A)-\epsilon')J_{j} -\epsilon'/k \right).
	\end{aligned}
	\end{equation}
	By assumption, for any $i \in A$, $j \notin A$, and for any $\epsilon > 0$, we have 
	\begin{equation}\label{condition}
	\begin{aligned}
	\sum_{n=1}^\infty \Pro_A\left( \pi^{R}_i(n) < c_i^*(A) -\epsilon \right) <\infty,\\
	\sum_{n=1}^\infty \Pro_A\left( \pi^{R}_j(n) < c_j^*(A) -\epsilon  \right) <\infty, 
	\end{aligned}
	\end{equation}
	and thus by \cite[Lemma A.2 (i), (ii)]{a_prob}, with $\rho_i= c_i^*(A)- \epsilon'$ and  
	$\rho_j=  c_j^*(A)-\epsilon'$, it follows that all the series in \eqref{series_bound} converge. Hence, letting $d \to \infty$ we obtain  \eqref{ubk}.
\end{IEEEproof}

\begin{IEEEproof}[Proof of Theorem \ref{asy_opt}]
In view of the asymptotic lower bound in Theorem \ref{lower_bound_bk}, it suffices to show that 
\begin{equation}\label{tsh2}
\Exp_{A}\left[ T^{R}_{leap} \right] \lesssim \frac{|\log(\alpha)|}{v_{A}(k_1,k_2,K,r)}, \quad \mbox{as }\; \alpha \to 0,
\end{equation}	
where  $T^R_{leap}$ is the stopping rule \eqref{tr_leap}. Without loss of generality, we focus on case $0 < |A| <M$ and $v_{A}(k_1,k_2,K,r)$ equals \eqref{v3} of Theorem \ref{lower_bound_bk}, as the other cases follow by the same approach. For $ a,\,b > 0$ selected according to \eqref{ab_thrs} it suffices to show that for any arbitrarily small $\epsilon >0$, we have
\begin{equation}\label{21e}
\Exp_{A}[T^{R}] \lesssim b \Big{/} \left( \sum_{i=1}^{k_1-l_A} (c^*_{<i>}(A) -\epsilon) I_{i}(A) -\epsilon \right), \quad \mbox{as }\, b \to \infty,
\end{equation}
where $c^{*}_{<i>}(A)$, $c^{*}_{\{i\}}(A)$ are defined as in \eqref{cv02}-\eqref{cv03}, and by letting $\epsilon \to 0$ we obtain \eqref{tsh2}. To this end, we fix $\epsilon >0$ small enough and we set
\begin{equation}\label{t1n}
L_{\epsilon}(a,b):= \max \left\{  \frac{b}{\sum_{i=1}^{k_1 -l_A} (c^{*}_{\langle i \rangle}(A)- \epsilon) I_{i}(A) - \epsilon}  ,  \frac{a}{\sum_{j=1}^{k_2} (c^{*}_{\{j\}}(A) -\epsilon) J_{l_A+j}(A)-\epsilon} \right\},
\end{equation}
and we observe that	
\begin{equation}\label{t2n}
\Exp_{A} \left[ T^R_{leap} \right] \leq L_{\epsilon}(a,b) + \sum_{n>L_{\epsilon}(a,b) } \Pro_{A}(T^R > n).
\end{equation}	
For any $n\in \mathbb{N}$, on the event $\{ T^R_{leap} > n \}$ there exist $U_{1} \subseteq A$ with $|U_1|=k_1 -l_A$, and $U_{2} \subseteq \boldsymbol{J_{l_A}}(A)$ with $|U_2|=k_2$, such that 
\begin{equation}
\sum_{i\in U_1} \hat{\Lambda}^{R}_{i}(n) < b, \quad \mbox{or} \quad \sum_{i \in U_2} \check{\Lambda}^{R}_{i}(n) < a.
\end{equation}
Moreover, for every $n > L_{\epsilon}(a,b)$ we have
\begin{equation}
b < n \left(  \sum_{i=1}^{k_1 - l_A}( c^{*}_{\langle i \rangle}(A)-\epsilon) I_{i}(A)-\epsilon \right),\qquad a < n \left( \sum_{j=1}^{k_2}(c^{*}_{\{j\}}(A) - \epsilon) J_{l_A +j}(A)-\epsilon \right).
\end{equation}	
For any set $U_1$, $U_2$, as described above, we have
\begin{equation*}
\sum_{i=1}^{k_1 - l_A} c^{*}_{\langle i \rangle}(A)\, I_{i}(A) \leq \sum_{i \in U_1} c^{*}_{i}(A)\, I_i, \qquad
\sum_{j=1}^{k_2} c^{*}_{\{j\}}(A)\, J_{j+l_A}(A) \leq  \sum_{j \in U_2} c^{*}_{j}(A)\, J_j.
\end{equation*}
Thus, there is an $\epsilon' > 0$ sufficiently smaller than $\epsilon$ such that 
\begin{equation}
b < n \left( \sum_{i \in U_1}( c^{*}_{i}(A)  - \epsilon') I_{i}-\epsilon' \right),\qquad a < n \left( \sum_{j \in U_2}(c^{*}_{j}(A) - \epsilon') J_{j} -\epsilon' \right).
\end{equation}
Therefore, for every $n > L_{\epsilon}(a,b)$ the event $\{ T^R_{leap} > n \}$ is included in the event
\begin{equation*}
\left\{\sum_{i\in U_1} \frac{\hat{\Lambda}^{R}_{i}(n)}{n} < \sum_{i \in U_1}( c^{*}_{i}(A)  - \epsilon') I_{i}-\epsilon'\right\} \bigcup \left\{\sum_{j \in U_2} \frac{\check{\Lambda}^{R}_{j}(n)}{n} < \sum_{j \in U_2}(c^{*}_{j}(A) - \epsilon') J_{j} -\epsilon'\right\}.
\end{equation*}	
In view of \eqref{t2n} and by application of Boole's inequality over all $U_{1} \subseteq A$ with $|U_1|=k_1 -l_A$, and all $U_{2} \subseteq \boldsymbol{J_{l_A}}(A)$ with $|U_2|=k_2$, we deduce that 
\begin{equation}\label{sum_ser}
\begin{aligned}
\Exp_{A} \left[ T^R_{leap} \right] \leq  L_{\epsilon}(a,b) &+ \sum_{U_1 \subseteq A :\, |U_1|=k_1-l_A} \sum_{n=1}^{\infty} \Pro_{A}\left(\sum_{i\in U_1} \frac{\hat{\Lambda}^{R}_{i}(n)}{n} < \sum_{i \in U_1}( c^{*}_{i}(A)  - \epsilon') I_{i}-\epsilon' \right)\\
&+ \sum_{U_2 \subseteq \boldsymbol{J_{l_A}}(A) : \, |U_2|=k_2} \sum_{n=1}^{\infty} \Pro_{A}\left(\sum_{j \in U_2} \frac{\check{\Lambda}^{R}_{j}(n)}{n} < \sum_{j \in U_2}(c^{*}_{j}(A) - \epsilon') J_{j} -\epsilon' \right).
\end{aligned}
\end{equation}	
	In order to prove the claim, it suffices to show the summability of the series in \eqref{sum_ser}. By Boole's inequality we have
	\begin{equation}\label{ttf}
	\begin{aligned}
	\Pro_{A}\left(\sum_{i\in U_1} \frac{\hat{\Lambda}^{R}_{i}(n)}{n} < \sum_{i \in U_1}( c^{*}_{i}(A)  - \epsilon') I_{i}-\epsilon' \right) &\leq \sum_{i \in U_1} \Pro_{A}\left(\frac{\hat{\Lambda}^{R}_{i}(n)}{n} < ( c^{*}_{i}(A)  - \epsilon') I_{i}-\epsilon'/k_1 \right),\\
	\Pro_{A}\left(\sum_{j \in U_2} \frac{\check{\Lambda}^{R}_{j}(n)}{n} < \sum_{j \in U_2}(c^{*}_{j}(A) - \epsilon') J_{j} -\epsilon' \right) & \leq \sum_{j \in U_2} \Pro_{A}\left( \frac{\check{\Lambda}^{R}_{j}(n)}{n} < (c^{*}_{j}(A) - \epsilon') J_{j}   -\epsilon'/k_2 \right).
	\end{aligned}
	\end{equation}
	By assumption, for any $i \in  A$, $j \notin A$, and for any $\epsilon >0$ we have
	\begin{equation}
	\begin{aligned}
	\sum_{n=1}^\infty \Pro_A\left( \pi^{R}_i(n) < c_i^*(A) -\epsilon \right) &<\infty,\\ 
	\sum_{n=1}^\infty \Pro_A\left( \pi^{R}_j(n) < c_j^*(A) -\epsilon  \right) &<\infty, 
	\end{aligned}
	\end{equation}
	and thus by \cite[Lemma A.2 (i) (ii)]{a_prob}, with $\rho_i= c_i^*(A)- \epsilon'$ and $\rho_j=  c_j^*(A)-\epsilon'$, it follows that all the series in  \eqref{ttf} converge. Thus, letting $b \to \infty$ we obtain \eqref{21e}.\\
\end{IEEEproof}

\begin{IEEEproof}[Proof of Theorem \ref{cons}]
	By definition of $\mathfrak{D}^{R}_n$, we have	
	\begin{equation}
	\begin{aligned}
	\Pro_{A}\left( \sigma_A^R > n \right) & \leq \sum_{i \in A} \Pro_{A}\left( \exists \, m \geq n : \Lambda^{R}_{i}(m) < 0\right)\\
	&+ \sum_{j \notin A} \Pro_{A}\left( \exists \, m \geq n : \Lambda^{R}_{j}(m) \geq 0\right),
	\end{aligned}
	\end{equation}
	which by Boole's inequality is further bounded by
	\begin{equation}
	\begin{aligned}
	\Pro_{A}\left( \sigma_A^R > n \right) \leq & \sum_{i \in A} \sum_{m=n}^{\infty} \Pro_{A}\left(\Lambda^{R}_{i}(m) < 0\right)\\
	&+ \sum_{j \notin A} \sum_{m=n}^{\infty} \Pro_{A}\left( \Lambda^{R}_{j}(m) \geq 0\right).
	\end{aligned}
	\end{equation}
	Therefore, in order to prove that $\Exp_{A}\left[\sigma_A^R\right] < \infty$, it suffices to show that
	\begin{equation}\label{e00}
	\begin{aligned}
	\sum_{n=1}^{\infty}\sum_{m=n}^{\infty} \Pro_{A}\left(\Lambda^{R}_{i}(m) < 0\right)=\sum_{n=1}^{\infty} n \, \Pro_{A}\left(\Lambda^{R}_{i}(n) < 0\right) &< \infty, \quad \forall \, i \in A,\\
	\sum_{n=1}^{\infty}\sum_{m=n}^{\infty} \Pro_{A}\left( \Lambda^{R}_{j}(m) \geq 0\right)=\sum_{n=1}^{\infty} n \, \Pro_{A}\left(\Lambda^{R}_{j}(n) \geq 0\right) &< \infty, \quad \forall \, j \in A^c.
	\end{aligned}
	\end{equation}
	We prove the case $i \in A$, as the case $i \in A^c$ follows in the same way. We note that
	\begin{equation}
	\begin{aligned} 
	\Pro_{A}\left(\Lambda^{R}_{i}(n) < 0\right)  \leq &\, \Pro_{A}\left(\Lambda^{R}_{i}(n) < 0,\, \pi^{R}_{i}(n) \geq C\, n^{-\delta}\right)\\
	&+\Pro_{A}\left( \pi^{R}_{i}(n)< C\, n^{-\delta}\right).
	\end{aligned} 
	\end{equation}
	In view of assumption \eqref{spar}, in order to prove \eqref{e00}, it suffices to show that
	\begin{equation}\label{tts1}
	\sum_{n=1}^{\infty} n\, \Pro_{A}\left(\Lambda^{R}_{i}(n) < 0,\, \pi^{R}_{i}(n) \geq C\, n^{-\delta}\right) < \infty,
	\end{equation}
	and since
	\begin{equation}
	\begin{aligned}
	&\Pro_{A}\left(\Lambda^{R}_{i}(n) < 0,\, \pi^{R}_{i}(n) \geq C\, n^{-\delta}\right) \\
	&= \Pro_{A}\left(\tilde{\Lambda}^{R}_{i}(n) < -I_i \, n \, \pi^{R}_{i}(n), \pi^{R}_{i}(n) \geq C\, n^{-\delta}\right)\\
	&\leq \Pro_{A}\left(|\tilde{\Lambda}^{R}_{i}(n)| > C \, I_i  \, n^{1-\delta}\right),
	\end{aligned}
	\end{equation}
	it suffices to show
	\begin{equation}\label{tts2}
	\sum_{n=1}^{\infty} n\,\Pro_{A}\left(|\tilde{\Lambda}^{R}_{i}(n)| > C I_i n^{1-\delta}\right) < \infty.
	\end{equation}
   By Markov's inequality
   \begin{equation}
   	\Pro_{A}\left(|\tilde{\Lambda}^{R}_{i}(n)| >C I_i n^{1-\delta}\right) \leq \frac{\Exp_{A}\left[|\tilde{\Lambda}^{R}_{i}(n)|^{\mathfrak{p}} \right]}{C^{\mathfrak{p}} I^{\mathfrak{p}}_{i} n^{(1-\delta)\mathfrak{p}}}.
   \end{equation}
   For each $i \in [M]$, $\{ \tilde{\Lambda}^{R}_{i}(n)\, :\, n \geq 0\}$ is a $\cF^{R}(n)$-martingale, thus by Rosenthal's inequality \cite[Theorem 2.12]{hall14} there is a constant $C_{0} > 0$, such that
   \begin{equation}
   	\Exp_{A}\left[|\tilde{\Lambda}^{R}_{i}(n)|^{\mathfrak{p}} \right] \leq C_{0} n^{\mathfrak{p}/2},
   \end{equation}
   and as a result,
   \begin{equation}
   	\Pro_{A}\left(|\tilde{\Lambda}^{R}_{i}(n)| >C I_i n^{1-\delta}\right) \leq \frac{C_0}{C^{\mathfrak{p}}I^{\mathfrak{p}}_{i}} \frac{n^{\mathfrak{p}/2}}{n^{(1-\delta)\mathfrak{p}}},
   \end{equation}
   For $\epsilon >0$ small enough such that $\delta < 1/2-(2+\epsilon)/\mathfrak{p}$, it holds
   \begin{equation}
   	\frac{n^{\mathfrak{p}/2}}{n^{(1-\delta)\mathfrak{p}}} < \frac{1}{n^{2+\epsilon}}, 
   \end{equation}
   which proves the claim \eqref{tts2}.
\end{IEEEproof}

\begin{IEEEproof}[Proof of Theorem \ref{c_st}(i)]
According to Theorems \ref{asy_opt1} and \ref{asy_opt}, it suffices to show that for each $i \in [M]$
\begin{equation}
\sum_{n=1}^\infty \Pro_{A} \left(\pi_i^{R}(n)< c^*_i(A)- \epsilon \right)< \infty, \quad \forall \; \epsilon > 0,
\end{equation} 
where $(c^*_{1}(A), \ldots, c^*_{M}(A))$ is defined according to Definitions \ref{def_c_star1} and \ref{def_c_star2}, respectively.

By the definition \eqref{q} of a probabilistic sampling rule, $R(n)$ is conditionally independent of $\cF^{R}_{n-1}$ given $\mathfrak{D}^{R}_{n-1}$. Thus, by \cite[Prop. 6.13]{kallenberg2002foundations} there is a measurable function $h:\bN \times 2^{[M]}\times [0,1] \to 2^{[M]}$ such that
\begin{equation}
R(n) = h(n,\mathfrak{D}^{R}_{n-1}, Z_{0}(n-1)), \quad n \in \bN,
\end{equation}
where $\{Z_{0}(n) \, :\, n \in \bN_0\}$ is a sequence of \textit{iid} random variables, uniformly distributed on $[0,1]$. Consequently, for each $i \in [M]$ there is a measurable function  $h_{i} :\bN \times 2^{[M]}\times [0,1] \to \{0,1\}$ such that
\begin{equation}\label{rh_st}
R_{i}(n) = h_{i}(n,\mathfrak{D}^{R}_{n-1}, Z_{0}(n-1)), \quad n \in \bN.
\end{equation}
We fix $\epsilon > 0$, and $i \in [M]$. For every $n \in \bN$, we have
\begin{equation}
\begin{aligned}
&\left\{\pi_i^{R}(n)< c^*_i(A)-\epsilon\right\} \\
&= \left\{\sum_{m=1}^{n}\left(R_{i}(m)-h_i(m,A, Z_{0}(m-1))\right) + \sum_{m=1}^{n}h_i(m,A,Z_{0}(m-1)) < n \, (c^*_i(A)- \epsilon) \right\},
\end{aligned}
\end{equation}
and as a result
\begin{equation}\label{ll_st}
\begin{aligned}
\Pro_{A} \left(\pi_i^{R}(n)< c^*_i(A)- \epsilon \right) &\leq \Pro_{A}\left(\sum_{m=1}^{n}\left(R_{i}(m)-h_i(m,A,Z_{0}(m-1))\right) < -n\, \epsilon/2\right)\\
&+\Pro_{A}\left(\sum_{m=1}^{n}h_i(m,A,Z_{0}(m-1))<n(c^*_i(A)-\epsilon/2)\right).
\end{aligned}
\end{equation}
For the first term on the right hand side of \eqref{ll_st} we have
\begin{equation}\label{on_st}
\begin{aligned}
&\Pro_{A}\left(\sum_{m=1}^{n}\left(R_{i}(m)-h_i(m,A,Z_{0}(m-1))\right) < -n\, \epsilon/2\right)\\
&\leq \Pro_{A}\left(\sum_{m=1}^{n} |R_{i}(m)-h_i(m,A,Z_{0}(m-1))|> n\, \epsilon/2, \sigma^{R}_A \leq n \right) + \Pro_{A}\left(\sigma^{R}_A \geq n\right)\\
&\leq \Pro_{A}\left(\sigma^{R}_A \geq n\, \epsilon/2\right) + \Pro_{A}\left(\sigma^{R}_A \geq n \right)
\end{aligned}
\end{equation}
where we used the fact that  $R_{i}(n)=h_i(n,A,Z_{0}(n-1))$
for all $n \geq \sigma^{R}_A$. Since $R$ is consistent both terms on the right hand side of \eqref{on_st} are summable.

For the second term on the right hand side of \eqref{ll_st} we have
\begin{equation}
\begin{aligned}
&\left\{\sum_{m=1}^{n}h_i(m,A,Z_{0}(m-1))<n(c^*_i(A)- \epsilon/2)\right\}\\
&=\left\{\sum_{m=1}^{n}\left(h_i(m,A,Z_{0}(m-1))-c_{i}^{R}(m,A)\right) <  \sum_{m=1}^{n} \left(c^*_i(A)-c_{i}^{R}(m,A)\right) -n \, \epsilon/2  \right\}.
\end{aligned}
\end{equation}

By the (general form of) Stolz-Cesaro theorem \cite{cesaro1888convergence},
\begin{equation}
\begin{aligned}
&\limsup_{n \to \infty} \frac{\sum_{m=1}^{n} \left(c^*_i(A)-c_{i}^{R}(m,A)\right)}{n} \\
&\leq  \limsup_{n \to \infty} \left(c^*_i(A)-c_{i}^{R}(n,A)\right) =c^*_i(A)- \liminf_{n \to \infty} c_{i}^{R}(n,A) \leq 0,
\end{aligned}
\end{equation}
where the last inequality is deduced by assumption \eqref{cr_geq}. Therefore, there exists an  $M >0$ such that for all $n \geq M$ it holds
\begin{equation}
\frac{1}{n} \sum_{m=1}^{n} \left(c^*_i(A)-c_{i}^{R}(m,A)\right) < \frac{\epsilon}{4}.
\end{equation}
Hence,
\begin{equation}\label{ts2}
\begin{aligned}
&\sum_{n=1}^{\infty} \Pro_{A}\left(\sum_{m=1}^{n}h_i(m,A, Z_{0}(m-1))<n \, (c^*_i(A)- \epsilon/2)\right)\\
&\leq M+\sum_{n=M}^{\infty} \Pro_{A}\left(\sum_{m=1}^{n}\left(h_i(m,A, Z_{0}(m-1))-c_{i}^{R}(m,A) \right) < -n \, \epsilon/4\right),
\end{aligned}
\end{equation}
and, thus, it suffices to show that
\begin{equation}\label{ts3}
\sum_{n=M}^{\infty} \Pro_{A}\left(\sum_{m=1}^{n}\left(h_i(m,A,Z_{0}(m-1))-c_{i}^{R}(m,A)\right) < -n \epsilon/4\right) < \infty.
\end{equation}
By \eqref{c} and \eqref{rh_st}, it holds
\begin{equation}
\Exp_{A}\left[h_i(n,A,Z_{0}(n-1)) \, | \, \cF^{R}_{n-1} \right] = c_{i}^{R}(n,A), \quad \forall \, n  \in \bN
\end{equation}
which shows that $\{h_i(n,A,Z_{0}(n-1))-c_{i}^{R}(n,A) \, :\, n\geq 1\}$ is a martingale difference sequence, whose absolute value is uniformly bounded by $1$. By application of Azuma-Hoeffding inequality, we deduce that
\begin{equation}
\Pro_{A}\left(\sum_{m=1}^{n}\left(h_i(m,A,Z_{0}(m-1))-c_{i}^{R}(m,A)\right) < -n \, \epsilon/2\right) \leq e^{-\gamma n}
\end{equation}
where $\gamma >0$ is a constant, which proves \eqref{ts3}.\\
\end{IEEEproof}

\begin{IEEEproof}[Proof of Proposition \ref{cst_rp}(ii)]
We show that the suggested sampling rule $R$ is consistent for part (ii). According to Theorem \ref{cons}, it suffices to show that for each $i \in [M]$,
\begin{equation}\label{ts_st}
\sum_{n=1}^{\infty} n\, \Pro_{A}\left(\pi^{R}_{i}(n) < C\, n^{-\delta}\right) < \infty,
\end{equation} 
for some $C>0$ and $\delta \in \left(0,\frac{1}{2}-\frac{2}{\mathfrak{p}}\right)$. We choose $C < C_p$. We fix $i \in [M]$, and for every $n \in \bN$, we notice that
\begin{equation}
\begin{aligned}
\left\{ \pi_{i}^{R}(n) < C\, n^{-\delta}\right\} &=\left\{ \sum_{m=1}^{n} R_{i}(m) < C n^{1-\delta} \right\}\\
&=\left\{ \sum_{m=1}^{n} (R_{i}(m)-c^{R}_{i}(m,\mathfrak{D}^{R}_{m-1})) +\sum_{m=1}^{n} c^{R}_{i}(m,\mathfrak{D}^{R}_{m-1}) < C n^{1-\delta}\right\}.
\end{aligned}
\end{equation}
By \eqref{c_gd}, we have $c^{R}_{i}(m,\mathfrak{D}^{R}_{m-1}) \geq C_{p}/m^{\delta}$, for all $m \in \{1,\ldots,n\}$, and we deduce that
\begin{equation}
\begin{aligned}
\sum_{m=1}^{n} c^{R}_{i}(m,\mathfrak{D}^{R}_{m-1}) &\geq \sum_{m=1}^{n} \frac{C_p}{m^{\delta}}  \geq C_{p} \, n^{1-\delta},
\end{aligned}
\end{equation}
which further implies
\begin{equation}
\begin{aligned}
&\Pro_{A}\left(\pi^{R}_{i}(n) < C\, n^{-\delta}\right) \\
&\leq \Pro_{A} \left(\sum_{m=1}^{n} (R_{i}(m)-c^{R}_{i}(m,\mathfrak{D}^{R}_{m-1})) < -(C_p-C) \, n^{1-\delta}\right).
\end{aligned}
\end{equation}
In view of \eqref{c}, it holds
\begin{equation}
\Exp_{A}\left[ R_{i}(n)|\cF^{R}_{n-1}\right] = c^{R}_{i}(n,\mathfrak{D}^{R}_{n-1}), \quad \forall \, n \in \bN,
\end{equation}
which shows that $\left\{R_{i}(n)-c^{R}_{i}(n,\mathfrak{D}^{R}_{n-1}) \, :\, n \in \bN \right\}$ is a martingale difference, whose absolute value is also uniformly bounded by $1$. Thus, by application of the Azuma-Hoeffding inequality, we deduce that
\begin{equation*}
\Pro_{A}\left(\pi^{R}_{i}(n) < C\, n^{-\delta}\right) \leq e^{-\zeta n^{1-2\delta}},
\end{equation*}
where $\zeta >0$ is a constant. Therefore, we deduce \eqref{ts_st} by noting that there is $M>0$ such that for all $n \geq M$,
\begin{equation*}
n^{1-2\delta} \geq \frac{3}{\zeta} \ln(n).
\end{equation*}
\end{IEEEproof}
\end{document}